\documentclass[a4paper, 11pt]{article}
\usepackage{amssymb, amsmath,latexsym,amsfonts,amsbsy, amsthm,mathtools,graphicx,CJKutf8,CJKnumb,CJKulem,color,geometry}
\usepackage{verbatim}
\usepackage{hyperref}
\usepackage{newtxtext, newtxmath}
\usepackage[misc]{ifsym}
\usepackage[marginal]{footmisc}
\usepackage[titletoc,title]{appendix}
\def\l{\left}
\def\r{\right}
\def\f{\frac}

\def\az{\alpha}
\def\lz{\lambda}
\def\Lz{\Lambda}
\def\az{\alpha}
\def\rz{\rho}

\def\ez{\epsilon}
\def\bz{\beta}
\def\dz{\delta}
\def\gz{\gamma}
\def\tz{\theta}
\def\sz{\sigma}

\def\beq{\begin{equation}}
\def\eeq{\end{equation}}
\def\be{\begin{equation*}}
\def\ee{\end{equation*}}
\def\beqn{\begin{eqnarray}}
\def\eeqn{\end{eqnarray}}
\def\ben{\begin{eqnarray*}}
\def\een{\end{eqnarray*}}

\theoremstyle{plain}
\theoremstyle{plain}\newtheorem{thm}{Theorem}[section]
\theoremstyle{plain}\newtheorem{prop}{Proposition}[section]
\theoremstyle{plain}
\theoremstyle{plain}\newtheorem{lem}{Lemma}[section]
\theoremstyle{plain}

\numberwithin{equation}{section}

\makeatletter
\def\thanks#1{\protected@xdef\@thanks{\@thanks
        \protect\footnotetext{#1}}}
\makeatother

\begin{document}

\title{Stability of a coupled wave-Klein-Gordon system with non-compactly supported initial data}
\author{Qian Zhang}

\date{}

\maketitle

\noindent {\bf{Abstract}}\ \ In this paper we study global nonlinear stability for a system of semilinear wave and Klein-Gordon equations with quadratic nonlinearities. We consider nonlinearities of the type of wave-Klein-Gordon interactions where there are no derivatives on the wave component. The initial data are assumed to have a suitable polynomial decay at infinity, but are not necessarily compactly supported. We prove small data global existence, sharp pointwise decay estimates and linear scattering results for the solutions. We give a sufficient condition which combines the Sobolev norms of the nonlinearities on flat time slices outside the light cone and the norms of them on truncated hyperboloids inside the light cone, for the linear scattering of the solutions.

\bigskip


\bigskip


\section{Introduction}\label{s1}

We consider the following coupled wave and Klein-Gordon system in $\mathbb{R}^{1+3}$
\beq\label{s3: equv}
\begin{split}
-\Box u&=F_u:=P_{0}uv+P_{1}^\az\partial_\az uv+P_{0}^\az u\partial_\az v+P_1^{\az\bz}\partial_\az u\partial_\bz v\\
-\Box v+v&=F_v:=Q_0uv+Q_1^\az\partial_\az uv
\end{split}
\eeq
with prescribed initial data on the time slice $t=2$
\beq\label{s1: ini}
(u,\partial_tu,v,\partial_tv)|_{t=2}=(u_0,u_1,v_0,v_1).
\eeq
Here $P_0$, $P_1^\az$, $P_0^\az$, $P_1^{\az\bz}$, $Q_0$, $Q_1^\az$ are constants, $\az,\bz\in\{0,1,2,3\}$, $\partial_0=\partial_t$, $\partial_a=\partial_{x_a}$ for $a\in\{1,2,3\}$. No null conditions on the nonlinearities are assumed. Throughout, Einstein summation convention is adopted. As usual, $\Box=g^{\az\bz}\partial_\az\partial_\bz=-\partial_t^2+\Delta$ denotes the wave operator, where $g=\mathrm{diag}(-1,1,1,1)$ denotes the Minkowski metric in $\mathbb{R}^{1+3}$. The initial data are assumed to have a suitable polynomial decay at infinity, but are not necessarily compactly supported. We consider nonlinearities of the type of quadratic wave-Klein-Gordon interactions, with bad terms where there are no derivatives on the wave component in both equations. We take into account all of the possible semilinear nonlinearities of wave-Klein-Gordon interactions in the wave equation. We prove global nonlinear stability for \eqref{s3: equv}-\eqref{s1: ini}, including small data global existence, sharp pointwise decay estimates and linear scattering results.

Let us remind that John \cite{J81} proved that solutions to certain wave equations in $\mathbb{R}^{1+3}$ with quadratic nonlinearities may blow up in finite time. On the other hand, Klainerman \cite{K86} introduced the famous null condition to prove global existence of solutions for a broad class of wave equations with nonlinearities quadratic in derivatives (see also Christodoulou \cite{Ch}). Later global well-posedness for coupled wave and Klein-Gordon type systems has been widely studied,  with many important physical models behind, like the Dirac-Klein-Gordon model, the Klein-Gordon-Zakharov model and the Dirac-Proca model; see Bachelot \cite{Bac}, Georgiev \cite{Ge}, Ozawa, Tsutaya and Tsutsumi  \cite{OT, Ts}. For certain types of bad nonlinearities, such systems have also been studied by Psarelli \cite{Ps}, Katayama \cite{Ka, Ka2}, and LeFloch and Ma \cite{LM18}, Wang \cite{W1} and Ionescu and Pausader \cite{IP} as a model for the full Einstein-Klein-Gordon equation. For compactly supported initial data, Dong and Wyatt \cite{DW} (see also \cite{D2}) proved global stability results for \eqref{s3: equv}-\eqref{s1: ini} using the method of hyperboloidal foliation of space-time. The hyperboloidal method used in \cite{DW, D2} was introduced by Klainerman in \cite{K2} to prove global existence results for nonlinear Klein-Gordon equations by using commuting vector fields. Later this method was developed by Klainerman, Wang and Yang \cite{KWY,W}, LeFloch and Ma \cite{LM14}, and by Dong \cite{Dk} to study various coupled wave and Klein-Gordon type systems, mostly with the restriction that the initial data are compactly supported.

Motivated by the existing work on coupled wave and Klein-Gordon equations, we are interested in establishing global existence results with sharp time decay estimates for \eqref{s3: equv}-\eqref{s1: ini}, for non-compactly supported initial data which are only assumed to decay at suitable rates at infinity. Hence we remove the compactness assumption on the support of the initial data in \cite{DW}. We also prove linear scattering results for the solutions, which are completely new compared with \cite{DW, D2}. We give a sufficient condition which combines the Sobolev norms of the nonlinearities on flat time slices in the exterior region (i.e., outside the light cone) and the norms of them on truncated interior hyperboloids (i.e., the portions of hyperboloids inside the light cone), for the linear scattering of the solutions. This paper is inspired by a recent work by Huneau and Stingo \cite{HS} where the authors studied global existence of solutions to a certain class of quasilinear systems of wave equations on a product space with null nonlinearities. We will show in a subsequent work that the idea used in this paper is also applicable to other coupled wave and Klein-Gordon type systems.

In the sequel, we use $C$ to denote a universal constant whose value may change from line to line. As usual, $A\lesssim B$ means that $A\le CB$ for some constant $C$. Given a vector or a scalar $w$, we use Japanese bracket to denote $\langle w\rangle:=(1+|w|^2)^{1/2}$. 

\bigskip

\noindent${\mathbf{Major\ difficulties\ and\ key\ ideas.}}$ We set the initial data at time $t=2$ to simplify the computations and this is not essential by the invariance of \eqref{s3: equv} under time translations. For the system \eqref{s3: equv}-\eqref{s1: ini} where the quadratic nonlinearities do not satisfy the null condition, the main difficulties in proving global existence of solutions include: $i)$ the nonlinearities are critical non-integrable, i.e.,
\be
\||F_u|+|F_v|\|_{L^2(\mathbb{R}^3)}\sim t^{-1},
\ee
and in some terms of the nonlinearities there are no derivatives on the wave component $u$ and the standard energy cannot control norms like $\|u\|_{L^2(\mathbb{R}^3)}$; $ii)$ there is no compactness assumption on the support of the initial data, hence we cannot use the hyperboloidal method as in \cite{DW, D2} directly, and the difficulty $i)$ above causes more trouble in using the vector field method on flat time slices in $\mathbb{R}^{1+3}$. For this, we adopt an idea by Huneau and Stingo \cite{HS}, i.e., we divide the space-time into exterior and interior regions (i.e., outside the light cone and inside the light cone respectively) and use weighted energy and Sobolev estimates on flat time slices in the exterior region, while employing hyperboloidal foliation of space-time in the interior region. With the bad nonlinearities in \eqref{s3: equv} our case is more complicated due to the growth of the highest order energy of the wave component resulting from the terms $u\partial v$ and $\partial u\partial v$. More importantly, to control $L^2$ norms of the wave component $u$, we need to conduct a nonlinear transform $U=u+F_u$ (which was due to \cite{Ts} and was also used in \cite{DW, D2}) to control the conformal energy of $u$ on truncated exterior hyperboloids via corresponding estimates for $U$. The conformal energy of higher order derivatives of $U$ on flat time slices in the exterior region has a growth because of the terms $u\partial v,\partial u\partial v$ appearing in the nonlinearities in \eqref{s3: equv}. Hence we need to treat the conformal energy of $U$ for both the higher order and the lower order cases.

The main ideas are described below. We decompose our space-time into two regions: the exterior region $\mathscr{D}^{\rm{ex}}=\{(t,x): t\ge 2, |x|\ge t-1\}$ and the interior region $\mathscr{D}^{\rm{in}}=\{(t,x): t\ge 2, |x|<t-1\}$. We denote by $\Sigma^{\rm{ex}}_t:=\{x: |x|\ge t-1\}$ the truncated time slice in the exterior region and $\mathscr{l}_{[2,t]}:=\{(\tau,x): |x|=\tau-1, 2\le \tau\le t\}$ the portion of the boundary of the light cone in the time strip $[2,t]$. We introduce two weight functions $t^{-\dz}$ and $\omega(z):=(2+z)^{1+\lz}$ ($\dz\ge 0, \lz>0$ are constants), and define the weighted energy on $\Sigma^{\rm{ex}}_t$ as 
$$E^{{\rm{ex}},\lz}_{m,\dz}(t,u):=\int_{\Sigma^{\rm{ex}}_t}t^{-\dz}\omega(r-t)\big(|\partial u|^2+m^2u^2\big){\rm{d}}x$$ 
(see \eqref{dexfs}), where $r=|x|$, $m=0$ (for the wave component) or $1$ (for the Klein-Gordon component). In deriving estimates of $E^{{\rm{ex}},\lz}_{m,\dz}(t,u)$, we can also bound a space-time $L^2_{tx}$ norm and a boundary integral, i.e., we obtain estimates of    
\beq\label{s1: exnorm}
\|u\|_{X^{{\rm{ex}},\lz}_{m,\dz,t}}:=[E^{{\rm{ex}},\lz}_{m,\dz}(t,u)]^{\f{1}{2}}+\l(\int_2^t\!\!\int_{\Sigma^{\rm{ex}}_\tau}\tau^{-\dz}\omega'(r-\tau)\big(|Gu|^2+m^2u^2\big){\rm{d}}x{\rm{d}}\tau\r)^{\f{1}{2}}+\l(\int_{\mathscr{l}_{[2,t]}}\!\!\tau^{-\dz}\big(|Gu|^2+m^2u^2\big){\rm{d}}\sz\r)^{\f{1}{2}}
\eeq
(see Proposition \ref{s2: exfs}), where $|Gu|:=\l(\sum_{a=1}^3|G_au|^2\r)^{1/2}$, and $G_a=\partial_a+(x_a/r)\partial_t$ for $a\in\{1,2,3\}$ denote the good derivatives. Let $\mathscr{H}_s:=\{(t,x): t^2=s^2+|x|^2\}$ denote a hyperboloid at hyperbolic time $s\ge 2$, and $\mathscr{H}^{\rm{ex}}_s:=\mathscr{H}_s\cap\{|x|\ge(s^2-1)/2\}$ and $\mathscr{H}^{\rm{in}}_s:=\mathscr{H}_s\cap\{|x|<(s^2-1)/2\}$ denote the portions of $\mathscr{H}_s$ in the exterior region and the interior region respectively. The weighted energy functionals $E^{{\rm{ex}},h}_{m,\dz}(s,u)$ and $E^{{\rm{in}},h}_{m,\dz}(s,u)$ (see \eqref{dexhs} and \eqref{dinhs}) on these truncated hyperboloids are controlled  in terms of the weighted energy $E^{{\rm{ex}}}_{m,\dz}(t,u)$ (see \eqref{dexfsdz}) on flat time slices in the exterior region; see Propositions \ref{s2: exhs} and \ref{s2: inhs}. To control $L^2$ norms of the wave component, we also introduce the weighted conformal energy functionals $E^{c,\rm{ex}}_\dz(t,u)$ (see \eqref{s2: E^{c,ex}}) on flat time slices $\Sigma^{\rm{ex}}_t$ in the exterior region, $E^{c,{\rm{ex}},h}_\dz(s,u)$ (see \eqref{s2: E^c,ex,h_dz}) on truncated exterior hyperboloids $\mathscr{H}^{\rm{ex}}_s$ and $E^{c,\rm{in}}_\dz(s,u)$ (see \eqref{s2: E^c,in_dz}) on truncated interior hyperboloids.

To prove global-in-time existence of the solution $(u,v)$ to \eqref{s3: equv}-\eqref{s1: ini} in the exterior region, we introduce a bootstrap assumption for 
\beq\label{s1: bsN}
\begin{split}
&\quad\sum_{|I|\le N}\big\{\|\Gamma^Iu\|_{X^{{\rm{ex}},\lz}_{0,\dz,t}}+\|\Gamma^Iv\|_{X^{{\rm{ex}},\lz}_{1,0,t}}\big\}+\sum_{|I|\le N-1}\|\Gamma^Iu\|_{X^{{\rm{ex}},\lz}_{0,0,t}}\\
&+\sum_{|I|\le N}\|t^{-\f{\dz}{2}}(2+r-t)^{\f{\lz-1}{2}}\Gamma^Iu\|_{L^2_x(\Sigma^{\rm{ex}}_t)}+\sum_{|I|\le N-1}\|(2+r-t)^{\f{\lz-1}{2}}\Gamma^Iu\|_{L^2_x(\Sigma^{\rm{ex}}_t)},
\end{split}
\eeq
where $N\ge 9$ is an integer, $\lz\ge 3, 0<\dz\ll 1$, $\Gamma\in\l\{\partial, L,\Omega\r\}$ (where $\partial=(\partial_\az)_{0\le\az\le 3},L=(L_a)_{1\le a\le 3}, \Omega=(\Omega_{ab})_{1\le a<b\le 3}$ denote the vector fields of translations, Lorentz boosts and rotations respectively) and $I$ denotes a multi-index. By \eqref{s1: exnorm}, this not only controls the weighted $L^2_x$ norms of derivatives of $u$ and $v$ (i.e., $\|(2+r-t)^{\f{1+\lz}{2}}\big(t^{-\f{\dz}{2}}|\partial\Gamma^Iu|+|\partial\Gamma^Iv|+|\Gamma^Iv|\big)\|_{L^2_x(\Sigma^{\rm{ex}}_t)}$ for $|I|\le N$), but also the space-time $L^2_{tx}$ norm 
\beq\label{s1: L^2_tx}
\big\|\ \|(2+r-\tau)^{\f{\lz}{2}}\Gamma^Iv\|_{L^2_x(\Sigma^{\rm{ex}}_\tau)}\ \big\|_{L^2_\tau([2,t])}, \quad|I|\le N. 
\eeq
To close the bootstrap estimate of $\sum_{|I|\le N}\|\Gamma^Iu\|_{X^{{\rm{ex}},\lz}_{0,\dz,t}}$, we need to bound  $\int_2^t\|\tau^{-\f{\dz}{2}}(2+r-\tau)^{\f{1+\lz}{2}}\Gamma^IF_u\|_{L^2_x(\Sigma^{\rm{ex}}_\tau)}{\rm{d}}\tau$ for $|I|\le N$. We write $\Gamma^IF_u$ roughly as
\be\label{s1: gzF_u}
(u+\partial u)(\Gamma^Iv+\partial\Gamma^Iv)+(\Gamma^Iu+\partial\Gamma^Iu)(v+\partial v).
\ee
On the truncated time slice $\Sigma^{\rm{ex}}_\tau$ ($\tau\in[2,t]$), the lower order energy and $L^2$-type estimates of $u$ and the energy estimate of $v$ in \eqref{s1: bsN}, together with weighted Sobolev inequality in the exterior region, give the pointwise decay of $u$ (with a weight $(2+r-\tau)^{\f{\lz}{2}}$), $\partial u$ and $v$ (with a weight $(2+r-\tau)^{\f{1+\lz}{2}}$) at the rate of $\tau^{-1}$. We also have an extra decay property of $v$ (with a weight smaller than $(2+r-\tau)^{\f{1+\lz}{2}}$). Hence we obtain that the weighted $L^1_tL^2_x$ norm $\int_2^t\|\tau^{-\f{\dz}{2}}(2+r-\tau)^{\f{1+\lz}{2}}\Gamma^IF_u\|_{L^2_x(\Sigma^{\rm{ex}}_\tau)}{\rm{d}}\tau$ is bounded. Since the expression of $F_v$ does not contain the terms $u\partial v$, $\partial u\partial v$, we use the weighted $L^2_{tx}$ estimate of derivatives of $v$ in \eqref{s1: L^2_tx}, together with the weighted pointwise estimates of  $u$ and $\partial u$ to close the highest order energy estimate of $v$ (i.e., $\sum_{|I|\le N}\|\Gamma^Iv\|_{X^{{\rm{ex}},\lz}_{1,0,t}}$ in \eqref{s1: bsN}).  The energy estimate of lower order derivatives of $u$ can be refined as in the highest order energy estimate of $v$, using the weighted space-time $L^2_{tx}$ estimate \eqref{s1: L^2_tx}. Next, the estimates of the weighted $L^2_x$ norms of $u$ in the second line in \eqref{s1: bsN} can be closed by the improvements of the energy estimates of $u$ above and weighted Hardy inequality in the exterior region. Once we obtain refined estimates of \eqref{s1: bsN}, global existence of the solution $(u,v)$ in the exterior region follows.

Using the refined global energy and $L^2$-type estimates of $u$, we obtain the estimate of conformal energy of $u$ on flat time slices in the exterior region. For example, we can bound 
\be
\sum_{|I|\le N-1}\|t^{-\f{\dz}{2}}(2+r-t)^{\f{\lz-1}{2}}L_0\Gamma^Iu\|_{L^2_x(\Sigma^{\rm{ex}}_t)}+\sum_{|I|\le N-2}\|(2+r-t)^{\f{\lz-1}{2}}L_0\Gamma^Iu\|_{L^2_x(\Sigma^{\rm{ex}}_t)}
\ee
($L_0$ is the scaling vector field) by weighted Hardy inequality in the exterior region. Also, the refined global energy estimates of $(u,v)$ on flat time slices in the exterior region imply the estimates on truncated exterior hyperboloids. To give the estimate of conformal energy of $u$ on truncated exterior hyperboloids, we conduct a nonlinear transform $U=u+F_u$. The function $U$ solves the equation $-\Box U=S+W$, where 
$$S\approx\partial_\az u\partial^\az v+\partial_\az\partial u\partial^\az v+\partial_\az u\partial^\az\partial v+\partial_\az\partial u\partial^\az\partial v$$
is a sum of null forms and $W$ is a linear combination of cubic terms. Because of the growth of the conformal energy of derivatives of $U$ of order $N-1$ on flat time slices in the exterior region, we need to control the weighted conformal energy of $\Gamma^IU$ for $|I|\le N-1$  (unweighted for lower order case), i.e., we give estimates of
\beq\label{s1: bsN-1}
\sum_{|I|\le N-1}[E^{c,\rm{ex}}_\dz(t,\Gamma^IU)]^{\f{1}{2}}+\sum_{|I|\le N-2}[E^{c,\rm{ex}}(t,\Gamma^IU)]^{\f{1}{2}}.
\eeq 
For this, we need to bound $\int_2^t\|\tau^{-\f{\dz}{2}}r\Gamma^I(S+W)\|_{L^2_x(\Sigma^{\rm{ex}}_\tau)}{\rm{d}}\tau$ for $|I|\le N-1$ and $\int_2^t\|r\Gamma^I(S+W)\|_{L^2_x(\Sigma^{\rm{ex}}_\tau)}{\rm{d}}\tau$ for $|I|\le N-2$. We use the following estimate of the null form $\partial_\az\phi\partial^\az\psi$:
\be
|\partial_\az\phi\partial^\az\psi|\lesssim t^{-1}\sum_a\big\{|L_a\phi|\cdot|\partial\psi|+|\partial\phi|\cdot |L_a\psi|+|t-r||\partial\phi|\cdot|\partial\psi|\big\},
\ee
which means that we can bound $|\tau^{-\f{\dz}{2}}r\Gamma^IS|$ by
\beq\label{s1: rgzS}
\tau^{-\f{\dz}{2}}r\tau^{-1}\langle\tau-r\rangle\Big\{\big(|\Gamma u|+|\Gamma\partial u|\big)\big(|\Gamma\Gamma^{I}v|+|\Gamma\Gamma^{I}\partial v|\big)+\big(|\Gamma\Gamma^Iu|+|\Gamma\Gamma^I\partial u|\big)\big(|\Gamma v|+|\Gamma\partial v|\big)\Big\}.
\eeq
We observe that $\langle\tau-r\rangle\sim(2+r-\tau)$ on $\Sigma^{\rm{ex}}_\tau$ ($\tau\in[2,t]$). Hence using the weighted energy and $L^2$-type estimates of $u$ and $v$ in \eqref{s1: bsN} obtained above, corresponding weighted pointwise decay estimates and the extra decay property of $v$ in the exterior region, we obtain that the weighted $L^1_tL^2_x$ norm $\int_2^t\|\tau^{-\f{\dz}{2}}r\Gamma^IS\|_{L^2_x(\Sigma^{\rm{ex}}_\tau)}{\rm{d}}\tau$ is bounded for $|I|\le N-1$. For $|I|\le N-2$, the weighted $L^2_{tx}$ estimate of derivatives of $v$ (i.e. \eqref{s1: L^2_tx}) implies that the unweighted norm $\int_2^t\|r\Gamma^IS\|_{L^2_x(\Sigma^{\rm{ex}}_\tau)}{\rm{d}}\tau$ is bounded. 

Once we have the conformal energy estimates of $U$ on flat time slices in the exterior region, we obtain corresponding conformal energy bounds of $U$ on truncated exterior hyperboloids. The estimates of conformal energy and $L^2$ norms of $u$ on truncated exterior hyperboloids follow from the estimates of $U$ (recall that $U-u=F_u$). The lower order cases of these estimates, combined with the bootstrap estimates in the interior region, give pointwise decay estimates for the solution on truncated interior hyperboloids, via Sobolev inequality on hyperboloids. We also recall that in deriving exterior energy and $L^2$-type estimates of $u$ and $v$ and conformal energy estimates of $U$, we also obtain the estimates of corresponding integrals on the boundary of the light cone. These boundary integral estimates, together with the pointwise decay estimates on interior hyperboloids, are used to refine the bootstrap assumptions in the interior region.

After obtaining the global existence of the solution to \eqref{s3: equv}-\eqref{s1: ini}, we also prove linear scattering results for the solution. These cannot be obtained from the pointwise estimates since the Sobolev norms of the nonlinearities on flat time slices in $\mathbb{R}^{1+3}$ are not integrable in time. We give a sufficient condition which combines the Sobolev norms of the nonlinearities on the flat time slice $\Sigma^{\rm{ex}}_t$ in the exterior region and the norms of them on truncated interior hyperboloids $\mathscr{H}^{\rm{in}}_s$, for the linear scattering of the solution.

Next we state our main results. We denote the initial data in \eqref{s1: ini} by $\vec{u}_0=(u_0,u_1)$, $\vec{v}_0=(v_0,v_1)$ and
\beq\label{s1: nau_0}
|\vec{\nabla}^k\vec{u}_0|:=|\nabla^{k+1}u_0|+|\nabla^ku_1|,\quad\quad|\vec{\nabla}^k\vec{v}_0|:=|\nabla^{k+1}v_0|+|\nabla^kv_1|,
\eeq
for $k\in\mathbb{N}$, with $\nabla=(\partial_1,\partial_2,\partial_3)$.

\begin{thm}\label{thm1}
Let $\lz\ge 2$, and $N\in\mathbb{N}$ with $N\ge 9$. Then there exists a constant $\ez_0>0$ such that for any $0<\ez<\ez_0$ and all initial data $(u_0,u_1,v_0,v_1)$ satisfying the smallness condition 
\beq\label{s1: inismall}
\l\|\langle |x|\rangle^{\lz}|u_0|+\langle |x|\rangle^{\lz_N}|v_0|\r\|_{L^2(\mathbb{R}^3)}+\sum_{k=0}^N\big\|\langle|x|\rangle^{\lz_k}|\vec{\nabla}^k\vec{u}_0|+\langle|x|\rangle^{\lz_N}|\vec{\nabla}^k\vec{v}_0|\big\|_{L^2(\mathbb{R}^3)}\le\ez
\eeq
with $\lz_k:=k+1+\lz$ for $0\le k\le N$ and $|\vec{\nabla}^k\vec{u}_0|$, $|\vec{\nabla}^k\vec{v}_0|$ defined as in \eqref{s1: nau_0}, the Cauchy problem \eqref{s3: equv}-\eqref{s1: ini} admits a global-in-time solution $(u,v)$, which enjoys the following sharp pointwise decay estimates
\ben
|u(t,x)|&\lesssim&\ez\langle t+|x|\rangle^{-1}\langle t-|x|\rangle^{-\f{1}{2}},\\
|v(t,x)|&\lesssim&\ez\langle t+|x|\rangle^{-\f{3}{2}}.
\een
\end{thm}

Let $l\in\mathbb{N}$. We denote ${\bf H}^l(\mathbb{R}^3):=H^{l+1}(\mathbb{R}^3)\times H^l(\mathbb{R}^3)$ and ${\bf{H}}_{l}(\mathbb{R}^3):=\big(\dot{H}^{l+1}(\mathbb{R}^3)\cap\dot{H}^1(\mathbb{R}^3)\big)\times H^l(\mathbb{R}^3)$ (where $H^k(\mathbb{R}^3), \dot{H}^{k}(\mathbb{R}^3), k\in\mathbb{N}$ denote the Sobolev spaces and homogeneous Sobolev spaces respectively).

\begin{thm}\label{thm2}
Under the assumptions of Theorem \ref{thm1}, let $(u,v)$ be the global solution obtained there with initial data satisfying \eqref{s1: inismall} and denote $\vec{u}=(u,\partial_tu)$ and $\vec{v}=(v,\partial_tv)$. Then $(\vec{u},\vec{v})$ scatters to a free solution in ${\bf{H}}_{N-1}(\mathbb{R}^3)\times{\bf{H}}^{N-1}(\mathbb{R}^3)$, i.e., there exist $\vec{u}^*_0=(u^*_0,u^*_1)\in{\bf{H}}_{N-1}(\mathbb{R}^3)$ and $\vec{v}^*_0=(v^*_0,v^*_1)\in{\bf{H}}^{N-1}(\mathbb{R}^3)$ such that
\ben
&&\lim_{t\to+\infty}\|(u,\partial_tu)-(u^*,\partial_tu^*)\|_{{\bf{H}}_{N-1}(\mathbb{R}^3)}=0,\\
&&\lim_{t\to+\infty}\|(v,\partial_tv)-(v^*,\partial_tv^*)\|_{{\bf{H}}^{N-1}(\mathbb{R}^3)}=0,
\een
where $(u^*,v^*)$ is the solution to the $3D$ linear homogeneous wave-Klein-Gordon system with the initial data $(u^*_0,u^*_1,v^*_0,v^*_1)$. If we assume further that $F_u=P_0uv$ in \eqref{s3: equv}, then the solution $(\vec{u},\vec{v})$ scatters to a free solution in the highest energy space ${\bf{H}}_{N}(\mathbb{R}^3)\times{\bf{H}}^{N}(\mathbb{R}^3)$.
\end{thm}

The organization of this paper is as follows. In Section \ref{sP} we introduce some notations, give energy and $L^2$-type estimates and weighted Sobolev and Hardy inequalities both on flat time slices in the exterior region and on truncated hyperboloids. In Sections \ref{s3} - \ref{s4} we prove global existence results for \eqref{s3: equv}-\eqref{s1: ini} with pointwise decay estimates in the exterior region and the interior region respectively. In Section \ref{s5} we prove linear scattering results for the global solution obtained above.

\section{Preliminaries}\label{sP}

\subsection{Notations}\label{sno}
We work in the $(1+3)$ dimensional space-time $\mathbb{R}^{1+3}$ with Minkowski metric $g=(-1,1,1,1)$, which is used to raise or lower indices. The space indices are denoted by Roman letters $a,b,\cdots\in\{1,2,3\}$, and the space-time indices are denoted by Greek letters $\az,\bz,\gz,\cdots\in\{0,1,2,3\}$. Einstein summation convention for repeated upper and lower indices is adopted throughout the paper. We denote a point in $\mathbb{R}^{1+3}$ by $(t,x)=(x_0,x_1,x_2,x_3)$ with $t=x_0,x=(x_1,x_2,x_3),x^a=x_a,a=1,2,3$, and its spatial radius is denoted by $r:=|x|=\sqrt{x_1^2+x_2^2+x_3^2}$. The following vector fields will be used frequently in the analysis:
\begin{itemize}
\item[(i)] Translations: $\partial_\az:=\partial_{x_\az}$, for $\az\in\{0,1,2,3\}$.
\item[(ii)] Lorentz boosts: $L_a:=x_a\partial_t+t\partial_a$, for $a\in\{1,2,3\}$.
\item[(iii)] Rotations: $\Omega_{ab}:=x_a\partial_b-x_b\partial_a$, for $1\le a<b\le 3$.
\item[(iv)] Scaling: $L_0=t\partial_t+x^a\partial_a$.
\end{itemize}

We decompose the space-time $\mathbb{R}^{1+3}$ into two regions:
$$\mathrm{interior\ region\ }\quad \mathscr{D}^{\rm{in}}:=\{(t,x): t\ge 2, r<t-1\},$$
$$\mathrm{exterior\ region\ }\quad \mathscr{D}^{\rm{ex}}:=\{(t,x): t\ge 2, r\ge t-1\}.$$
We denote the constant time slices which foliate $\mathscr{D}^{\rm{ex}}$ as
$$\Sigma^{\rm{ex}}_t:=\{x\in\mathbb{R}^3: |x|\ge t-1\}.$$
The portion of exterior region in the time strip $[2,T]$ for any fixed time $T$ is denoted by
$$\mathscr{D}^{\rm{ex}}_T:=\{(t,x)\in\mathscr{D}^{\rm{ex}}: 2\le t\le T\}.$$
For $2\le T_1<T_2$, we denote the portion of the boundary of the light cone $\{r=t-1\}$ in the time interval $[T_1,T_2]$ by 
$$\mathscr{l}_{[T_1,T_2]}:=\{(t,x): r=t-1, T_1\le t\le T_2\}.$$
For any $s\ge 2$, we use
$$\mathscr{H}_s:=\{(t,x): t^2=s^2+|x|^2\}$$
to denote the hyperboloid at hyperbolic time $s$. We denote by $\mathscr{H}^{\rm{in}}_s$ (resp. $\mathscr{H}^{\rm{ex}}_s$) the portion of $\mathscr{H}_s$ contained in the interior region $\mathscr{D}^{\rm{in}}$ (resp. in the exterior region $\mathscr{D}^{\rm{ex}}$), i.e.,
$$\mathscr{H}^{\rm{in}}_s:=\{(t,x)\in\mathscr{H}_s: |x|<(s^2-1)/2\},$$
$$\mathscr{H}^{\rm{ex}}_s:=\{(t,x)\in\mathscr{H}_s: |x|\ge(s^2-1)/2\}.$$
We set 
\beq\label{s2: tr(s)}
t(s):=\f{s^2+1}{2},\quad\quad r(s):=\f{s^2-1}{2}. 
\eeq
In addition, we denote by $\mathscr{H}^{\rm{in}}_{[2,s]}$ the hyperbolic interior region limited by two hyperboloids $\mathscr{H}_2$ and $\mathscr{H}_s$, and by $\mathscr{H}^{\rm{ex}}_{[2,s]}$ the portion of the exterior region below $\mathscr{H}^{\rm{ex}}_{s}$, i.e.,
$$\mathscr{H}^{\rm{in}}_{[2,s]}:=\{(t,x)\in\mathscr{D}^{\rm{in}}: 2^2\le t^2-|x|^2\le s^2\},$$
$$\mathscr{H}^{\rm{ex}}_{[2,s]}:=\{(t,x)\in\mathscr{D}^{\rm{ex}}: t^2-|x|^2\le s^2\}.$$

We denote 
\beq\label{s2: defNa}
|\partial u|=\big(|\partial_tu|^2+|\nabla u|^2\big)^{\f{1}{2}},\quad\quad|\nabla u|:=\l(\sum_{a=1}^3|\partial_au|^2\r)^{\f{1}{2}}.
\eeq
\beq\label{s2: defLO} 
|Lu|:=\l(\sum_{a=1}^3|L_au|^2\r)^{\f{1}{2}},\quad\quad|\Omega u|:=\l(\sum_{1\le a<b\le 3}|\Omega_{ab}u|^2\r)^{\f{1}{2}}.
\eeq
We also denote 
\beq\label{s2: G_a}
G_a=\partial_a+\f{x_a}{r}\partial_t, \quad\quad |Gu|=\l(\sum_{a=1}^3|G_au|^2\r)^{\f{1}{2}}.
\eeq
and
\beq\label{s2: barpar_a}
\bar{\partial}_a:=t^{-1}L_a=\partial_a+\f{x_a}{t}\partial_t,\ \ \quad\quad\bar{\partial}_r=\f{x^a}{r}\bar{\partial}_a,\ \ \quad\quad|\bar{\nabla}u|:=\l(\sum_{a=1}^3|\bar{\partial}_au|^2\r)^{\f{1}{2}}.
\eeq
For any sufficiently smooth function $g=g(t,x)$ defined on a hyperboloid $\mathscr{H}_s$, we denote 
$$\int_{\mathscr{H}_s}g{\rm{d}}x:=\int_{\mathbb{R}^3}g(\sqrt{s^2+|x|^2},x){\rm{d}}x,$$
$$\int_{\mathscr{H}^{\rm{in}}_s}g{\rm{d}}x:=\int_{|x|<r(s)}g(\sqrt{s^2+|x|^2},x){\rm{d}}x,\quad\quad \int_{\mathscr{H}^{\rm{ex}}_s}g{\rm{d}}x:=\int_{|x|\ge r(s)}g(\sqrt{s^2+|x|^2},x){\rm{d}}x,$$
where $r(s)$ is as in \eqref{s2: tr(s)}. For any $1\le p<\infty$, we denote
\beq\label{s2: L^pH_s}
\|g\|_{L^p(\mathscr{H}_s)}:=\l(\int_{\mathscr{H}_s}|g|^p{\rm{d}}x\r)^{\f{1}{p}}=\l(\int_{\mathbb{R}^3}|g(\sqrt{s^2+|x|^2},x)|^p{\rm{d}}x\r)^{\f{1}{p}},
\eeq
$$\|g\|_{L^p(\mathscr{H}^{\rm{in}}_s)}:=\l(\int_{\mathscr{H}^{\rm{in}}_s}|g|^p{\rm{d}}x\r)^{\f{1}{p}},\quad\quad \|g\|_{L^p(\mathscr{H}^{\rm{ex}}_s)}:=\l(\int_{\mathscr{H}^{\rm{ex}}_s}|g|^p{\rm{d}}x\r)^{\f{1}{p}}.$$

We denote the ordered set
\beq\label{s2: gamma_k}
\{\Gamma_k\}_{k=1}^{10}:=\{(\partial_\az)_{0\le\az\le 3}, (L_a)_{1\le a\le 3}, (\Omega_{ab})_{1\le a<b\le 3}\}.
\eeq
For any multi-index $I=(i_1,\cdots,i_{10})\in\mathbb{N}^{10}$ of length $|I|=\sum_{k=1}^{10}i_k$, we denote
\beq\label{s2: Ga^I}
\Gamma^I=\prod_{k=1}^{10}\Gamma_k^{i_k},\quad\mathrm{where}\quad\Gamma=(\Gamma_1,\dots,\Gamma_{10}).
\eeq
For any $I=(i_0,i_1,i_2,i_3)\in\mathbb{N}^{4}$, $J=(j_1,j_2,j_3)\in\mathbb{N}^3$, $K=(k_1,k_2,k_3)\in\mathbb{N}^3$, let 
\beq\label{s2: pOL^J}
\partial^I:=\prod_{\az=0}^{3}\partial_\az^{i_\az},\quad\quad L^J:=\prod_{a=1}^{3}L_a^{j_a},\quad\quad\Omega^K:=\Omega_{12}^{k_1}\ \Omega_{13}^{k_2}\ \Omega_{23}^{k_3}.
\eeq

Given an index set $\Theta$, a finite set $\{P_{\!\tz}: \tz\in\Theta\}$ of linear operators and a linear operator $Q$, we write
\be
Q=\underset{\tz\in\Theta}{\sum\nolimits^{\prime}}P_{\!\tz}
\ee
if there exist some constants $c_{\!\tz}$ such that we have
\be
Q=\sum_{\tz\in\Theta}c_{\!\tz} P_{\!\tz}.
\ee
The following relations on commutators between the vector fields in $\{\Gamma_k\}_{k=1}^{10}$ are well-known (see for example \cite{Sog}):
\beqn\label{s2: commutator}
&&[\partial_\az, L_a]=\underset{0\le\bz\le 3}{\sum\nolimits^{\prime}}\partial_\bz,\quad\quad [\partial_\az, \Omega_{ab}]=\underset{0\le\bz\le 3}{\sum\nolimits^{\prime}}\partial_\bz,\quad\quad [\partial_\az, L_0]=\partial_\az,\quad\quad 0\le\az\le 3, 1\le a<b\le 3,\nonumber\\
&&[Z_k,Z_l]=\underset{Z_j\in V}{\sum\nolimits^{\prime}}Z_j,\quad\quad Z_k, Z_l\in V:=\{(L_a)_{1\le a\le 3}, (\Omega_{ab})_{1\le a<b\le 3}\},\nonumber\\
&&[L_0,\Omega_{ab}]=[L_0,L_c]=0,\quad\quad 1\le a<b\le 3, 1\le c\le 3,
\eeqn
where the commutator $[A,B]$ is defined as 
$$[A,B]:=AB-BA.$$

Below we state weighted energy and Sobolev type inequalities both on the flat time slices in the exterior region and on truncated hyperboloids. These were proved in \cite{HS} for a system of wave equations on a product space with quadratic nonlinearities which are linear combinations of null forms. Since in our case there are bad terms in the nonlinearities leading to the growth of energy of the solutions, we include more weight functions in the definitions of the energy functionals, and also prove $L^2$-type and pointwise estimates for the wave component.

\subsection{Energy estimates for wave and Klein-Gordon equations}

For any fixed constants $\lz>0$ and $\dz\ge 0$, we define
\beq\label{dexfs}
E^{{\rm{ex}},\lz}_{m,\dz}(t,u):=\int_{\Sigma^{\rm{ex}}_t}t^{-\dz}(2+r-t)^{1+\lz}(|\partial u|^2+m^2u^2)(t,x){\rm{d}}x.
\eeq
When $\dz=0$, we denote for simplicity
\beq\label{dexfsaz} 
E^{{\rm{ex}},\lz}_m(t,u):=\int_{\Sigma^{\rm{ex}}_t}(2+r-t)^{1+\lz}(|\partial u|^2+m^2u^2)(t,x){\rm{d}}x.
\eeq
We also denote 
\beq\label{dexfsdz} 
E^{{\rm{ex}}}_{m,\dz}(t,u):=\int_{\Sigma^{\rm{ex}}_t}t^{-\dz}(|\partial u|^2+m^2u^2)(t,x){\rm{d}}x,
\eeq
In particular, we define the unweighted energy
\beq\label{dexfs0} 
E^{{\rm{ex}}}_m(t,u):=\int_{\Sigma^{\rm{ex}}_t}(|\partial u|^2+m^2u^2)(t,x){\rm{d}}x.
\eeq

We have the following weighted energy estimate in the exterior region.

\begin{prop}\label{s2: exfs}
Let $\lz>0$ and $\dz\ge 0$. For any $t\ge 2$ and any sufficiently smooth function $u$ which is defined in the region $\mathscr{D}^{\rm{ex}}_t$ and decays sufficiently fast at space infinity, we have
\ben
\|u\|_{X^{{\rm{ex}},\lz}_{m,\dz,t}}\lesssim\|u\|_{X^{{\rm{ex}},\lz}_{m,\dz,2}}+\int_2^t\|\tau^{-\f{\dz}{2}}(2+r-\tau)^{\f{1+\lz}{2}}(-\Box u+m^2u)\|_{L^2_x(\Sigma^{\rm{ex}}_\tau)}{\rm{d}}\tau,
\een
where (we recall the definition \eqref{s2: G_a})
\beqn\label{s2: exnorm}
\|u\|_{X^{{\rm{ex}},\lz}_{m,\dz,t}}:=[E^{{\rm{ex}},\lz}_{m,\dz}(t,u)]^{\f{1}{2}}&+&\l(\int_2^t\int_{\Sigma^{\rm{ex}}_\tau}\tau^{-\dz}(2+r-\tau)^\lz\big(|Gu|^2+m^2u^2\big){\rm{d}}x{\rm{d}}\tau\r)^{\f{1}{2}}\nonumber\\
&+&\l(\int_{\mathscr{l}_{[2,t]}}\tau^{-\dz}\big(|Gu|^2+m^2u^2\big){\rm{d}}\sz\r)^{\f{1}{2}}
\eeqn

\begin{proof}
Let $\omega$ be a positive weight. Denote $F:=-\Box u+m^2u$ and multiply on both sides of this equality by $t^{-\dz}\omega(r-t)\partial_t u$ and by direct calculation, and we obtain 
\beqn\label{s2: omd_tubo}
\partial_t\big(t^{-\dz}\omega(|\partial u|^2+m^2u^2)\big)-2\partial_a\big(t^{-\dz}\omega\partial_tu\partial^au\big)&+&t^{-\dz}\omega'\big(|Gu|^2+m^2u^2\big)+\dz t^{-\dz-1}\omega(|\partial u|^2+m^2u^2)\nonumber\\
&=&2t^{-\dz}\omega\partial_tu F.
\eeqn
We note that the upward unit normal to $\mathscr{l}_{[2,t]}$ is $2^{-\f{1}{2}}(1,-x/r)$. Taking $\omega(z)=(2+z)^{1+\lz}$, and integrating \eqref{s2: omd_tubo} over the region $\mathscr{D}^{\rm{ex}}_{t}$, we have
\beqn\label{s2: exfse}
E^{{\rm{ex}},\lz}_{m,\dz}(t,u)&+&(1+\lz)\int_2^t\int_{\Sigma^{\rm{ex}}_\tau}\tau^{-\dz}(2+r-\tau)^\lz\big(|Gu|^2+m^2u^2\big){\rm{d}}x{\rm{d}}\tau+\int_{\mathscr{l}_{[2,t]}}\tau^{-\dz}\big(|Gu|^2+m^2u^2\big)2^{-\f{1}{2}}{\rm{d}}\sz\nonumber\\
&+&\dz\int_2^t\int_{\Sigma^{\rm{ex}}_\tau} \tau^{-\dz-1}(2+r-\tau)^{1+\lz}(|\partial u|^2+m^2u^2){\rm{d}}x{\rm{d}}\tau\nonumber\\
&=&E^{{\rm{ex}},\lz}_{m,\dz}(2,u)+2\int_2^t\int_{\Sigma^{\rm{ex}}_\tau}\tau^{-\dz}(2+r-\tau)^{1+\lz}\partial_tuF {\rm{d}}x{\rm{d}}\tau.
\eeqn
Differentiating \eqref{s2: exfse} with respect to $t$ and using H\"older inequality, we derive
\ben
\partial_t\Bigg\{E^{{\rm{ex}},\lz}_{m,\dz}(t,u)&+&\int_2^t\int_{\Sigma^{\rm{ex}}_\tau}\tau^{-\dz}(2+r-\tau)^\lz\big(|Gu|^2+m^2u^2\big){\rm{d}}x{\rm{d}}\tau+\int_{\mathscr{l}_{[2,t]}}\tau^{-\dz}\big(|Gu|^2+m^2u^2\big)2^{-\f{1}{2}}{\rm{d}}\sz\Bigg\}\nonumber\\
&\le&2[E^{{\rm{ex}},\lz}_{m,\dz}(t,u)]^{\f{1}{2}}\cdot\|t^{-\f{\dz}{2}}(2+r-t)^{\f{1+\lz}{2}}F\|_{L^2_x(\Sigma^{\rm{ex}}_t)},
\een
which implies the conclusion of the proposition.
\end{proof}
\end{prop}

We define the energy identity $e^h_m[u]$ on hyperboloids as (we recall the definitions \eqref{s2: defNa} and \eqref{s2: defLO})
\beqn\label{e^h_m}
e^h_m[u]:&=&|\partial u|^2+m^2u^2+2\f{x^a}{t}\partial_tu\partial_au=\f{s^2}{t^2}(\partial_tu)^2+t^{-2}|Lu|^2+m^2u^2\nonumber\\
&=&\f{s^2}{t^2}|\nabla u|^2+t^{-2}|L_0u|^2+t^{-2}|\Omega u|^2+m^2u^2.
\eeqn
Let $\dz\ge 0$. The energy functionals on the truncated hyperboloids $\mathscr{H}^{\rm{in}}_s$ and $\mathscr{H}^{\rm{ex}}_s$ are defined as
\beq\label{dinhs}
E^{{\rm{in}},h}_{m,\dz}(s,u):=\int_{\mathscr{H}^{\rm{in}}_s}t^{-\dz}e^h_m[u]{\rm{d}}x,
\eeq
\beq\label{dexhs}
E^{{\rm{ex}},h}_{m,\dz}(s,u):=\int_{\mathscr{H}^{\rm{ex}}_s}t^{-\dz}e^h_m[u]{\rm{d}}x
\eeq
respectively. We denote for simplicity 
\beq\label{dhs0}
E^{{\rm{in}},h}_{m}(s,u):=E^{{\rm{in}},h}_{m,0}(s,u),\quad\quad E^{{\rm{ex}},h}_m(s,u):=E^{{\rm{ex}},h}_{m,0}(s,u).
\eeq

On the truncated interior hyperboloids $\mathscr{H}^{\rm{in}}_s$, we have the following energy estimate.

\begin{prop}\label{s2: inhs}
Let $\dz\ge 0$. For any $s\ge 2$ and any sufficiently smooth function $u$ defined in $\mathscr{H}^{\rm{in}}_{[2,s]}$, we have
\ben
E^{{\rm{in}},h}_{m,\dz}(s,u)&\lesssim&E^{{\rm{in}},h}_{m,\dz}(2,u)+\int_{\mathscr{l}_{[\f{5}{2},\f{s^2+1}{2}]}}t^{-\dz}\big(|Gu|^2+m^2u^2\big){\rm{d}}\sz\\
&+&\int_2^s[E^{{\rm{in}},h}_{m,\dz}(\tau,u)]^{\f{1}{2}}\cdot\|t^{-\f{\dz}{2}}(-\Box u+m^2u)\|_{L^2(\mathscr{H}^{\rm{in}}_\tau)}{\rm{d}}\tau.
\een

\begin{proof}
We note that on $\mathscr{H}_{s}$ we have ${\bf n}{\rm{d}}\sz=(1,-x/t){\rm{d}}x$, where ${\bf n}$ and ${\rm{d}}\sz$ denote the upward unit normal and the volume element of $\mathscr{H}_s$ respectively. Taking $\omega\equiv 1$ in \eqref{s2: omd_tubo}, integrating over the region $\mathscr{H}^{\rm{in}}_{[2,s]}$, and using the transformation $(t,x)\to (\tau,x)$, where $\tau=\sqrt{t^2-|x|^2}$, we obtain
\ben
E^{{\rm{in}},h}_{m,\dz}(s,u)+\dz\int_{\mathscr{H}^{\rm{in}}_{[2,s]}}t^{-\dz-1}(|\partial u|^2+m^2u^2){\rm{d}}x{\rm{d}}t&=&E^{{\rm{in}},h}_{m,\dz}(2,u)+\int_{\mathscr{l}_{[\f{5}{2},\f{s^2+1}{2}]}}t^{-\dz}\big(|Gu|^2+m^2u^2\big)2^{-\f{1}{2}}{\rm{d}}\sz\\
&+&2\int_2^s\int_{\mathscr{H}^{{\rm{in}}}_\tau}t^{-\dz}\partial_tu(-\Box u+m^2u)\f{\tau}{t}{\rm{d}}x{\rm{d}}\tau.
\een
Using H\"older inequality, the conclusion follows.
\end{proof}
\end{prop}

The energy estimate on truncated exterior hyperboloids $\mathscr{H}^{\rm{ex}}_s$ is stated as follows.

\begin{prop}\label{s2: exhs}
Let $\dz\ge 0$. For any $s\in[2,\infty)$ and any sufficiently smooth function $u$ which is defined in $\mathscr{H}^{\rm{ex}}_{[2,s]}$ and decays sufficiently fast at space infinity, we have
\ben
&&E^{{\rm{ex}},h}_{m,\dz}(s,u)+\int_{\mathscr{l}_{[2,\f{s^2+1}{2}]}}t^{-\dz}\big(|Gu|^2+m^2u^2\big){\rm{d}}\sz\\
&\lesssim&E^{\rm{ex}}_{m,\dz}(2,u)+\sup_{t\in[2,\infty)}[E^{\rm{ex}}_{m,\dz}(t,u)]^{\f{1}{2}}\cdot\int_2^\infty\|t^{-\f{\dz}{2}}(-\Box u+m^2u)\|_{L^2_x(\Sigma^s_t)}{\rm{d}}t,
\een
where we recall \eqref{dexfsdz} for the definition of $E^{\rm{ex}}_{m,\dz}(t,u)$, and
\beq\label{s2: Sigma^s_t}
\Sigma^s_t:=\l\{\begin{array}{lcl}
\big\{x\in\mathbb{R}^3: |x|\ge t-1\big\},&&t\le(s^2+1)/2;\\
\big\{x\in\mathbb{R}^3: |x|\ge\sqrt{t^2-s^2}\big\}, &&t>(s^2+1)/2.
\end{array}\r.
\eeq

\begin{proof}
Taking $\omega\equiv 1$ in \eqref{s2: omd_tubo} and integrating over the region $\mathscr{H}^{\rm{ex}}_{[2,s]}$, we obtain
\ben
E^{{\rm{ex}},h}_{m,\dz}(s,u)&+&\int_{\mathscr{l}_{[2,\f{s^2+1}{2}]}}t^{-\dz}\big(|Gu|^2+m^2u^2\big)2^{-\f{1}{2}}{\rm{d}}\sz+\dz\int_{\mathscr{H}^{\rm{ex}}_{[2,s]}}t^{-\dz-1}(|\partial u|^2+m^2u^2){\rm{d}}x{\rm{d}}t\\
&=&E^{\rm{ex}}_{m,\dz}(2,u)+2\int_2^\infty\int_{\Sigma^s_t}t^{-\dz}\partial_tu(-\Box u+m^2u){\rm{d}}x{\rm{d}}t.
\een
We note that $\Sigma^s_t\subset\Sigma^{\rm{ex}}_t$, hence the conclusion follows.
\end{proof}
\end{prop}

\subsection{Conformal energy estimates for wave equation}

In this subsection, we denote by $n$ the spatial dimension. Unless otherwise specified, the conclusions in this subsection hold for $n=2$ or $3$. We first give conformal energy estimates on flat time slices in the exterior region. Let 
\beq\label{s2: K_0}
K_0=(t^2+r^2)\partial_t+2tx^a\partial_a. 
\eeq
Given $\dz\ge 0$, we denote (here we recall the definition \eqref{s2: defLO})
\beqn\label{s2: E^{c,ex}}
E^{c,\rm{ex}}_\dz(t,u):\!\!\!&=&\!\!\!\int_{\Sigma^{\rm{ex}}_t}t^{-\dz}\l\{\f{|K_0u+(n-1)tu|^2}{t^2+r^2}+\f{1}{t^2+r^2}|(r^2-t^2)\partial_ru+(n-1)ru|^2+(t^2+r^2)r^{-2}|\Omega u|^2\r\}{\rm{d}}x\nonumber\\
\!\!\!&=&\!\!\!\int_{\Sigma^{\rm{ex}}_t}t^{-\dz}\l\{|L_0u+(n-1)u|^2+|Lu|^2+|\Omega u|^2\r\}{\rm{d}}x.
\eeqn
When $\dz=0$, we denote for simplicity 
\beq\label{s2: E^{c,ex}_0}
E^{c,\rm{ex}}(t,u)=E^{c,\rm{ex}}_\dz(t,u).
\eeq

\begin{prop}\label{s2: excon}
Let $\dz\ge 0$. For any $t\in[2,\infty)$ and any sufficiently smooth function $u$ which is defined in $\mathscr{D}^{\rm{ex}}_t$ and decays sufficiently fast at space infinity, we have
\ben
[E^{c,\rm{ex}}_\dz(t,u)]^{\f{1}{2}}&+&\l(\int_{\mathscr{l}_{[2,t]}}\tau^{-\dz}\l\{\l|(\tau+r)(\partial_tu+\partial_ru)+(n-1)u\r|^2+(\tau-r)^2r^{-2}|\Omega u|^2\r\}{\rm{d}}\sz\r)^{\f{1}{2}}\\
&\lesssim&[E^{c,\rm{ex}}_\dz(2,u)]^{\f{1}{2}}+\int_2^t\|\tau^{-\f{\dz}{2}}(\tau^2+r^2)^{\f{1}{2}}(-\Box u)\|_{L^2_x(\Sigma^{\rm{ex}}_\tau)}{\rm{d}}\tau.
\een
If the spatial dimension $n=3$, we also have 
\ben
\int_{\Sigma^{\rm{ex}}_t}u^2(t,x){\rm{d}}x+\int_{|x|=t-1}u^2(t,x){\rm{d}}\sz(x)\lesssim E^{c,{\rm{ex}}}(t,u).
\een

\begin{proof}
Let $K_0$ be as in \eqref{s2: K_0}. By direct calculation,
\beqn\label{s2: boxuK_0u}
&&t^{-\dz}\big(K_0u+(n-1)tu\big)(-\Box u)\nonumber\\
&=&\f{1}{2}\partial_t\l(t^{-\dz}\l\{(t^2+r^2)|\partial u|^2+4tr\partial_tu\partial_ru+2(n-1)u\big(t\partial_tu+r\partial_ru\big)+(n-1)^2u^2\r\}\r)\nonumber\\
&+&\f{1}{2}\dz t^{-\dz-1}\l\{(t^2+r^2)|\partial u|^2+4tr\partial_tu\partial_ru+2(n-1)u\big(t\partial_tu+r\partial_ru\big)+(n-1)^2u^2\r\}\nonumber\\
&+&\partial_a\l(t^{-\dz}\l\{tx^a\big(|\nabla u|^2-|\partial_tu|^2\big)-(t^2+r^2)\partial_tu\partial^au-2tr\partial_ru\partial^au-(n-1)tu\partial^au-\f{n-1}{2}\partial_t(x^au^2)\r\}\r).\nonumber\\
\eeqn
Integrating \eqref{s2: boxuK_0u} over the region $\mathscr{D}^{\rm{ex}}_t$, we obtain
\beqn\label{s2: E^c,ex}
\int_{\Sigma^{\rm{ex}}_t}t^{-\dz}e_{\rm{con}}[u](t,x){\rm{d}}x&-&\int_{\Sigma^{\rm{ex}}_2}2^{-\dz}e_{\rm{con}}[u](2,x){\rm{d}}x+\int_{\mathscr{l}_{[2,t]}}\tau^{-\dz}H(\tau,x)2^{-\f{1}{2}}{\rm{d}}\sz+A(t)\nonumber\\
&=&2\int_2^t\int_{\Sigma^{\rm{ex}}_\tau}\tau^{-\dz}\big(K_0u+(n-1)\tau u\big)(-\Box u){\rm{d}}x{\rm{d}}\tau,
\eeqn
where 
\be
A(t):=\dz \int_2^t\int_{\Sigma^{\rm{ex}}_\tau}\tau^{-\dz-1}\l\{(\tau^2+r^2)|\partial u|^2+4\tau r\partial_tu\partial_ru+2(n-1)u\big(\tau\partial_tu+r\partial_ru\big)+(n-1)^2u^2\r\}{\rm{d}}x{\rm{d}}\tau,
\ee
\beqn\label{s2: e_con}
\!\!\!\!\!\!\!e_{\rm{con}}[u]:\!\!&=&(t^2+r^2)|\partial u|^2+4tr\partial_tu\partial_ru+2(n-1)u\big(t\partial_tu+r\partial_ru\big)+(n-1)^2u^2\nonumber\\
&=&\f{|K_0u+(n-1)tu|^2}{t^2+r^2}+\f{1}{t^2+r^2}|(r^2-t^2)\partial_ru+(n-1)ru|^2+(t^2+r^2)r^{-2}|\Omega u|^2
\eeqn
and
\beqn\label{s2: H(t,x)}
H(\tau,x):\!\!&=&(\tau^2+r^2)|\partial u|^2+4\tau r\partial_tu\partial_ru+2(n-1)u\big(\tau\partial_tu+r\partial_ru\big)+(n-1)^2u^2\nonumber\\
&-&\!\!\!\!\f{2x_a}{r}\Big\{\tau x^a\big(|\nabla u|^2-|\partial_tu|^2\big)-(\tau ^2+r^2)\partial_tu\partial^au-2 \tau r\partial_ru\partial^au-(n-1)\tau u\partial^au-\f{n-1}{2}\partial_t(x^au^2)\Big\}\nonumber\\
&=&\l|(\tau+r)(\partial_tu+\partial_ru)+(n-1)u\r|^2+(\tau-r)^2r^{-2}|\Omega u|^2.
\eeqn
We can also write \eqref{s2: e_con} as
\beq\label{s2: econ}
e_{\rm{con}}[u]=|L_0u+(n-1)u|^2+|Lu|^2+|\Omega u|^2.
\eeq
Differentiating \eqref{s2: E^c,ex} with respect to $t$ and using H\"older inequality, we obtain
\ben
\partial_t\Bigg\{E^{c,\rm{ex}}_\dz(t,u)+\int_{\mathscr{l}_{[2,t]}}\tau^{-\dz}H(\tau,x)2^{-\f{1}{2}}{\rm{d}}\sz\Bigg\}\le 2 [E^{c,\rm{ex}}_\dz(t,u)]^{\f{1}{2}}\cdot\|t^{-\f{\dz}{2}}(t^2+r^2)^{\f{1}{2}}(-\Box u)\|_{L^2_x(\Sigma^{\rm{ex}}_t)}.
\een
The first conclusion follows.

We turn to the estimate of $\|u\|_{L^2(\Sigma^{\rm{ex}}_t)}$ for $n=3$. Let $w:=ru$. We have
\beq\label{s2: w=ru}
\f{1}{t^2+r^2}|(r^2-t^2)\partial_ru+2ru|^2=\f{1}{t^2+r^2}\f{(r^2-t^2)^2}{r^2}|\partial_r w|^2+\f{r^2+t^2}{r^4}w^2+2\f{r^2-t^2}{r^3}w\partial_rw.
\eeq
We observe that
\beqn\label{s2: wboun}
2\int_{\Sigma^{\rm{ex}}_t}\f{r^2-t^2}{r^3}w\partial_rw{\rm{d}}x=-\int_{|x|=t-1}\f{r^2-t^2}{r^3}w^2(t,x){\rm{d}}\sz(x)-\int_{\Sigma^{\rm{ex}}_t}\f{r^2+t^2}{r^4}w^2{\rm{d}}x,
\eeqn
hence using H\"older inequality, we derive
\ben
\int_{\Sigma^{\rm{ex}}_t}\f{r^2+t^2}{r^4}w^2{\rm{d}}x\le\int_{|x|=t-1}\f{t^2-r^2}{r^3}w^2(t,x){\rm{d}}\sz(x)+2\l(\int_{\Sigma^{\rm{ex}}_t}\f{r^2+t^2}{r^4}w^2{\rm{d}}x\r)^{\f{1}{2}}\l(\int_{\Sigma^{\rm{ex}}_t}\f{(r^2-t^2)^2}{r^2(r^2+t^2)}|\partial_rw|^2{\rm{d}}x\r)^{\f{1}{2}},
\een
which implies
\beq\label{s2: w^2}
\l(\int_{\Sigma^{\rm{ex}}_t}\f{r^2+t^2}{r^4}w^2{\rm{d}}x\r)^{\f{1}{2}}\le 2\l(\int_{\Sigma^{\rm{ex}}_t}\f{(r^2-t^2)^2}{r^2(r^2+t^2)}|\partial_rw|^2{\rm{d}}x\r)^{\f{1}{2}}+\l(\int_{|x|=t-1}\f{t^2-r^2}{r^3}w^2(t,x){\rm{d}}\sz(x)\r)^{\f{1}{2}}.
\eeq
We note that \eqref{s2: wboun} also gives
\ben
\int_{\Sigma^{\rm{ex}}_t}\l\{\f{r^2+t^2}{r^4}w^2(t,x)+2\f{r^2-t^2}{r^3}w\partial_rw\r\}{\rm{d}}x
=\int_{|x|=t-1}\f{t^2-r^2}{r^3}w^2(t,x){\rm{d}}\sz(x),
\een
which together with \eqref{s2: w=ru} and \eqref{s2: E^{c,ex}_0} yields
\ben
\int_{\Sigma^{\rm{ex}}_t}\f{1}{t^2+r^2}\f{(r^2-t^2)^2}{r^2}|\partial_r w|^2{\rm{d}}x+\int_{|x|=t-1}\f{t^2-r^2}{r^3}w^2(t,x){\rm{d}}\sz(x)\le E^{c,{\rm{ex}}}(t,u).
\een
Combining this and \eqref{s2: w^2}, we obtain
\ben
\int_{\Sigma^{\rm{ex}}_t}u^2\l(1+\f{t^2}{r^2}\r){\rm{d}}x+\int_{|x|=t-1}\f{2t-1}{t-1}u^2(t,x){\rm{d}}\sz(x)\lesssim E^{c,{\rm{ex}}}(t,u).
\een
\end{proof}
\end{prop}

We next state conformal energy estimates on truncated interior hyperboloids. Given $\dz\ge 0$, we define 
\beqn\label{s2: E^c,in_dz}
E^{c,\rm{in}}_\dz(s,u):&=&\int_{\mathscr{H}^{\rm{in}}_s}t^{-\dz}\l\{\f{1}{t^2}|K_0u+(n-1)tu|^2+s^2|\bar{\nabla}u|^2\r\}{\rm{d}}x\nonumber\\
&=&\int_{\mathscr{H}^{\rm{in}}_s}t^{-\dz}\l\{|L_0u+x^a\bar{\partial}_au+(n-1)u|^2+s^2|\bar{\nabla}u|^2\r\}{\rm{d}}x,
\eeqn
where $K_0$ is as in \eqref{s2: K_0} and we recall \eqref{s2: barpar_a}. When $\dz=0$, we write for simplicity 
\beq\label{s2: E^c,in_0}
E^{c,\rm{in}}(s,u)=E^{c,\rm{in}}_0(s,u).
\eeq

\begin{prop}\label{s2: incon}
Let $0\le\dz\le 1/2$. For any $s\in[2,\infty)$ and any sufficiently smooth function $u$ defined in $\mathscr{H}^{\rm{in}}_{[2,s]}$, we have
\ben
E^{c,\rm{in}}_\dz(s,u)\lesssim E^{c,\rm{in}}_\dz(2,u)&+&\int_{\mathscr{l}_{[\f{5}{2},\f{s^2+1}{2}]}}t^{-\dz}\l\{\l|(t+r)(\partial_tu+\partial_ru)+(n-1)u\r|^2+(t-r)^2r^{-2}|\Omega u|^2\r\}{\rm{d}}\sz\\
&+&\int_2^s[E^{c,\rm{in}}_\dz(\tau,u)]^{\f{1}{2}}\cdot\|t^{-\f{\dz}{2}}\tau(-\Box u)\|_{L^2(\mathscr{H}^{\rm{in}}_\tau)}{\rm{d}}\tau.
\een
Let $t(s)$ and $r(s)$ be as in \eqref{s2: tr(s)}. When the spatial dimension $n=3$, we also have
\be
\l\|t^{-\f{\dz}{2}}\f{s}{r}u\r\|_{L^2(\mathscr{H}^{\rm{in}}_s)}\lesssim \l(\int_{|x|=r(s)}(t(s))^{-\dz}\f{s^2}{r}u^2(t(s),x){\rm{d}}\sz\r)^{\f{1}{2}}+[E^{c,\rm{in}}_\dz(s,u)]^{\f{1}{2}}.
\ee

\begin{proof}
Integrating \eqref{s2: boxuK_0u} over $\mathscr{H}^{\rm{in}}_{[2,s]}$, we obtain
\ben
\int_{\mathscr{H}^{\rm{in}}_s}t^{-\dz}e_{\rm{con}}^h[u]{\rm{d}}x-\int_{\mathscr{H}^{\rm{in}}_2}t^{-\dz}e_{\rm{con}}^h[u]{\rm{d}}x&-&\int_{\mathscr{l}_{[\f{5}{2},\f{s^2+1}{2}]}}t^{-\dz}H(t,x)2^{-\f{1}{2}}{\rm{d}}\sz+B(s)\\
&=&2\int_{\mathscr{H}^{\rm{in}}_{[2,s]}}t^{-\dz}\big(K_0u+(n-1)tu\big)(-\Box u){\rm{d}}x{\rm{d}}t,
\een
where
\be
B(s):=\dz \int_{\mathscr{H}^{\rm{in}}_{[2,s]}}t^{-\dz-1}\l\{(t^2+r^2)|\partial u|^2+4tr\partial_tu\partial_ru+2(n-1)u\big(t\partial_tu+r\partial_ru\big)+(n-1)^2u^2\r\}{\rm{d}}x{\rm{d}}t,
\ee
\ben
e_{\rm{con}}^h[u]:&=&(t^2+r^2)|\partial u|^2+4tr\partial_tu\partial_ru+2(n-1)u\big(t\partial_tu+r\partial_ru\big)+(n-1)^2u^2\\
&&-\f{2x_a}{t}\l\{tx^a\big(|\nabla u|^2-|\partial_tu|^2\big)-(t^2+r^2)\partial_tu\partial^au-2tr\partial_ru\partial^au-(n-1)tu\partial^au-\f{n-1}{2}\partial_t(x^au^2)\r\}\\
&=&(t^2+3r^2)|\partial_tu|^2+(t^2-r^2)|\nabla u|^2+\f{2r}{t}\big(3t^2+r^2\big)\partial_tu\partial_ru+4r^2|\partial_ru|^2+\f{2(n-1)}{t}\big(t^2+r^2\big)u\partial_tu\\
&&+4(n-1)ru\partial_ru+(n-1)^2u^2
\een
and $H(t,x)$ is as in \eqref{s2: H(t,x)}. Using the definition \eqref{s2: barpar_a} we have $\partial_a=\bar{\partial}_a-(x_a/t)\partial_t$, $\partial_r=\bar{\partial}_r-(r/t)\partial_t$ and 
\be
|\nabla u|^2-|\partial_ru|^2=|\bar{\nabla}u|^2-|\bar{\partial}_ru|^2=r^{-2}|\Omega u|^2.
\ee
Substituting these identities into the expression of $e^h_{\rm{con}}[u]$ and using that $s^2=t^2-r^2$, we obtain
\beqn\label{s2: e_con^h}
e_{\rm{con}}^h[u]=\l|\f{s^2}{t}\partial_tu+2r\bar{\partial}_ru+(n-1)u\r|^2+s^2|\bar{\nabla}u|^2=\f{1}{t^2}\l|K_0u+(n-1)tu\r|^2+s^2|\bar{\nabla}u|^2.
\eeqn
We can also write $e_{\rm{con}}^h[u]$ as
\beq\label{s2: e_con^h'}
e_{\rm{con}}^h[u]=\l|L_0u+x^a\bar{\partial}_au+(n-1)u\r|^2+s^2|\bar{\nabla}u|^2.
\eeq
Using the transform $(t,x)\to(\tau,x)$, $\tau=\sqrt{t^2-|x|^2}$ and \eqref{s2: H(t,x)}, we obtain
\ben
\int_{\mathscr{H}^{\rm{in}}_s}t^{-\dz}e_{\rm{con}}^h[u]{\rm{d}}x&\le&\int_{\mathscr{H}^{\rm{in}}_2}t^{-\dz}e_{\rm{con}}^h[u]{\rm{d}}x\\
&+&\int_{\mathscr{l}_{[\f{5}{2},\f{s^2+1}{2}]}}t^{-\dz}\l\{\l|(t+r)(\partial_tu+\partial_ru)+(n-1)u\r|^2+(t-r)^2r^{-2}|\Omega u|^2\r\}2^{-\f{1}{2}}{\rm{d}}\sz\\
&+&2\int_2^s\int_{\mathscr{H}^{\rm{in}}_{\tau}}t^{-\dz}\tau(-\Box u)\f{K_0u+(n-1)tu}{t}{\rm{d}}x{\rm{d}}\tau.
\een
Using \eqref{s2: e_con^h} and H\"older inequality, we obtain the first conclusion. 

We turn to the proof of the second conclusion. For $n=3$, let $\tilde{u}_s(x):=u\big(\sqrt{s^2+|x|^2},x\big)$. We have
\be
\int_{|x|< r(s)}(s^2+r^2)^{-\f{\dz}{2}}\f{x^a\partial_a\tilde{u}^2}{r^2}{\rm{d}}x=\int_{|x|=r(s)}(t(s))^{-\dz}\f{\tilde{u}^2}{r}{\rm{d}}\sz-\int_{|x|< r(s)}\tilde{u}^2\partial_a\l((s^2+r^2)^{-\f{\dz}{2}}\f{x^a}{r^2}\r){\rm{d}}x.
\ee
We observe that
\ben
\partial_a\l((s^2+r^2)^{-\f{\dz}{2}}\f{x^a}{r^2}\r)=(s^2+r^2)^{-\f{\dz}{2}}\f{1}{r^2}\l(1-\dz\f{r^2}{s^2+r^2}\r)\ge(s^2+r^2)^{-\f{\dz}{2}}\f{1}{2r^2},
\een 
which implies
\be
\int_{|x|< r(s)}(s^2+r^2)^{-\f{\dz}{2}}\f{x^a\partial_a\tilde{u}^2}{r^2}{\rm{d}}x\le\int_{|x|=r(s)}(t(s))^{-\dz}\f{\tilde{u}^2}{r}{\rm{d}}\sz-\int_{|x|\le r(s)}(s^2+r^2)^{-\f{\dz}{2}}\f{\tilde{u}^2}{2r^2}{\rm{d}}x.
\ee
Then H\"older inequality gives
\ben
\int_{|x|< r(s)}(s^2+r^2)^{-\f{\dz}{2}}\f{\tilde{u}^2}{2r^2}{\rm{d}}x&\le&\int_{|x|=r(s)}(t(s))^{-\dz}\f{\tilde{u}^2}{r}{\rm{d}}\sz\\
&+&2\l(\int_{|x|< r(s)}(s^2+r^2)^{-\f{\dz}{2}}\f{\tilde{u}^2}{r^2}{\rm{d}}x\r)^{\f{1}{2}}\l(\int_{|x|< r(s)}(s^2+r^2)^{-\f{\dz}{2}}|\partial_r\tilde{u}|^2{\rm{d}}x\r)^{\f{1}{2}}.
\een
It follows that
\be
\int_{|x|< r(s)}(s^2+r^2)^{-\f{\dz}{2}}\f{\tilde{u}^2}{r^2}{\rm{d}}x\lesssim\int_{|x|=r(s)}(t(s))^{-\dz}\f{\tilde{u}^2}{r}{\rm{d}}\sz+\int_{|x|< r(s)}(s^2+r^2)^{-\f{\dz}{2}}|\partial_r\tilde{u}|^2{\rm{d}}x
\ee
or
\be
\int_{\mathscr{H}^{\rm{in}}_s}t^{-\dz}\f{u^2}{r^2}{\rm{d}}x\lesssim\int_{|x|=r(s)}(t(s))^{-\dz}\f{u^2(t(s),x)}{r}{\rm{d}}\sz+\int_{\mathscr{H}^{\rm{in}}_s}t^{-\dz}|\bar{\partial}_ru|^2{\rm{d}}x.
\ee
\end{proof}
\end{prop}

The conformal energy estimate on truncated exterior hyperboloids is given by the proposition below. 

For $\dz\ge 0$, we denote 
\beq\label{s2: E^c,ex,h_dz}
E^{c,{\rm{ex}},h}_\dz(s,u):=\int_{\mathscr{H}^{\rm{ex}}_s}t^{-\dz}\big(|L_0u+x^a\bar{\partial}_au+(n-1)u|^2+s^2|\bar{\nabla}u|^2\big){\rm{d}}x,
\eeq 
\beq\label{s2: E^c,ex,h_0}
E^{c,{\rm{ex}},h}(s,u)=E^{c,{\rm{ex}},h}_0(s,u).
\eeq

\begin{prop}\label{s2: exconh}
Let $0\le\dz\le 1/2$. For any $s\ge 2$ and any $t>t_1>s$, let $\mathscr{H}_s\cap[t_1,t]$ denote the portion of $\mathscr{H}_s$ in the time strip $[t_1,t]$. Denote by $\tilde{\mathscr{D}}^s_{[t_1,t]}$ the region bounded by $\mathscr{H}_s\cap[t_1,t]$, the time slices $t_1$ and $t$. For any sufficiently smooth function $u$ which is defined in $\tilde{\mathscr{D}}^s_{[t_1,t]}$ and decays sufficiently fast at space infinity, we have
\ben
&&\l(\int_{\mathscr{H}_s\cap[t_1,t]}\tau^{-\dz}\big(|L_0u+x^a\bar{\partial}_au+(n-1)u|^2+s^2|\bar{\nabla}u|^2\big){\rm{d}}x\r)^{\f{1}{2}}\\
&\lesssim&\l(\int_{\Gamma^s_{t_1}}t_1^{-\dz}\big(|L_0u+(n-1)u|^2+|Lu|^2+|\Omega u|^2\big){\rm{d}}x\r)^{\f{1}{2}}+\int_{t_1}^t\|\tau^{-\f{\dz}{2}}(\tau^2+r^2)^{\f{1}{2}}(-\Box u)\|_{L^2_x(\Gamma^s_\tau)}{\rm{d}}\tau,
\een
where $\Gamma^s_\tau:=\{x\in\mathbb{R}^n: |x|\ge (\tau^2-s^2)^{\f{1}{2}}\}$, for $\tau\in[t_1,t]$. In particular,
\be
[E^{c,{\rm{ex}},h}_\dz(s,u)]^{\f{1}{2}}\lesssim [E^{c,\rm{ex}}_\dz(t(s),u)]^{\f{1}{2}}+\int_{t(s)}^\infty\|\tau^{-\f{\dz}{2}}r(-\Box u)\|_{L^2_x(\Gamma^s_\tau)}{\rm{d}}\tau,
\ee
where we recall the definitions \eqref{s2: E^{c,ex}} and \eqref{s2: tr(s)}. When the spatial dimension $n=3$, we also have
\be
\l\|t^{-\f{\dz}{2}}\f{s}{r}u\r\|_{L^2(\mathscr{H}^{\rm{ex}}_s)}+\l(\int_{|x|=r(s)}(t(s))^{-\dz}\f{s^2}{r}u^2(t(s), x){\rm{d}}\sz\r)^{\f{1}{2}}\lesssim\|t^{-\f{\dz}{2}}s|\bar{\nabla}u|\|_{L^2(\mathscr{H}^{\rm{ex}}_s)},
\ee
where $t(s)=\f{s^2+1}{2}$, $r(s)=\f{s^2-1}{2}$ (see \eqref{s2: tr(s)}) and we recall \eqref{s2: barpar_a}.

\begin{proof}
Integrating \eqref{s2: boxuK_0u} over the region $\tilde{\mathscr{D}}^s_{[t_1,t]}$, we obtain
\beqn\label{s2: econGst}
\int_{\Gamma^s_{t}}t^{-\dz}e_{\rm{con}}[u]{\rm{d}}x&-&\int_{\Gamma^s_{t_1}}t_1^{-\dz}e_{\rm{con}}[u]{\rm{d}}x+\int_{\mathscr{H}_s\cap[t_1,t]}\tau^{-\dz}e_{\rm{con}}^h[u]{\rm{d}}x+\tilde{A}(t)\nonumber\\
&=&2\int_{t_1}^{t}\int_{\Gamma^s_\tau}\tau^{-\dz}\big(K_0u+(n-1)\tau u\big)(-\Box u){\rm{d}}x{\rm{d}}\tau,
\eeqn
where $e_{\rm{con}}[u]$ and $e_{\rm{con}}^h[u]$ are as in \eqref{s2: e_con} (i.e., \eqref{s2: econ}) and \eqref{s2: e_con^h} (i.e.,\eqref{s2: e_con^h'}) respectively, and
\be
\tilde{A}(t):=\dz \int_{t_1}^t\int_{\Gamma^s_\tau}\tau^{-\dz-1}\l\{(\tau^2+r^2)|\partial u|^2+4\tau r\partial_tu\partial_ru+2(n-1)u\big(\tau\partial_tu+r\partial_ru\big)+(n-1)^2u^2\r\}{\rm{d}}x{\rm{d}}\tau.
\ee
Differentiating \eqref{s2: econGst} with respect to $t$ and using H\"older inequality, we obtain 
\ben
\partial_t\Bigg\{\int_{\Gamma^s_{t}}t^{-\dz}e_{\rm{con}}[u]{\rm{d}}x+\int_{\mathscr{H}_s\cap[t_1,t]}\tau^{-\dz}e_{\rm{con}}^h[u]{\rm{d}}x\Bigg\}\le 2\|t^{-\f{\dz}{2}}(t^2+r^2)^{\f{1}{2}}(-\Box u)\|_{L^2_x(\Gamma^s_t)}\cdot\l(\int_{\Gamma^s_{t}}t^{-\dz}e_{\rm{con}}[u]{\rm{d}}x\r)^{\f{1}{2}}.
\een
This gives the first conclusion. 

Let $n=3$. We set $\tilde{u}_s(x):=u\big(\sqrt{s^2+|x|^2},x\big)$ and note that
\be
\int_{|x|\ge r(s)}(s^2+r^2)^{-\f{\dz}{2}}\f{x^a\partial_a\tilde{u}^2}{r^2}{\rm{d}}x=-\int_{|x|=r(s)}(t(s))^{-\dz}\f{\tilde{u}^2}{r}{\rm{d}}\sz-\int_{|x|\ge r(s)}\tilde{u}^2\partial_a\l((s^2+r^2)^{-\f{\dz}{2}}\f{x^a}{r^2}\r){\rm{d}}x.
\ee
Hence similar to the proof of Proposition \ref{s2: incon}, we obtain
\be
\int_{|x|\ge r(s)}(s^2+r^2)^{-\f{\dz}{2}}\f{\tilde{u}^2}{r^2}{\rm{d}}x+\int_{|x|=r(s)}(t(s))^{-\dz}\f{\tilde{u}^2}{r}{\rm{d}}\sz\lesssim\int_{|x|\ge r(s)}(s^2+r^2)^{-\f{\dz}{2}}|\partial_r\tilde{u}|^2{\rm{d}}x,
\ee
i.e.,
\be
\int_{\mathscr{H}^{\rm{ex}}_s}t^{-\dz}\f{u^2}{r^2}{\rm{d}}x+\int_{|x|=r(s)}(t(s))^{-\dz}\f{u^2(t(s), x)}{r}{\rm{d}}\sz\lesssim\int_{\mathscr{H}^{\rm{ex}}_s}t^{-\dz}|\bar{\partial}_ru|^2{\rm{d}}x.
\ee
\end{proof}
\end{prop}

\subsection{Sobolev and Hardy inequalities}

The following weighted Sobolev inequality in the exterior region was given in \cite[Lemmas 4.1 and 4.2]{HS}. For readers' convenience, we include the proof in Appendix \ref{sA}.

\begin{lem}\label{s2: Sobex}
Let $\Lz\in\mathbb{R}$. For any $t\in[2,\infty)$ and any sufficiently smooth function $u$ defined on $\Sigma^{\rm{ex}}_t$, we have
\be
\sup_{\Sigma^{\rm{ex}}_t}(2+r-t)^\Lz r^2|u(t,x)|^2\lesssim\sum_{|J|\le 2}\int_{\Sigma^{\rm{ex}}_t}(2+r-t)^\Lz(|\partial_r\Omega^Ju|^2+|\Omega^Ju|^2){\rm{d}}x,
\ee
\be
\sup_{\Sigma^{\rm{ex}}_t}(2+r-t)^\Lz r^2|u(t,x)|^2\lesssim\sum_{|J|\le 2}\int_{\Sigma^{\rm{ex}}_t}\l\{(2+r-t)^{\Lz+1}|\partial_r\Omega^Ju|^2+(2+r-t)^{\Lz-1}|\Omega^Ju|^2\r\}{\rm{d}}x,
\ee
where we recall the definition \eqref{s2: pOL^J}.
\end{lem}

\begin{lem}\label{s2: Hardyex}
Let $\Lz\ge 2$. For any $t\in [2,\infty)$ and any sufficiently smooth function $u$ which is defined on $\Sigma^{\rm{ex}}_t$ and satisfies 
$$\int_{\Sigma^{\rm{ex}}_t}(2+r-t)^\Lz u^2(t,x){\rm{d}}x<\infty,$$ 
we have
\be
\int_{\Sigma^{\rm{ex}}_t}(2+r-t)^\Lz u^2{\rm{d}}x\lesssim\int_{\Sigma^{\rm{ex}}_t}(2+r-t)^{\Lz+2}|\partial_ru|^2{\rm{d}}x.
\ee
In the case that $u$ is compactly supported in $x$, we only need the bound $\Lz>-1$.
\begin{proof}
We first assume that $u$ is compactly supported in $x$. This case was proved in \cite[Lemmas 4.5]{HS}. Using that
\be
\partial_r\big((2+r-t)^{\Lz+1}r^2\big)\ge(\Lz+1)(2+r-t)^\Lz r^2,
\ee
we obtain
\ben
\int_{\Sigma^{\rm{ex}}_t}(2+r-t)^\Lz u^2(t,x){\rm{d}}x&\le&\f{1}{\Lz+1}\int_{\Sigma^{\rm{ex}}_t}\f{x^a}{r^3}\partial_a\big((2+r-t)^{\Lz+1}r^2\big)u^2(t,x){\rm{d}}x\\
&=&-\f{1}{\Lz+1}\int_{|x|=t-1}u^2(t,x){\rm{d}}\sz(x)-\f{1}{\Lz+1}\int_{\Sigma^{\rm{ex}}_t}(2+r-t)^{\Lz+1}r^2\partial_a\l(\f{x^a}{r^3}u^2\r){\rm{d}}x\\
&\le&-\f{2}{\Lz+1}\int_{\Sigma^{\rm{ex}}_t}(2+r-t)^{\Lz+1}u\partial_ru{\rm{d}}x\\
&\le&\f{2}{\Lz+1}\l(\int_{\Sigma^{\rm{ex}}_t}(2+r-t)^{\Lz}u^2{\rm{d}}x\r)^{\f{1}{2}}\cdot\l(\int_{\Sigma^{\rm{ex}}_t}(2+r-t)^{\Lz+2}|\partial_ru|^2{\rm{d}}x\r)^{\f{1}{2}},
\een
which implies 
\ben
\l(\int_{\Sigma^{\rm{ex}}_t}(2+r-t)^\Lz u^2{\rm{d}}x\r)^{\f{1}{2}}
\le\f{2}{\Lz+1}\l(\int_{\Sigma^{\rm{ex}}_t}(2+r-t)^{\Lz+2}|\partial_ru|^2{\rm{d}}x\r)^{\f{1}{2}}.
\een
In the general case that $u$ is not compactly supported in $x$, we choose a cut-off function $\chi\in C^\infty_0(\mathbb{R})$ with $r|\chi'(r)|\le 5/4$. Applying the last estimate to $\chi(\ez r)u$ for any $\ez>0$, we obtain
\ben
\int_{\Sigma^{\rm{ex}}_t}(2+r-t)^\Lz |\chi(\ez r)u|^2{\rm{d}}x
&\le&\f{4}{(\Lz+1)^2}\int_{\Sigma^{\rm{ex}}_t}(2+r-t)^{\Lz+2}\big\{|\chi(\ez r)\partial_ru|^2+|\ez\chi'(\ez r)|^2u^2\big\}{\rm{d}}x\\
&\le&\f{4}{(\Lz+1)^2}\int_{\Sigma^{\rm{ex}}_t}\Big\{(2+r-t)^{\Lz+2}|\partial_ru|^2+\f{25}{16}(2+r-t)^{\Lz}u^2\Big\}{\rm{d}}x.
\een
Let $\ez\to 0$ and recall that $\Lz\ge 2$ (which gives $\f{25}{4}\f{1}{(\Lz+1)^2}\le\f{25}{36}<1$), and we obtain
\be
\int_{\Sigma^{\rm{ex}}_t}(2+r-t)^\Lz u^2{\rm{d}}x\lesssim\int_{\Sigma^{\rm{ex}}_t}(2+r-t)^{\Lz+2}|\partial_ru|^2{\rm{d}}x.
\ee
\end{proof}
\end{lem}

We denote $\{Z_k\}_{k=1}^6:=\{(L_a)_{1\le a\le 3},(\Omega_{ab})_{1\le a<b\le 3}\}$. For any multi-index $J=(j_1,\cdots,j_{6})\in\mathbb{N}^{6}$ of length $|J|=\sum_{k=1}^{6}j_k$, we denote 
\beq\label{s2: Z^J}
Z^J=\prod_{k=1}^{6}Z_k^{j_k}.
\eeq 

\begin{lem}\label{s2: exL_auL^2}
Let $\Lz>-1$ and $k\in\{1,\cdots,6\}$. For any $t\in[2,\infty)$ and any sufficiently smooth function $u$ defined on $\Sigma^{\rm{ex}}_t$, we have
\ben
&&\int_{\Sigma^{\rm{ex}}_t}(2+r-t)^\Lz|Z_ku|^2{\rm{d}}x\lesssim\sum_{|J|\le 1}\int_{\Sigma^{\rm{ex}}_t}(2+r-t)^{\Lz+2}|\partial Z^{J}u|^2{\rm{d}}x,\\
&&\int_{\Sigma^{\rm{ex}}_t}(2+r-t)^\Lz|L_0u|^2{\rm{d}}x\lesssim\sum_{|J|\le 1}\int_{\Sigma^{\rm{ex}}_t}(2+r-t)^{\Lz+2}|\partial Z^Ju|^2{\rm{d}}x.
\een
\begin{proof}
The second inequality follows from the first one, using the relation 
\beq\label{L_0}
L_0=(t-r)\partial_t+(r-t)\partial_r+(x^a/r)L_a
\eeq
and the observation that 
\beq\label{s2: Edecay0}
|t-r|\lesssim(2+r-t)\quad\mathrm{in}\ \mathscr{D}^{\rm{ex}}.
\eeq 
The proof of the first inequality was given in \cite[Corollary 4.6]{HS}. We choose a cut-off function $\chi\in C^\infty_0(\mathbb{R})$, apply Lemma \ref{s2: Hardyex} to $\chi(\ez r)Z_ku$ for any $\ez>0$, and obtain
\ben
\int_{\Sigma^{\rm{ex}}_t}(2+r-t)^\Lz|\chi(\ez r)Z_ku|^2{\rm{d}}x&\lesssim&\int_{\Sigma^{\rm{ex}}_t}(2+r-t)^{\Lz+2}\l\{|\chi(\ez r)|^2|\partial_rZ_ku|^2+\ez^2|\chi'(\ez r)|^2|Z_ku|^2\r\}{\rm{d}}x\\
&\lesssim&\int_{\Sigma^{\rm{ex}}_t}(2+r-t)^{\Lz+2}\l\{|\partial_rZ_ku|^2+|\partial u|^2\r\}{\rm{d}}x,
\een
where we use that $|Z_ku|\lesssim r|\partial u|$ on $\Sigma^{\rm{ex}}_t$ and that $|\ez r\chi'(\ez r)|^2\lesssim 1$. Let $\ez\to 0$, the conclusion follows.
\end{proof}
\end{lem}

\begin{lem}\label{s2: exL_auL^infty}
Let $\lz>0$ and $k\in\{1,\cdots,6\}$. For any $t\in[2,\infty)$ and any sufficiently smooth function $u$ defined on $\Sigma^{\rm{ex}}_t$, we have
\ben
&&\sup_{\Sigma^{\rm{ex}}_t}(2+r-t)^\lz r^2|Z_ku|^2\lesssim\sum_{|K|\le 3}\int_{\Sigma^{\rm{ex}}_t}(2+r-t)^{\lz+1}|\partial Z^{K}u|^2{\rm{d}}x,\\
&&\sup_{\Sigma^{\rm{ex}}_t}(2+r-t)^{\lz-1} r^2|L_0u|^2\lesssim\sum_{|K|\le 4}\int_{\Sigma^{\rm{ex}}_t}(2+r-t)^{\lz+1}|\partial \Gamma^{K}u|^2{\rm{d}}x,
\een
where we recall \eqref{s2: Ga^I} and \eqref{s2: Z^J}.
\begin{proof}
By Lemmas \ref{s2: Sobex} and \ref{s2: exL_auL^2}, we have
\ben
\sup_{\Sigma^{\rm{ex}}_t}(2+r-t)^\lz r^2|Z_ku|^2&\lesssim&\sum_{|J|\le 2}\int_{\Sigma^{\rm{ex}}_t}\l\{(2+r-t)^{\lz+1}|\partial_r\Omega^JZ_ku|^2+(2+r-t)^{\lz-1}|\Omega^JZ_ku|^2\r\}{\rm{d}}x\\
&\lesssim&\sum_{|J|\le 2}\int_{\Sigma^{\rm{ex}}_t}(2+r-t)^{\lz+1}|\partial_r\Omega^JZ_ku|^2{\rm{d}}x+\sum_{|K|\le 3}\int_{\Sigma^{\rm{ex}}_t}(2+r-t)^{\lz+1}|\partial Z^{K}u|^2{\rm{d}}x\\
&\lesssim&\sum_{|K|\le 3}\int_{\Sigma^{\rm{ex}}_t}(2+r-t)^{\lz+1}|\partial Z^{K}u|^2{\rm{d}}x,
\een
\ben
\sup_{\Sigma^{\rm{ex}}_t}(2+r-t)^{\lz-1} r^2|L_0u|^2&\lesssim&\sum_{|J|\le 2}\int_{\Sigma^{\rm{ex}}_t}(2+r-t)^{\lz-1}\big\{|\partial_r\Omega^JL_0u|^2+|\Omega^JL_0u|^2\big\}{\rm{d}}x\\
&\lesssim&\sum_{|J|\le 2}\int_{\Sigma^{\rm{ex}}_t}(2+r-t)^{\lz-1}\big\{|L_0\partial\Omega^Ju|^2+|\partial\Omega^Ju|^2+|L_0\Omega^Ju|^2\big\}{\rm{d}}x\\
&\lesssim&\sum_{|K|\le 4}\int_{\Sigma^{\rm{ex}}_t}(2+r-t)^{\lz+1}|\partial \Gamma^{K}u|^2{\rm{d}}x,
\een
where we use the estimates on commutators in \eqref{s2: commutator}.
\end{proof}
\end{lem}

We give the following Sobolev inequality on hyperboloids; see \cite[Proposition 2.2]{LM18}.

\begin{lem}\label{s2: Sobin}
For any $s\in[2,\infty)$, any sufficiently smooth function $u$ defined in the region $\mathscr{H}_{[2,s]}:=\{(t,x)\in[2,\infty)\times\mathbb{R}^3: 2^2\le t^2-|x|^2\le s^2\}$ and any $\tau\in[2,s]$, we have
\be
\sup_{\mathscr{H}^{\rm{in}}_\tau}t^{3/2}|u(t,x)|\lesssim \sum_{|J|\le 2}\|L^Ju\|_{L^2(\mathscr{H}_\tau)},
\ee
where we recall \eqref{s2: pOL^J} and \eqref{s2: L^pH_s}.
\end{lem}

\subsection{Estimates on null forms and extra decay of the Klein-Gordon component}

We next give estimates of the null form 
$$\partial_\az \phi\partial^\az \psi=-\partial_t\phi\partial_t\psi+\partial_a\phi\partial^a\psi.$$

\begin{lem}\label{s2: Q_0}
For any sufficiently smooth functions $\phi$ and $\psi$, we have 
\ben
&&|\partial_\az \phi\partial^\az \psi|\lesssim\f{1}{t}\big(|\partial_t\phi|\cdot|L_0\psi|+|L_a\phi|\cdot|\partial^a\psi|\big),\\
&&|\partial_\az \phi\partial^\az \psi|\lesssim\f{|t-r|}{t}|\partial \phi|\cdot|\partial\psi|+\sum_{a=1}^3\big(|\bar{\partial}_a \phi|\cdot|\partial\psi|+|\partial\phi|\cdot|\bar{\partial}_a\psi|\big),\\
&&|\Gamma_k\big(\partial_\az \phi\partial^\az \psi\big)|\lesssim|(\partial_\az\Gamma_k \phi)(\partial^\az\psi)|+|(\partial_\az\phi)(\partial^\az\Gamma_k\psi)|,\quad k\in\{1,\cdots,10\},
\een
where we recall the definition \eqref{s2: gamma_k}.
\begin{proof}
We only prove the first and second inequalities. For the last one, see for example \cite{Sog}. Using the relation $\partial_a=\bar{\partial}_a-(x_a/t)\partial_t$ where $\bar{\partial}_a=t^{-1}L_a$, we have
\ben
-\partial_\az \phi\partial^\az \psi&=&\f{1}{t}\partial_t\phi\l(t\partial_t\psi+x^a\partial_a\psi-x^a\partial_a\psi\r)-\l(\bar{\partial}_a\phi-\f{x_a}{t}\partial_t\phi\r)\partial^a\psi\\
&=&\f{1}{t}\partial_t\phi\l((t-r)\partial_t\psi+(r-t)\partial_r\psi+\f{x^a}{r}L_a\psi\r)-\bar{\partial}_a\phi\partial^a\psi\\
&=&\f{t-r}{t}\partial_t\phi\partial_t\psi+\f{r-t}{t}\partial_t\phi\partial_r\psi+\f{x^a}{r}\partial_t\phi\bar{\partial}_a\psi-\bar{\partial}_a\phi\partial^a\psi.
\een
\end{proof}
\end{lem}

\begin{lem}\label{s2: Edecay}
For any sufficiently smooth function $v$, we have
\ben
&&|v|\lesssim\f{t-r}{t}\sum_{|I|\le 1}|\partial\Gamma^Iv|+|-\Box v+v|\quad\quad\mathrm{in}\;\mathscr{D}^{\rm{in}},\\
&&|v|\lesssim\f{2+r-t}{t}\sum_{|I|\le 1}|\partial\Gamma^Iv|+|-\Box v+v|\quad\quad\mathrm{in}\;\mathscr{D}^{\rm{ex}}\cap\{r\le 2t\}.
\een

\begin{proof}
We write the d'Alembert operator $-\Box$ as
\be
-\Box=\f{(t-r)(t+r)}{t^2}\partial_t\partial_t+\f{x^a}{t^2}\partial_tL_a-\f{1}{t}\partial^aL_a+\f{3}{t}\partial_t-\f{x^a}{t^2}\partial_a.
\ee
In $\mathscr{D}^{\rm{in}}$ we have $1\le t-r$, while in the region $\mathscr{D}^{\rm{ex}}$, we recall \eqref{s2: Edecay0}. Hence the lemma follows.
\end{proof}
\end{lem}

\section{Global existence in the exterior region}\label{s3}

In this section we prove global existence of solutions to \eqref{s3: equv}-\eqref{s1: ini} in the exterior region.

\subsection{Bootstrap assumptions}\label{ss3.1}
Fix $N\in\mathbb{N}$ with $N\ge 9$, $0<\dz\ll 1$ and $\lz\ge 3$. Let $C_1\gg1$ and $0<\ez\ll C_1^{-1}$ be chosen later. We assume the following bootstrap setting for the solution $(u,v)$ to \eqref{s3: equv}-\eqref{s1: ini}:
\beq\label{s3: bsex}
\sum_{|I|\le N}\big\{\|\Gamma^Iu\|_{X^{{\rm{ex}},\lz}_{0,\dz,t}}+\|\Gamma^Iv\|_{X^{{\rm{ex}},\lz}_{1,0,t}}\big\}+\sum_{|I|\le N-1}\|\Gamma^Iu\|_{X^{{\rm{ex}},\lz}_{0,0,t}}\le C_1\ez,
\eeq
\beq\label{s3: bsex'}
\sum_{|I|\le N}\|t^{-\f{\dz}{2}}(2+r-t)^{\f{\lz-1}{2}}\Gamma^Iu\|_{L^2_x(\Sigma^{\rm{ex}}_t)}+\sum_{|I|\le N-1}\|(2+r-t)^{\f{\lz-1}{2}}\Gamma^Iu\|_{L^2_x(\Sigma^{\rm{ex}}_t)}\le C_1\ez,
\eeq
where (see \eqref{s2: exnorm})
\beq\label{s3: gzudz}
\|\Gamma^Iu\|_{X^{{\rm{ex}},\lz}_{0,\dz,t}}=[E^{{\rm{ex}},\lz}_{0,\dz}(t,\Gamma^Iu)]^{\f{1}{2}}+\l(\int_2^t\int_{\Sigma^{\rm{ex}}_\tau}\tau^{-\dz}(2+r-\tau)^\lz|G\Gamma^Iu|^2{\rm{d}}x{\rm{d}}\tau\r)^{\f{1}{2}}+\l(\int_{\mathscr{l}_{[2,t]}}\tau^{-\dz}|G\Gamma^Iu|^2{\rm{d}}\sz\r)^{\f{1}{2}},
\eeq
\beqn\label{s3: gzv0}
\|\Gamma^Iv\|_{X^{{\rm{ex}},\lz}_{1,0,t}}=[E^{{\rm{ex}},\lz}_{1}(t,\Gamma^Iv)]^{\f{1}{2}}&+&\l(\int_2^t\int_{\Sigma^{\rm{ex}}_\tau}(2+r-\tau)^\lz\big(|G\Gamma^Iv|^2+|\Gamma^Iv|^2\big){\rm{d}}x{\rm{d}}\tau\r)^{\f{1}{2}}\nonumber\\
&+&\l(\int_{\mathscr{l}_{[2,t]}}\big(|G\Gamma^Iv|^2+|\Gamma^Iv|^2\big){\rm{d}}\sz\r)^{\f{1}{2}},
\eeqn
\beq\label{s3: gzu0}
\|\Gamma^Iu\|_{X^{{\rm{ex}},\lz}_{0,0,t}}=[E^{{\rm{ex}},\lz}_{0}(t,\Gamma^Iu)]^{\f{1}{2}}+\l(\int_2^t\int_{\Sigma^{\rm{ex}}_\tau}(2+r-\tau)^\lz|G\Gamma^Iu|^2{\rm{d}}x{\rm{d}}\tau\r)^{\f{1}{2}}
+\l(\int_{\mathscr{l}_{[2,t]}}|G\Gamma^Iu|^2{\rm{d}}\sz\r)^{\f{1}{2}},
\eeq
and we recall \eqref{dexfs} and \eqref{dexfsaz} for the definitions of the energy functionals above.

Let 
\beq\label{s3: T^*}
T^*:=\sup\{T>2: \eqref{s3: bsex}-\eqref{s3: bsex'} \mathrm{\ hold\ for\ }t\in[2,T]\}.
\eeq

\begin{prop}\label{s3: maxT}
There exists some constants $C_1>0$ sufficiently large and $0<\ez_0\ll C_1^{-1}$ sufficiently small such that, for any $0<\ez<\ez_0$, if $(u,v)$ is a solution to \eqref{s3: equv}-\eqref{s1: ini} in a time interval $[2,T]$ and satisfies \eqref{s3: bsex}-\eqref{s3: bsex'}, then for $t\in[2,T]$ we have
\ben\label{s3: imbsex}
&&\sum_{|I|\le N}\big\{\|\Gamma^Iu\|_{X^{{\rm{ex}},\lz}_{0,\dz,t}}+\|\Gamma^Iv\|_{X^{{\rm{ex}},\lz}_{1,0,t}}\big\}+\sum_{|I|\le N-1}\|\Gamma^Iu\|_{X^{{\rm{ex}},\lz}_{0,0,t}}\le \f{1}{2}C_1\ez,\nonumber\\
&&\sum_{|I|\le N}\|t^{-\f{\dz}{2}}(2+r-t)^{\f{\lz-1}{2}}\Gamma^Iu\|_{L^2_x(\Sigma^{\rm{ex}}_t)}+\sum_{|I|\le N-1}\|(2+r-t)^{\f{\lz-1}{2}}\Gamma^Iu\|_{L^2_x(\Sigma^{\rm{ex}}_t)}\le \f{1}{2}C_1\ez.
\een
\end{prop}

In the above proposition $T$ is arbitrary, hence the solution $(u,v)$ exists globally in time in $\mathscr{D}^{\rm{ex}}$ (i.e., $T^*=\infty$ where $T^*$ is as in \eqref{s3: T^*}) and satisfies \eqref{s3: bsex}-\eqref{s3: bsex'} for all $t\in[2,\infty)$. Below we give the proof of Proposition \ref{s3: maxT}. Let 
\ben
l(\tau):=(C_1\ez)^{-2}\sum_{|I|\le N}\int_{\Sigma^{\rm{ex}}_\tau}(2+r-\tau)^\lz|\Gamma^Iv|^2{\rm{d}}x.
\een
By \eqref{s3: bsex} and \eqref{s3: gzv0}, we have $l\in L^1([2,t])$ and $\|l\|_{L^1([2,t])}\lesssim 1$. Using \eqref{s3: bsex}-\eqref{s3: gzu0}, \eqref{dexfs}, \eqref{dexfsaz} and Lemma \ref{s2: exL_auL^2}, we obtain the following $L^2$-type estimates for the solution $(u,v)$ on the time interval $[2,T]$:
\beq\label{s3: exuvL^2N}
\sum_{|I|\le N}\|(2+r-t)^{\f{1+\lz}{2}}\big(t^{-\f{\dz}{2}}|\partial\Gamma^Iu|+|\partial\Gamma^Iv|+|\Gamma^Iv|\big)\|_{L^2_x(\Sigma^{\rm{ex}}_t)}\lesssim C_1\ez,
\eeq
\beq\label{s3: exuvL^2N'}
\sum_{|I|\le N}\|(2+r-\tau)^{\f{\lz}{2}}\Gamma^Iv\|_{L^2_x(\Sigma^{\rm{ex}}_\tau)}\lesssim C_1\ez\sqrt{l(\tau)},\quad \tau\in[2,t],
\eeq
\beq\label{s3: exuvL^2N''}
\sum_{|I|\le N}\|t^{-\f{\dz}{2}}(2+r-t)^{\f{\lz-1}{2}}\Gamma^Iu\|_{L^2_x(\Sigma^{\rm{ex}}_t)}+\sum_{|I|\le N-1}\|t^{-\f{\dz}{2}}(2+r-t)^{\f{\lz-1}{2}}L_0\Gamma^Iu\|_{L^2_x(\Sigma^{\rm{ex}}_t)}\lesssim C_1\ez,
\eeq
\beq\label{s3: exuvL^2N-1}
\sum_{|I|\le N-1}\|(2+r-t)^{\f{1+\lz}{2}}|\partial\Gamma^Iu|+(2+r-t)^{\f{\lz-1}{2}}|\Gamma^Iu|\|_{L^2_x(\Sigma^{\rm{ex}}_t)}\lesssim C_1\ez,
\eeq
\beq\label{s3: exuvL^2N-1'}
\sum_{|I|\le N-2}\|(2+r-t)^{\f{\lz-1}{2}}L_0\Gamma^Iu\|_{L^2_x(\Sigma^{\rm{ex}}_t)}\lesssim C_1\ez.
\eeq
By Lemmas \ref{s2: Sobex} and \ref{s2: exL_auL^infty}, we derive the following pointwise estimates for $t\in [2,T]$:
\beq\label{s3: exuvL^infty}
\sum_{|I|\le N-4}\sup_{\Sigma^{\rm{ex}}_t}(2+r-t)^{\f{1+\lz}{2}}r|\partial\Gamma^Iu|+\sum_{|I|\le N-3}\sup_{\Sigma^{\rm{ex}}_t}(2+r-t)^{\f{1+\lz}{2}}r\big(|\partial\Gamma^Iv|+|\Gamma^Iv|\big)\lesssim C_1\ez,
\eeq
\beq\label{s3: exuvL^infty''}
\sum_{|I|\le N-3}\sup_{\Sigma^{\rm{ex}}_\tau}(2+r-\tau)^{\f{\lz}{2}}r|\Gamma^Iv|\lesssim C_1\ez\sqrt{l(\tau)},\quad \tau\in[2,t],
\eeq
\beq\label{s3: exuvL^infty'}
\sum_{|I|\le N-3}\sup_{\Sigma^{\rm{ex}}_t}(2+r-t)^{\f{\lz}{2}}r|\Gamma^Iu|+\sum_{|I|\le N-5}\sup_{\Sigma^{\rm{ex}}_t}(2+r-t)^{\f{\lz-1}{2}}r|L_0\Gamma^Iu|\lesssim C_1\ez,
\eeq
where we use \eqref{s2: commutator}. We recall the second equation in \eqref{s3: equv}. By Lemma \ref{s2: Edecay}, in $\Sigma^{\rm{ex}}_t\cap\{r\le 2t\}$ we have
\ben
\sum_{|I|\le N-4}|\Gamma^Iv|&\lesssim&\f{2+r-t}{t}\sum_{|I|\le N-3}|\partial\Gamma^Iv|+\sum_{|I|\le N-4}|\Gamma^I(Q_0uv+Q_1^\az\partial_\az uv)|\\
&\lesssim&\f{2+r-t}{r}\sum_{|I|\le N-3}|\partial\Gamma^Iv|+C_1\ez t^{-1}\sum_{|I|\le N-4}|\Gamma^Iv|,
\een
where we use \eqref{s3: exuvL^infty'} and \eqref{s3: exuvL^infty}, which implies 
\be
\sum_{|I|\le N-4}|\Gamma^Iv|\lesssim\f{2+r-t}{r}\sum_{|I|\le N-3}|\partial\Gamma^Iv|\quad\quad\mathrm{in}\;\Sigma^{\rm{ex}}_t\cap\{r\le 2t\}.
\ee
It follows that
\beq\label{s3: exvextra}
\sum_{|I|\le N-4}\sup_{\Sigma^{\rm{ex}}_t\cap\{r\le 2t\}}|(2+r-t)r^2\Gamma^Iv|\lesssim\sum_{|I|\le N-3}\sup_{\Sigma^{\rm{ex}}_t\cap\{r\le 2t\}}|(2+r-t)^2r\partial\Gamma^Iv|\lesssim C_1\ez,
\eeq
where we use \eqref{s3: exuvL^infty} and that $\f{1+\lz}{2}\ge 2$ (we recall the assumption that $\lz\ge 3$). On the other hand, in the region $\{r\ge 2t\}$, we have $2+r-t\sim r$, hence using \eqref{s3: exuvL^infty} again, we obtain
\beq\label{s3: exvextra'}
\sum_{|I|\le N-4}\sup_{\Sigma^{\rm{ex}}_t\cap\{r\ge 2t\}}|(2+r-t)r^2\Gamma^Iv|\lesssim\sum_{|I|\le N-4}\sup_{\Sigma^{\rm{ex}}_t\cap\{r\ge 2t\}}|(2+r-t)^2r\Gamma^Iv|\lesssim C_1\ez.
\eeq
Combining \eqref{s3: exvextra} and \eqref{s3: exvextra'}, we derive
\beq\label{s3: exvextra''}
\sum_{|I|\le N-4}\sup_{\Sigma^{\rm{ex}}_t}|(2+r-t)r^2\Gamma^Iv|\lesssim C_1\ez.
\eeq

\subsection{Improved estimates for the solution $(u,v)$ in the exterior region}\label{ss3.2}

By Proposition \ref{s2: exfs},
\ben
\|\Gamma^Iu\|_{X^{{\rm{ex}},\lz}_{0,\dz,t}}&\lesssim&\|\Gamma^Iu\|_{X^{{\rm{ex}},\lz}_{0,\dz,2}}+\int_2^t\|\tau^{-\f{\dz}{2}}(2+r-\tau)^{\f{1+\lz}{2}}\Gamma^IF_u\|_{L^2_x(\Sigma^{\rm{ex}}_\tau)}{\rm{d}}\tau,\quad |I|\le N,\\
\|\Gamma^Iv\|_{X^{{\rm{ex}},\lz}_{1,0,t}}&\lesssim&\|\Gamma^Iv\|_{X^{{\rm{ex}},\lz}_{1,0,2}}+\int_2^t\|(2+r-\tau)^{\f{1+\lz}{2}}\Gamma^IF_v\|_{L^2_x(\Sigma^{\rm{ex}}_\tau)}{\rm{d}}\tau,\quad |I|\le N,\\
\|\Gamma^Iu\|_{X^{{\rm{ex}},\lz}_{0,0,t}}&\lesssim&\|\Gamma^Iu\|_{X^{{\rm{ex}},\lz}_{0,0,2}}+\int_2^t\|(2+r-\tau)^{\f{1+\lz}{2}}\Gamma^IF_u\|_{L^2_x(\Sigma^{\rm{ex}}_\tau)}{\rm{d}}\tau,\quad |I|\le N-1,
\een
where $F_u$ and $F_v$ are as in \eqref{s3: equv}. For $|I|\le N$ and $\tau\in[2,t]$, we have
\beqn\label{ScaexvN}
&&\|(2+r-\tau)^{\f{1+\lz}{2}}\big(|\Gamma^I(uv)|+|\Gamma^I(\partial uv)|\big)\|_{L^2_x(\Sigma^{\rm{ex}}_\tau)}\nonumber\\
&\lesssim&\sum_{\substack{|I_1|\le N-4\\|I_2|\le N}}\|(2+r-\tau)^{\f{\lz}{2}}r\big(|\Gamma^{I_1}u|+|\partial\Gamma^{I_1}u|\big)\|_{L^\infty_x(\Sigma^{\rm{ex}}_\tau)}\cdot\|(2+r-\tau)^{\f{1}{2}}\Gamma^{I_2}v\|_{L^2_x(\Sigma^{\rm{ex}}_\tau)}\cdot\tau^{-1}\nonumber\\
&+&\sum_{\substack{|I_2|\le N-4\\|I_1|\le N}}\|\tau^{-\f{\dz}{2}}(2+r-\tau)^{\f{\lz-1}{2}}\big(|\Gamma^{I_1}u|+|\partial\Gamma^{I_1}u|\big)\|_{L^2_x(\Sigma^{\rm{ex}}_\tau)}\cdot\|(2+r-\tau)r\tau\Gamma^{I_2}v\|_{L^\infty_x(\Sigma^{\rm{ex}}_\tau)}\cdot\tau^{-2+\f{\dz}{2}}\nonumber\\
&\lesssim&(C_1\ez)^2\big\{\sqrt{l(\tau)}\tau^{-1}+\tau^{-2+\f{\dz}{2}}\big\}
\eeqn
and
\ben
&&\|\tau^{-\f{\dz}{2}}(2+r-\tau)^{\f{1+\lz}{2}}\big(|\Gamma^I(u\partial v)|+|\Gamma^I(\partial u\partial v)|\big)\|_{L^2_x(\Sigma^{\rm{ex}}_\tau)}\nonumber\\
&\lesssim&\sum_{\substack{|I_1|\le N-4\\|I_2|\le N}}\|r\big(|\Gamma^{I_1}u|+|\partial\Gamma^{I_1}u|\big)\|_{L^\infty_x(\Sigma^{\rm{ex}}_\tau)}\cdot\|(2+r-\tau)^{\f{1+\lz}{2}}\partial\Gamma^{I_2}v\|_{L^2_x(\Sigma^{\rm{ex}}_\tau)}\cdot\tau^{-1-\f{\dz}{2}}\nonumber\\
&+&\sum_{\substack{|I_2|\le N-5\\|I_1|\le N}}\|\tau^{-\f{\dz}{2}}(2+r-\tau)^{\f{\lz-1}{2}}\big(|\Gamma^{I_1}u|+|\partial\Gamma^{I_1}u|\big)\|_{L^2_x(\Sigma^{\rm{ex}}_\tau)}\cdot\|(2+r-\tau)r\tau\partial\Gamma^{I_2}v\|_{L^\infty_x(\Sigma^{\rm{ex}}_\tau)}\cdot\tau^{-2}\nonumber\\
&\lesssim&(C_1\ez)^2\tau^{-1-\f{\dz}{2}},
\een
where we use \eqref{s3: exuvL^infty'}, \eqref{s3: exuvL^infty}, \eqref{s3: exuvL^2N'}, \eqref{s3: exuvL^2N''}, \eqref{s3: exuvL^2N} and \eqref{s3: exvextra''}. For $|I|\le N-1$ and $\tau\in[2,t]$, we have
\beqn\label{ScaexdvN-1}
&&\|(2+r-\tau)^{\f{1+\lz}{2}}\big(|\Gamma^I(u\partial v)|+|\Gamma^I(\partial u\partial v)|\big)\|_{L^2_x(\Sigma^{\rm{ex}}_\tau)}\nonumber\\
&\lesssim&\sum_{\substack{|I_1|\le N-4\\|I_2|\le N}}\|(2+r-t)^{\f{\lz}{2}}r\big(|\Gamma^{I_1}u|+|\partial\Gamma^{I_1}u|\big)\|_{L^\infty_x(\Sigma^{\rm{ex}}_\tau)}\cdot\|(2+r-\tau)^{\f{1}{2}}\Gamma^{I_2}v\|_{L^2_x(\Sigma^{\rm{ex}}_\tau)}\cdot\tau^{-1}\nonumber\\
&+&\sum_{\substack{|I_2|\le N-5\\|I_1|\le N-1}}\|(2+r-\tau)^{\f{\lz-1}{2}}\big(|\Gamma^{I_1}u|+|\partial\Gamma^{I_1}u|\big)\|_{L^2_x(\Sigma^{\rm{ex}}_\tau)}\cdot\|(2+r-\tau)r\tau\partial\Gamma^{I_2}v\|_{L^\infty_x(\Sigma^{\rm{ex}}_\tau)}\cdot\tau^{-2}\nonumber\\
&\lesssim&(C_1\ez)^2\big\{\sqrt{l(\tau)}\tau^{-1}+\tau^{-2}\big\},
\eeqn
where we use \eqref{s3: exuvL^infty'}, \eqref{s3: exuvL^infty}, \eqref{s3: exuvL^2N'}, \eqref{s3: exuvL^2N-1} and \eqref{s3: exvextra''}. It follows that
\beqn\label{s3: S1gzF_u}
\int_2^t\Bigg\{\sum_{|I|\le N}\|(2+r-\tau)^{\f{1+\lz}{2}}\big(\tau^{-\f{\dz}{2}}|\Gamma^IF_u|+|\Gamma^IF_v|\big)\|_{L^2_x(\Sigma^{\rm{ex}}_\tau)}\!\!&+&\!\!\sum_{|I|\le N-1}\|(2+r-\tau)^{\f{1+\lz}{2}}\Gamma^IF_u\|_{L^2_x(\Sigma^{\rm{ex}}_\tau)}\Bigg\}{\rm{d}}\tau\nonumber\\
&\lesssim&(C_1\ez)^2
\eeqn
and therefore
\beq\label{s3: S1gzuN}
\sum_{|I|\le N}\big\{\|\Gamma^Iu\|_{X^{{\rm{ex}},\lz}_{0,\dz,t}}+\|\Gamma^Iv\|_{X^{{\rm{ex}},\lz}_{1,0,t}}\big\}+\sum_{|I|\le N-1}\|\Gamma^Iu\|_{X^{{\rm{ex}},\lz}_{0,0,t}}\lesssim\ez+(C_1\ez)^2.
\eeq
By Lemma \ref{s2: Hardyex} and the definitions \eqref{s3: gzudz}, \eqref{s3: gzu0}, \eqref{dexfs} and \eqref{dexfsaz},
\beq\label{s3: uL^2imN}
\|t^{-\f{\dz}{2}}(2+r-t)^{\f{\lz-1}{2}}\Gamma^Iu\|_{L^2_x(\Sigma^{\rm{ex}}_t)}\lesssim\|t^{-\f{\dz}{2}}(2+r-t)^{\f{1+\lz}{2}}\partial\Gamma^Iu\|_{L^2_x(\Sigma^{\rm{ex}}_t)}\lesssim\ez+(C_1\ez)^2,\quad |I|\le N,
\eeq
\beq\label{s3: uL^2im}
\|(2+r-t)^{\f{\lz-1}{2}}\Gamma^Iu\|_{L^2_x(\Sigma^{\rm{ex}}_t)}\lesssim\|(2+r-t)^{\f{1+\lz}{2}}\partial\Gamma^Iu\|_{L^2_x(\Sigma^{\rm{ex}}_t)}\lesssim\ez+(C_1\ez)^2,\quad |I|\le N-1.
\eeq
Hence we have strictly improved the bootstrap estimates \eqref{s3: bsex}-\eqref{s3: bsex'} for $C_1$ sufficiently large and $0<\ez\ll C_1^{-1}$ sufficiently small. Therefore the proof of Proposition \ref{s3: maxT} is completed, and we obtain the global existence of the solution $(u,v)$ in the exterior region $\mathscr{D}^{\rm{ex}}$ with the global-in-time estimates \eqref{s3: bsex}-\eqref{s3: bsex'}.

We recall the definitions \eqref{s2: E^{c,ex}} and \eqref{s2: E^{c,ex}_0}. Using Lemma \ref{s2: exL_auL^2},  \eqref{s3: uL^2imN} and \eqref{s3: uL^2im}, we obtain
\beq\label{s3: S4cgzuN-12}
\sum_{|I|\le N-1}[E^{c,\rm{ex}}_\dz(t,\Gamma^Iu)]^{\f{1}{2}}+\sum_{|I|\le N-2}[E^{c,\rm{ex}}(t,\Gamma^Iu)]^{\f{1}{2}}\lesssim\ez+(C_1\ez)^2.
\eeq

\subsection{Energy estimates on exterior hyperboloids}\label{ss3.3}

We next derive energy estimates on truncated exterior hyperboloids. 

${\bf Step\ 1.}$ Energy estimates of $(u,v)$ on exterior hyperboloids: By Proposition \ref{s2: exhs}, for all $s\in[2,\infty)$, we have
\beqn\label{s3: exhuN}
\!\!\!\!\!\!\!\!\!\!\!\!&&\!\!\!\sum_{|I|\le N}\Bigg\{E^{{\rm{ex}},h}_{0,\dz}(s,\Gamma^Iu)+\int_{\mathscr{l}_{[2,\f{s^2+1}{2}]}}t^{-\dz}|G\Gamma^Iu|^2{\rm{d}}\sz\Bigg\}\nonumber\\
\!\!\!\!\!\!\!\!\!\!\!\!&\lesssim&\!\!\!\sum_{|I|\le N}\Bigg\{E^{\rm{ex}}_{0,\dz}(2,\Gamma^Iu)+\sup_{t\in[2,\infty)}[E^{\rm{ex}}_{0,\dz}(t,\Gamma^Iu)]^{\f{1}{2}}\cdot\int_2^\infty\|t^{-\f{\dz}{2}}\Gamma^IF_u\|_{L^2_x(\Sigma^s_t)}{\rm{d}}t\Bigg\}\lesssim\ez^2+(C_1\ez)^3,
\eeqn
where we use \eqref{s3: S1gzF_u} and \eqref{s3: S1gzuN} (we recall \eqref{s3: gzudz}) and that $\Sigma^s_t\subset\Sigma^{\rm{ex}}_t$ by \eqref{s2: Sigma^s_t}, and recall \eqref{dexhs} and \eqref{dexfsdz} for the definitions of $E^{{\rm{ex}},h}_{0,\dz}(s,\Gamma^Iu)$ and $E^{\rm{ex}}_{0,\dz}(t,\Gamma^Iu)$ respectively. Similarly, by \eqref{s3: S1gzF_u} and \eqref{s3: S1gzuN} (here we recall \eqref{s3: gzv0} and \eqref{s3: gzu0}), we also have
\beq\label{s3: exhvN}
\begin{split}
&E^{{\rm{ex}},h}_{1}(s,\Gamma^Iv)+\int_{\mathscr{l}_{[2,\f{s^2+1}{2}]}}\big(|G\Gamma^Iv|^2+|\Gamma^Iv|^2\big){\rm{d}}\sz\\
\lesssim&\ \  E^{\rm{ex}}_1(2,\Gamma^Iv)+\sup_{t\in[2,\infty)}[E^{\rm{ex}}_1(t,\Gamma^Iv)]^{\f{1}{2}}\cdot\int_2^\infty\|\Gamma^IF_v\|_{L^2_x(\Sigma^s_t)}{\rm{d}}t\lesssim\ez^2+(C_1\ez)^3,\quad |I|\le N,
\end{split}
\eeq
\beq\label{s3: exhuN-1}
\begin{split}
&E^{{\rm{ex}},h}_{0}(s,\Gamma^Iu)+\int_{\mathscr{l}_{[2,\f{s^2+1}{2}]}}|G\Gamma^Iu|^2{\rm{d}}\sz\\
\lesssim&\ \  E^{\rm{ex}}_{0}(2,\Gamma^Iu)+\sup_{t\in[2,\infty)}[E^{\rm{ex}}_{0}(t,\Gamma^Iu)]^{\f{1}{2}}\cdot\int_2^\infty\|\Gamma^IF_u\|_{L^2_x(\Sigma^s_t)}{\rm{d}}t\lesssim\ez^2+(C_1\ez)^3,\quad |I|\le N-1,
\end{split}
\eeq
where we recall \eqref{dhs0} and \eqref{dexfs0} for the definitions of $E^{{\rm{ex}},h}_{1}(s,\Gamma^Iv)$, $E^{{\rm{ex}},h}_{0}(s,\Gamma^Iu)$, $E^{\rm{ex}}_{1}(t,\Gamma^Iv)$ and $E^{\rm{ex}}_{0}(t,\Gamma^Iu)$.

${\bf Step\ 2.}$ Conformal energy estimates of $u$ on exterior hyperboloids: 

We first point out that, by the global existence result in Section \ref{ss3.2}, the estimates \eqref{s3: exuvL^2N}-\eqref{s3: exuvL^infty'} and \eqref{s3: exvextra''} hold for all $t\in[2,\infty)$.

${\bf Step\ 2.1.}$ A nonlinear transform: Let 
\beq\label{s3: defU}
U:=u+F_u=u+P_{0}uv+P_{1}^\az\partial_\az uv+P_{0}^\az u\partial_\az v+P_1^{\az\bz}\partial_\az u\partial_\bz v.
\eeq
Then $U$ satisfies (see \eqref{s3: equv})
\beq\label{s3: defBoxU}
\begin{split}
&-\Box U=P_{0}uv+P_{1}^\az\partial_\az uv+P_{0}^\az u\partial_\az v+P_1^{\az\bz}\partial_\az u\partial_\bz v\\
&+\  P_0\big\{P_{0}uv+P_{1}^\az\partial_\az uv+P_{0}^\az u\partial_\az v+P_1^{\az\bz}\partial_\az u\partial_\bz v\big\}v+P_0u\big\{\!\!-v+Q_0uv+Q_1^\az\partial_\az uv\big\}-2P_0\partial_\gz u\partial^\gz v\\
&+\  P_{1}^\gz\partial_\gz\big(P_{0}uv+P_{1}^\az\partial_\az uv+P_{0}^\az u\partial_\az v+P_1^{\az\bz}\partial_\az u\partial_\bz v\big)v+P_{1}^\gz\partial_\gz u\big\{\!-v+Q_0uv+Q_1^\az\partial_\az uv\big\}-2P_1^\az\partial_\gz\partial_\az u\partial^\gz v\\
&+\  P_{0}^\gz\big\{P_{0}uv+P_{1}^\az\partial_\az uv+P_{0}^\az u\partial_\az v+P_1^{\az\bz}\partial_\az u\partial_\bz v\big\}\partial_\gz v+P_{0}^\gz u\partial_\gz\big(\!\!-v+Q_0uv+Q_1^\az\partial_\az uv\big)-2P_0^\az\partial_\gz u\partial^\gz\partial_\az v\\
&+\  P_1^{\az'\bz'}\partial_{\az'}\big(P_{0}uv+P_{1}^\az\partial_\az uv+P_{0}^\az u\partial_\az v+P_1^{\az\bz}\partial_\az u\partial_\bz v\big)\partial_{\bz'} v\\
&+\  P_1^{\az'\bz'}\partial_{\az'} u\partial_{\bz'}\big(\!\!-v+Q_0uv+Q_1^\az\partial_\az uv\big)-2P_1^{\az\bz}\partial_\gz\partial_\az u\partial^\gz\partial_\bz v\\
&=\  S+W,
\end{split}
\eeq
where $\az',\bz'\in\{0,1,2,3\}$, $S$ denotes the sum of the null forms, i.e.
\beq\label{s3: defS}
S:=-2\big\{P_0\partial_\gz u\partial^\gz v+P_1^\az\partial_\gz\partial_\az u\partial^\gz v+P_0^\az\partial_\gz u\partial^\gz\partial_\az v+P_1^{\az\bz}\partial_\gz\partial_\az u\partial^\gz\partial_\bz v\big\}
\eeq
and $W$ is a linear combination (with constant coefficients) of cubic terms, which is of the form  
\beqn\label{s3: defT}
W&\approx&uv^2+\partial uv^2+u\partial vv+\partial u\partial vv+u^2v+u\partial uv\nonumber\\
&+&\partial^2 uv^2+u\partial^2vv+\partial^2 u\partial vv+\partial u\partial^2vv+\partial u\partial uv\nonumber\\
&+&u\partial v\partial v+\partial u\partial v\partial v+ u^2\partial v+u\partial^2 uv+u\partial u\partial v\nonumber\\
&+&u\partial^2 v\partial v+\partial^2 u\partial v\partial v+\partial u\partial^2v\partial v+\partial u\partial^2 uv+\partial u\partial u\partial v,
\eeqn
(here "$\approx$" means that we omit the coefficients of these terms), where we denote by $\partial$ and $\partial^2$ any of the derivatives $\partial_\az,\az\in\{0,1,2,3\}$ and $\partial_\az\partial_\bz,\az,\bz\in\{0,1,2,3\}$ respectively.

By Proposition \ref{s2: excon}, for all $t\in[2,\infty)$, we have
\beq\label{s3: E_conU0}
\begin{split}
&[E^{c,\rm{ex}}_\dz(t,\Gamma^IU)]^{\f{1}{2}}+\mathcal{R}_\dz(t,\Gamma^IU)\lesssim[E^{c,\rm{ex}}_\dz(2,\Gamma^IU)]^{\f{1}{2}}+\int_2^t\|\tau^{-\f{\dz}{2}}r\Gamma^I(S+W)\|_{L^2_x(\Sigma^{\rm{ex}}_\tau)}{\rm{d}}\tau,\quad|I|\le N-1,\\
&[E^{c,\rm{ex}}(t,\Gamma^IU)]^{\f{1}{2}}+\mathcal{R}_0(t,\Gamma^IU)\lesssim[E^{c,\rm{ex}}(2,\Gamma^IU)]^{\f{1}{2}}+\int_2^t\|r\Gamma^I(S+W)\|_{L^2_x(\Sigma^{\rm{ex}}_\tau)}{\rm{d}}\tau,\quad|I|\le N-2,
\end{split}
\eeq
where we recall \eqref{s2: E^{c,ex}} and \eqref{s2: E^{c,ex}_0} for the definitions of $E^{c,\rm{ex}}_\dz(t,\Gamma^IU)$ and $E^{c,\rm{ex}}(t,\Gamma^IU)$, and
\beq\label{s3: S4H(t,x)}
\begin{split}
&\mathcal{R}_\dz(t,\Gamma^IU):=\l(\int_{\mathscr{l}_{[2,t]}}\tau^{-\dz}H[\Gamma^IU](\tau,x){\rm{d}}\sz\r)^{\f{1}{2}},\quad\quad\mathcal{R}_0(t,\Gamma^IU):=\l(\int_{\mathscr{l}_{[2,t]}}H[\Gamma^IU](\tau,x){\rm{d}}\sz\r)^{\f{1}{2}},\\
&H[\Gamma^IU](\tau,x):=|(\tau+r)(\partial_t\Gamma^IU+\partial_r\Gamma^IU)+2\Gamma^IU|^2+(\tau-r)^2r^{-2}|\Omega \Gamma^IU|^2.
\end{split}
\eeq
For $|I|\le N-1$ and $\tau\in[2,t]$, by Lemma \ref{s2: Q_0}, we have
\beq\label{s3: cgzudv}
\begin{split}
&\quad\quad\|\tau^{-\f{\dz}{2}}r\big(|\Gamma^I(\partial_\gz u\partial^\gz v)|+|\Gamma^I(\partial_\gz u\partial^\gz \partial v)|+|\Gamma^I(\partial_\gz\partial u\partial^\gz v)|+|\Gamma^I(\partial_\gz\partial u\partial^\gz \partial v)|\big)\|_{L^2_x(\Sigma^{\rm{ex}}_\tau)}\\
&\lesssim\ \sum_{\substack{|I_1|+|I_2|\le N-1\\|J|,|J'|=1}}\|\tau^{-\f{\dz}{2}}r\tau^{-1}\langle\tau-r\rangle\cdot\big(|\Gamma^J\Gamma^{I_1}u|+|\Gamma^J\Gamma^{I_1}\partial u|\big)\cdot\big(|\Gamma^{J'}\Gamma^{I_2}v|+|\Gamma^{J'}\Gamma^{I_2}\partial v|\big)\|_{L^2_x(\Sigma^{\rm{ex}}_\tau)}\\
&\lesssim\ \sum_{\substack{|I_1|\le N-5,|J|=1\\|I_2|\le N}}\|r\big(|\Gamma^J\Gamma^{I_1}u|+|\Gamma^J\Gamma^{I_1}\partial u|\big)\|_{L^\infty_x(\Sigma^{\rm{ex}}_\tau)}\cdot\|(2+r-\tau)\big(|\Gamma^{I_2}v|+|\partial\Gamma^{I_2}v|\big)\|_{L^2_x(\Sigma^{\rm{ex}}_\tau)}\cdot\tau^{-1-\f{\dz}{2}}\\
&+\ \sum_{\substack{|I_2|\le N-6,|J'|=1\\|I_1|\le N}}\|\tau^{-\f{\dz}{2}}\big(|\Gamma^{I_1}u|+|\partial\Gamma^{I_1}u|\big)\|_{L^2_x(\Sigma^{\rm{ex}}_\tau)}\cdot\|(2+r-\tau)r\tau\big(|\Gamma^{J'}\Gamma^{I_2}v|+|\Gamma^{J'}\Gamma^{I_2}\partial v|\big)\|_{L^\infty_x(\Sigma^{\rm{ex}}_\tau)}\cdot\tau^{-2}\\
&\lesssim\  (C_1\ez)^2\tau^{-1-\f{\dz}{2}},
\end{split}
\eeq
where we use \eqref{s2: Edecay0}, \eqref{s3: exuvL^infty'}, \eqref{s3: exuvL^infty}, \eqref{s3: exuvL^2N}, \eqref{s3: exuvL^2N''} and \eqref{s3: exvextra''} (recall that $\lz\ge 3$ and $N\ge 9$). It follows that 
\beq\label{s3: cgzS}
\int_2^t\|\tau^{-\f{\dz}{2}}r\Gamma^IS\|_{L^2_x(\Sigma^{\rm{ex}}_\tau)}{\rm{d}}\tau\lesssim (C_1\ez)^2,\quad |I|\le N-1.
\eeq
For $|I|\le N-2$, we argue as in \eqref{s3: cgzudv}, use \eqref{s3: exuvL^2N'} in the estimate of 
$$\sum_{|I_2|\le N-1}\|(2+r-\tau)\big(|\Gamma^{I_2}v|+|\partial\Gamma^{I_2}v|\big)\|_{L^2_x(\Sigma^{\rm{ex}}_\tau)},$$ 
and use \eqref{s3: exuvL^2N-1} in the estimate of 
$$\sum_{|I_1|\le N-1}\||\Gamma^{I_1}u|+|\partial\Gamma^{I_1}u|\|_{L^2_x(\Sigma^{\rm{ex}}_\tau)},$$ 
and derive
\be
\|r\Gamma^IS\|_{L^2_x(\Sigma^{\rm{ex}}_\tau)}\lesssim (C_1\ez)^2\big\{\sqrt{l(\tau)}\tau^{-1}+\tau^{-2}\big\},\quad |I|\le N-2,\tau\in[2,t],
\ee
which implies
\beq\label{s3: cgzSN-2}
\int_2^t\|r\Gamma^IS\|_{L^2_x(\Sigma^{\rm{ex}}_\tau)}{\rm{d}}\tau\lesssim (C_1\ez)^2,\quad |I|\le N-2.
\eeq

Next we turn to the estimate of $\int_2^t\|r\Gamma^IW\|_{L^2_x(\Sigma^{\rm{ex}}_\tau)}{\rm{d}}\tau$ for $|I|\le N-1$. We divide the terms appearing in the expression of $W$ (see \eqref{s3: defT}) into four classes:
\begin{itemize}
\item[$i)$] $uv^2$, $\partial uv^2$, $u\partial vv$, $\partial u\partial vv$, $u\partial^2vv$, $\partial u\partial^2vv$, $u\partial v\partial v$, $\partial u\partial v\partial v$, $\partial u\partial^2v\partial v$, $u\partial^2v\partial v$;
\item[$ii)$] $\partial^2uv^2$, $\partial^2u\partial vv$, $\partial^2u\partial v\partial v$;
\item[$iii)$] $u^2v$, $u\partial uv$, $\partial u\partial uv$, $u\partial u\partial v$, $\partial u\partial u\partial v$, $u^2\partial v$;
\item[$iv)$] $\partial u\partial^2uv$, $u\partial^2uv$.
\end{itemize}
In each class we pick the most difficult term, that is, 
\beq\label{s3: pic}
u\partial^2v\partial v, \partial^2u\partial v\partial v, u^2\partial v, u\partial^2uv
\eeq
respectively. We will see below that the estimates of the other terms in the corresponding classes are similar to or simpler than these chosen terms. For $|I|\le N-1$ and $\tau\in[2,t]$, we have
\ben
&&\|r\Gamma^I(u\partial^2 v\partial v)\|_{L^2_x(\Sigma^{\rm{ex}}_\tau)}\lesssim\sum_{|I_1|+|I_2|+|I_3|\le N-1}\|r|\Gamma^{I_1}u|\cdot|\partial\Gamma^{I_2}\partial v|\cdot|\partial\Gamma^{I_3}v|\|_{L^2_x(\Sigma^{\rm{ex}}_\tau)}\\
&\lesssim&\sum_{\substack{|I_1|,|I_2|\le N-6\\|I_3|\le N-1}}\|r\Gamma^{I_1}u\|_{L^\infty_x(\Sigma^{\rm{ex}}_\tau)}\cdot\|r\tau\partial\Gamma^{I_2}\partial v\|_{L^\infty_x(\Sigma^{\rm{ex}}_\tau)}\cdot\|\partial\Gamma^{I_3}v\|_{L^2_x(\Sigma^{\rm{ex}}_\tau)}\cdot\tau^{-2}\\
&+&\sum_{\substack{|I_1|,|I_3|\le N-5\\|I_2|\le N-1}}\|r\Gamma^{I_1}u\|_{L^\infty_x(\Sigma^{\rm{ex}}_\tau)}\cdot\|\partial\Gamma^{I_2}\partial v\|_{L^2_x(\Sigma^{\rm{ex}}_\tau)}\cdot\|r\tau\partial\Gamma^{I_3}v\|_{L^\infty_x(\Sigma^{\rm{ex}}_\tau)}\cdot\tau^{-2}\\
&+&\sum_{\substack{|I_2|,|I_3|\le N-5\\|I_1|\le N-1}}\|\Gamma^{I_1}u\|_{L^2_x(\Sigma^{\rm{ex}}_\tau)}\cdot\|r\partial\Gamma^{I_2}\partial v\|_{L^\infty_x(\Sigma^{\rm{ex}}_\tau)}\cdot\|r\tau\partial\Gamma^{I_3}v\|_{L^\infty_x(\Sigma^{\rm{ex}}_\tau)}\cdot\tau^{-2}\\
&\lesssim&(C_1\ez)^3\tau^{-2},
\een
where we use \eqref{s3: exuvL^infty'}, \eqref{s3: exvextra''}, \eqref{s3: exuvL^2N}, \eqref{s3: exuvL^2N-1} and \eqref{s3: exuvL^infty} and recall that $N\ge 9$. Similarly,
\ben
&&\|r\Gamma^I(\partial^2u\partial v\partial v)\|_{L^2_x(\Sigma^{\rm{ex}}_\tau)}\lesssim\sum_{|I_1|+|I_2|+|I_3|\le N-1}\|r|\partial\Gamma^{I_1}\partial u|\cdot|\partial\Gamma^{I_2}v|\cdot|\partial\Gamma^{I_3}v|\|_{L^2_x(\Sigma^{\rm{ex}}_\tau)}\\
&\lesssim&\sum_{\substack{|I_1|,|I_2|\le N-5\\|I_3|\le N-1}}\|r\partial\Gamma^{I_1}\partial u\|_{L^\infty_x(\Sigma^{\rm{ex}}_\tau)}\cdot\|r\tau\partial\Gamma^{I_2}v\|_{L^\infty_x(\Sigma^{\rm{ex}}_\tau)}\cdot\|\partial\Gamma^{I_3}v\|_{L^2_x(\Sigma^{\rm{ex}}_\tau)}\cdot\tau^{-2}\\
&+&\sum_{\substack{|I_2|,|I_3|\le N-5\\|I_1|\le N}}\|\tau^{-\f{\dz}{2}}\partial\Gamma^{I_1}u\|_{L^2_x(\Sigma^{\rm{ex}}_\tau)}\cdot\|r\tau\partial\Gamma^{I_2}v\|_{L^\infty_x(\Sigma^{\rm{ex}}_\tau)}\cdot\|r\tau\partial\Gamma^{I_3}v\|_{L^\infty_x(\Sigma^{\rm{ex}}_\tau)}\cdot\tau^{-3+\f{\dz}{2}}\\
&\lesssim&(C_1\ez)^3\tau^{-2},
\een
where we use \eqref{s3: exuvL^infty}, \eqref{s3: exvextra''} and \eqref{s3: exuvL^2N}. We also have 
\ben
&&\|r\Gamma^I(u^2\partial v)\|_{L^2_x(\Sigma^{\rm{ex}}_\tau)}\lesssim\sum_{|I_1|+|I_2|+|I_3|\le N-1}\|r|\Gamma^{I_1}u|\cdot|\Gamma^{I_2}u|\cdot|\partial\Gamma^{I_3}v|\|_{L^2_x(\Sigma^{\rm{ex}}_\tau)}\\
&\lesssim&\sum_{\substack{|I_1|,|I_2|\le N-3\\|I_3|\le N}}\|r\Gamma^{I_1}u\|_{L^\infty_x(\Sigma^{\rm{ex}}_\tau)}\cdot\|r\Gamma^{I_2}u\|_{L^\infty_x(\Sigma^{\rm{ex}}_\tau)}\cdot\|\Gamma^{I_3}v\|_{L^2_x(\Sigma^{\rm{ex}}_\tau)}\cdot \tau^{-1}\\
&+&\sum_{\substack{|I_2|,|I_3|\le N-5\\|I_1|\le N-1}}\|\Gamma^{I_1}u\|_{L^2_x(\Sigma^{\rm{ex}}_\tau)}\cdot\|r\Gamma^{I_2}u\|_{L^\infty_x(\Sigma^{\rm{ex}}_\tau)}\cdot\|r\tau\partial\Gamma^{I_3}v\|_{L^\infty_x(\Sigma^{\rm{ex}}_\tau)}\cdot \tau^{-2}\\
&\lesssim&(C_1\ez)^3\big\{\sqrt{l(\tau)}\tau^{-1}+\tau^{-2}\big\}
\een
and
\ben
&&\|r\Gamma^I(u\partial^2u v)\|_{L^2_x(\Sigma^{\rm{ex}}_\tau)}\lesssim\sum_{|I_1|+|I_2|+|I_3|\le N-1}\|r|\Gamma^{I_1}u|\cdot|\partial\Gamma^{I_2}\partial u|\cdot|\Gamma^{I_3}v|\|_{L^2_x(\Sigma^{\rm{ex}}_\tau)}\\
&\lesssim&\sum_{\substack{|I_1|,|I_2|\le N-5\\|I_3|\le N-1}}\|r\Gamma^{I_1}u\|_{L^\infty_x(\Sigma^{\rm{ex}}_\tau)}\cdot\|r\partial\Gamma^{I_2}\partial u\|_{L^\infty_x(\Sigma^{\rm{ex}}_\tau)}\cdot\|\Gamma^{I_3}v\|_{L^2_x(\Sigma^{\rm{ex}}_\tau)}\cdot \tau^{-1}\\
&+&\sum_{\substack{|I_1|,|I_3|\le N-4\\|I_2|\le N-1}}\|r\Gamma^{I_1}u\|_{L^\infty_x(\Sigma^{\rm{ex}}_\tau)}\cdot\|\tau^{-\f{\dz}{2}}\partial\Gamma^{I_2}\partial u\|_{L^2_x(\Sigma^{\rm{ex}}_\tau)}\cdot\|r\tau\Gamma^{I_3}v\|_{L^\infty_x(\Sigma^{\rm{ex}}_\tau)}\cdot \tau^{-2+\f{\dz}{2}}\\
&+&\sum_{\substack{|I_2|,|I_3|\le N-5\\|I_1|\le N-1}}\|\Gamma^{I_1}u\|_{L^2_x(\Sigma^{\rm{ex}}_\tau)}\cdot\|r\partial\Gamma^{I_2}\partial u\|_{L^\infty_x(\Sigma^{\rm{ex}}_\tau)}\cdot\|r\tau\Gamma^{I_3}v\|_{L^\infty_x(\Sigma^{\rm{ex}}_\tau)}\cdot \tau^{-2}\\
&\lesssim&(C_1\ez)^3\big\{\sqrt{l(\tau)}\tau^{-1}+\tau^{-2+\f{\dz}{2}}\big\},
\een
where we use \eqref{s3: exuvL^infty'}, \eqref{s3: exuvL^2N'}, \eqref{s3: exuvL^2N-1}, \eqref{s3: exvextra''}, \eqref{s3: exuvL^infty} and \eqref{s3: exuvL^2N}. Combining these estimates, we obtain 
\beq\label{s3: cgzTN-1}
\int_2^t\|r\Gamma^IW\|_{L^2_x(\Sigma^{\rm{ex}}_\tau)}{\rm{d}}\tau\lesssim (C_1\ez)^3,\quad |I|\le N-1.
\eeq
This together with \eqref{s3: E_conU0}, \eqref{s3: cgzS} and \eqref{s3: cgzSN-2} give 
\beq\label{s3: S4cgzUN-1}
[E^{c,\rm{ex}}_\dz(t,\Gamma^IU)]^{\f{1}{2}}+\mathcal{R}_\dz(t,\Gamma^IU)\lesssim\ez+(C_1\ez)^2,\quad |I|\le N-1,
\eeq
\beq\label{s3: S4cgzUN-2}
[E^{c,\rm{ex}}(t,\Gamma^IU)]^{\f{1}{2}}+\mathcal{R}_0(t,\Gamma^IU)\lesssim\ez+(C_1\ez)^2,\quad |I|\le N-2,
\eeq
for all $t\in[2,\infty)$, where we recall the definitions of $\mathcal{R}_\dz(t,\Gamma^IU)$ and $\mathcal{R}_0(t,\Gamma^IU)$ in \eqref{s3: S4H(t,x)}.

${\bf Step\ 2.2.}$ Estimates on exterior hyperboloids: For any $s\in[2,\infty)$, we denote $t(s)=\f{s^2+1}{2}$ (see \eqref{s2: tr(s)}). Using Proposition \ref{s2: exconh}, and the estimates \eqref{s3: cgzS}, \eqref{s3: cgzSN-2}, \eqref{s3: cgzTN-1}, \eqref{s3: S4cgzUN-1} and \eqref{s3: S4cgzUN-2}, we have
\beq\label{s3: S2.2U}
\begin{split}
[E^{c,{\rm{ex}},h}_{\dz}(s,\Gamma^IU)]^{\f{1}{2}}&\lesssim[E^{c,\rm{ex}}_{\dz}(t(s),\Gamma^IU)]^{\f{1}{2}}+\int_{t(s)}^\infty\|\tau^{-\f{\dz}{2}}r\Gamma^I(S+W)\|_{L^2_x(\Sigma^{\rm{ex}}_\tau)}{\rm{d}}\tau\\
&\lesssim\ez+(C_1\ez)^2,\quad |I|\le N-1,\\
[E^{c,{\rm{ex}},h}(s,\Gamma^IU)]^{\f{1}{2}}&\lesssim[E^{c,\rm{ex}}(t(s),\Gamma^IU)]^{\f{1}{2}}+\int_{t(s)}^\infty\|r\Gamma^I(S+W)\|_{L^2_x(\Sigma^{\rm{ex}}_\tau)}{\rm{d}}\tau\\
&\lesssim\ez+(C_1\ez)^2,\quad |I|\le N-2,
\end{split}
\eeq
where (see \eqref{s2: E^c,ex,h_dz} and \eqref{s2: E^c,ex,h_0})
\beq\label{s3: exhcU}
\begin{split}
E^{c,{\rm{ex}},h}_{\dz}(s,\Gamma^IU):&=\int_{\mathscr{H}^{\rm{ex}}_s}t^{-\dz}\big(|L_0\Gamma^IU+x^a\bar{\partial}_a\Gamma^IU+2\Gamma^IU|^2+s^2|\bar{\nabla}\Gamma^IU|^2\big){\rm{d}}x,\\
E^{c,{\rm{ex}},h}(s,\Gamma^IU):&=E^{c,{\rm{ex}},h}_{0}(s,\Gamma^IU).
\end{split}
\eeq
We recall that $U-u=F_u$ (see \eqref{s3: defU}). We also note that on $\mathscr{H}^{\rm{ex}}_s$ we have (recall \eqref{s2: barpar_a} and \eqref{L_0}) 
\ben
&&\f{t}{2}\le t-1\le r\le t,\quad\quad s^2\lesssim\f{s^2-1}{2}\le r,\\
&&|L_0w|+|x^a\bar{\partial}_aw|+|w|+s|\bar{\nabla}w|\lesssim \f{s^2}{t}|\partial w|+|Lw|+|w|\lesssim|\partial w|+|Lw|+|w|,\\
&&\f{1}{r}\big(|L_0w|+|x^a\bar{\partial}_aw|+|w|+s|\bar{\nabla}w|\big)\lesssim \f{1}{t}|\partial w|+|\bar{\nabla}w|+\f{1}{t}|w|
\een
for any smooth function $w$. Hence,
\beqn\label{s3: exchudv}
&&[E^{c,{\rm{ex}},h}(s,\Gamma^I(u\partial v))]^{\f{1}{2}}\nonumber\\
&\lesssim&II:=\sum_{|I_1|+|I_2|\le N-1}\Bigg\{\l(\int_{\mathscr{H}^{\rm{ex}}_s}\big\{|(L_0+x^a\bar{\partial}_a+2)\Gamma^{I_1}u|^2+s^2|\bar{\nabla}\Gamma^{I_1}u|^2\big\}|\Gamma^{I_2}\partial v|^2{\rm{d}}x\r)^{\f{1}{2}}\nonumber\\
&&\quad\quad\quad\quad\quad\quad\quad+\l(\int_{\mathscr{H}^{\rm{ex}}_s}|\Gamma^{I_1}u|^2\big\{|(L_0+x^a\bar{\partial}_a+2)\Gamma^{I_2}\partial v|^2+s^2|\bar{\nabla}\Gamma^{I_2}\partial v|^2\big\}{\rm{d}}x\r)^{\f{1}{2}}\Bigg\},
\eeqn
where
\beqn\label{s3: exchudv'}
II&\lesssim&\sum_{\substack{|I_1|\le N-4\\|I_2|\le N-1}}\||\partial\Gamma^{I_1}u|+|L\Gamma^{I_1}u|+|\Gamma^{I_1}u|\|_{L^\infty(\mathscr{H}^{\rm{ex}}_s)}\cdot\l(\int_{\mathscr{H}^{\rm{ex}}_s}|\Gamma^{I_2}\partial v|^2{\rm{d}}x\r)^{\f{1}{2}}\nonumber\\
&+&\sum_{\substack{|I_2|\le N-5\\1\le |I_1|\le N-1}}\|t^{\f{\dz}{2}}r\Gamma^{I_2}\partial v\|_{L^\infty(\mathscr{H}^{\rm{ex}}_s)}\cdot\l(\int_{\mathscr{H}^{\rm{ex}}_s}t^{-\dz}\big(t^{-1}|\partial\Gamma^{I_1}u|+|\bar{\nabla}\Gamma^{I_1}u|+t^{-1}|\Gamma^{I_1}u|\big)^2{\rm{d}}x\r)^{\f{1}{2}}\nonumber\\
&+&\sum_{\substack{|I_1|\le N-3\\|I_2|\le N-1}}\|r\Gamma^{I_1}u\|_{L^\infty(\mathscr{H}^{\rm{ex}}_s)}\cdot\l(\int_{\mathscr{H}^{\rm{ex}}_s}\big(t^{-1}|\partial\Gamma^{I_2}\partial v|+|\bar{\nabla}\Gamma^{I_2}\partial v|+t^{-1}|\Gamma^{I_2}\partial v|\big)^2{\rm{d}}x\r)^{\f{1}{2}}\nonumber\\
&+&\sum_{\substack{|I_2|\le N-6\\1\le|I_1|\le N-1}}\|t^{\f{\dz}{2}}t\big(|\partial\Gamma^{I_2}\partial v|+|L\Gamma^{I_2}\partial v|+|\Gamma^{I_2}\partial v|\big)\|_{L^\infty(\mathscr{H}^{\rm{ex}}_s)}\cdot\l(\int_{\mathscr{H}^{\rm{ex}}_s}|t^{-\f{\dz}{2}}t^{-1}\Gamma^{I_1}u|^2{\rm{d}}x\r)^{\f{1}{2}}\nonumber\\
&\lesssim&(C_1\ez)^2,
\eeqn
where we use \eqref{s3: exuvL^infty}, \eqref{s3: exuvL^infty'}, \eqref{s3: exhvN}, \eqref{s3: exvextra''} and \eqref{s3: exhuN}, and that
\ben
\sum_{1\le |I_1|\le N}\l(\int_{\mathscr{H}^{\rm{ex}}_s}|t^{-\f{\dz}{2}}t^{-1}\Gamma^{I_1}u|^2{\rm{d}}x\r)^{\f{1}{2}}\lesssim\sum_{|I'_1|\le N}\l(\int_{\mathscr{H}^{\rm{ex}}_s}t^{-\dz}t^{-2}\big(|\partial\Gamma^{I'_1}u|+|L\Gamma^{I'_1}u|\big)^2{\rm{d}}x\r)^{\f{1}{2}}\lesssim\sum_{|I|\le N}[E^{{\rm{ex}},h}_{0,\dz}(s,\Gamma^Iu)]^{\f{1}{2}}.
\een
Here we note that on $\mathscr{H}_s$ we have
\beq\label{s3: Omab}
|\Omega_{ab}w|=\l|\f{x_a}{t}L_bw-\f{x_b}{t}L_aw\r|\lesssim\sum_{c=1}^3|L_cw|,\quad 1\le a<b\le 3,
\eeq
for any smooth function $w$. The estimates \eqref{s3: exchudv}-\eqref{s3: exchudv'} also hold if we replace $u\partial v$ on the left hand side of \eqref{s3: exchudv} by $uv$ (which is simpler), or by $\partial u\partial v$ and $\partial uv$ (for the estimates of them we do not need \eqref{s3: exuvL^infty'}). Combining these estimates with \eqref{s3: S2.2U}, we obtain that for all $s\in[2,\infty)$
\beq\label{s3: exEch}
\begin{split}
\l(\int_{\mathscr{H}^{\rm{ex}}_s}t^{-\dz}\big(|(L_0+x^a\bar{\partial}_a+2)\Gamma^Iu|^2+s^2|\bar{\nabla}\Gamma^Iu|^2\big){\rm{d}}x\r)^{\f{1}{2}}&\lesssim\ez+(C_1\ez)^2,\quad |I|\le N-1,\\
\l(\int_{\mathscr{H}^{\rm{ex}}_s}\big(|(L_0+x^a\bar{\partial}_a+2)\Gamma^Iu|^2+s^2|\bar{\nabla}\Gamma^Iu|^2\big){\rm{d}}x\r)^{\f{1}{2}}&\lesssim\ez+(C_1\ez)^2,\quad |I|\le N-2.
\end{split}
\eeq
By Proposition \ref{s2: exconh}, we obtain
\beq\label{s3: exhuL^2}
\begin{split}
&\sum_{|I|\le N-1}\Bigg\{\|t^{-\f{\dz}{2}}st^{-1}\Gamma^Iu\|_{L^2(\mathscr{H}^{\rm{ex}}_s)}+\l(\int_{|x|=r(s)}(t(s))^{-\dz}s^2r^{-1}|\Gamma^Iu|^2(t(s), x){\rm{d}}\sz\r)^{\f{1}{2}}\Bigg\}\\
&\lesssim\sum_{|I|\le N-1}\|t^{-\f{\dz}{2}}s|\bar{\nabla}\Gamma^Iu|\|_{L^2(\mathscr{H}^{\rm{ex}}_s)}
\lesssim\ez+(C_1\ez)^2,\\
&\sum_{|I|\le N-2}\Bigg\{\|st^{-1}\Gamma^Iu\|_{L^2(\mathscr{H}^{\rm{ex}}_s)}+\l(\int_{|x|=r(s)}s^2r^{-1}|\Gamma^Iu|^2(t(s), x){\rm{d}}\sz\r)^{\f{1}{2}}\Bigg\}\\
&\lesssim\sum_{|I|\le N-2}\|s|\bar{\nabla}\Gamma^Iu|\|_{L^2(\mathscr{H}^{\rm{ex}}_s)}\lesssim\ez+(C_1\ez)^2,
\end{split}
\eeq
where $t(s)=\f{s^2+1}{2}$ and $r(s)=\f{s^2-1}{2}$ (see \eqref{s2: tr(s)}). We note that on $\mathscr{H}_s$ we have 
\beq\label{s3: st^-1L_0}
st^{-1}|L_0w|\le|(L_0+x^a\bar{\partial}_a+2)w|+s|\bar{\nabla}w|+st^{-1}|w|
\eeq
for any smooth function $w$. Combining \eqref{s3: exEch}, \eqref{s3: exhuL^2}, \eqref{s3: st^-1L_0} and \eqref{s3: Omab}, we obtain that for all $s\in[2,\infty)$ 
\beqn\label{s3: exhL_0uL^2}
&&\sum_{|I|\le N-1}\|t^{-\f{\dz}{2}}st^{-1}\big(|\Gamma^Iu|+|L_0\Gamma^Iu|+|L\Gamma^Iu|+|\Omega\Gamma^Iu|\big)\|_{L^2(\mathscr{H}^{\rm{ex}}_s)}\nonumber\\
&+&\sum_{|I|\le N-2}\|st^{-1}\big(|\Gamma^Iu|+|L_0\Gamma^Iu|+|L\Gamma^Iu|+|\Omega\Gamma^Iu|\big)\|_{L^2(\mathscr{H}^{\rm{ex}}_s)}\lesssim\ez+(C_1\ez)^2.
\eeqn

\section{Global existence in the interior region}\label{s4}

In this section we prove global existence of solutions to \eqref{s3: equv}-\eqref{s1: ini} in the interior region.

\subsection{Bootstrap setting}

The bootstrap assumptions in the interior region are the following energy bounds for the solution $(u,v)$ to \eqref{s3: equv}-\eqref{s1: ini}:
\beq\label{s4: inuvL^2N}
\sum_{|I|\le N}\big\{[E^{{\rm{in}},h}_0(s,\Gamma^Iu)]^{\f{1}{2}}+[E^{{\rm{in}},h}_1(s,\Gamma^Iv)]^{\f{1}{2}}\big\}\le C_2\ez s^{\dz},
\eeq
\beq\label{s4: inuvL^2N-1}
\sum_{|I|\le N-1}\big\{[E^{{\rm{in}},h}_0(s,\Gamma^Iu)]^{\f{1}{2}}+[E^{{\rm{in}},h}_1(s,\Gamma^Iv)]^{\f{1}{2}}\big\}\le C_2\ez,
\eeq
\beq\label{s4: inucN-1}
\sum_{|I|\le N-1}[E^{c,\rm{in}}(s,\Gamma^Iu)]^{\f{1}{2}}\le C_2\ez s^{\f{1}{2}+\dz},
\eeq
\beq\label{s4: inucN-2}
\sum_{|I|\le N-2}[E^{c,\rm{in}}(s,\Gamma^Iu)]^{\f{1}{2}}\le C_2\ez s^{\dz},\quad\quad\sum_{|I|\le N-3}[E^{c,\rm{in}}(s,\Gamma^Iu)]^{\f{1}{2}}\le C_2\ez,
\eeq
where $N\ge 9$, $0<\dz\ll 1$ are as in Section \ref{ss3.1}, $C_2\gg C_1, 0<\ez\ll C_2^{-1}$ are constants to be chosen below ($C_1$ is as in Section \ref{ss3.1}), and we recall \eqref{dinhs}, \eqref{dhs0}, \eqref{s2: E^c,in_dz} and \eqref{s2: E^c,in_0} for the definitions of the energy functionals above.

Let $\tilde{T}^*$ be defined as
\beq\label{s4: tildeT^*}
\tilde{T}^*:=\sup\{\tilde{T}>2: \eqref{s4: inuvL^2N}-\eqref{s4: inucN-2} \mathrm{\ hold\ for\ } s\in[2,\tilde{T}]\}.
\eeq

\begin{prop}\label{s4: maxtildeT}
There exists some constants $C_2>0$ sufficiently large and $0<\ez_1\ll C_2^{-1}$ sufficiently small such that for any $0<\ez<\ez_1$, if $(u,v)$ is a solution to \eqref{s3: equv}-\eqref{s1: ini} and satisfies the exterior estimates \eqref{s3: bsex}-\eqref{s3: bsex'} globally in time as well as the interior bounds \eqref{s4: inuvL^2N}-\eqref{s4: inucN-2} for all $s\in[2,\tilde{T}]$, then we have the following improved interior estimates
\be\label{s4: iminuvL^2N}
\begin{split}
&\sum_{|I|\le N}\big\{[E^{{\rm{in}},h}_0(s,\Gamma^Iu)]^{\f{1}{2}}+[E^{{\rm{in}},h}_1(s,\Gamma^Iv)]^{\f{1}{2}}\big\}\le \f{1}{2}C_2\ez s^{\dz},\\ 
&\sum_{|I|\le N-1}\big\{[E^{{\rm{in}},h}_0(s,\Gamma^Iu)]^{\f{1}{2}}+[E^{{\rm{in}},h}_1(s,\Gamma^Iv)]^{\f{1}{2}}\big\}\le \f{1}{2}C_2\ez,\\
&\sum_{|I|\le N-1}[E^{c,\rm{in}}(s,\Gamma^Iu)]^{\f{1}{2}}\le \f{1}{2}C_2\ez s^{\f{1}{2}+\dz},\\
&\sum_{|I|\le N-2}[E^{c,\rm{in}}(s,\Gamma^Iu)]^{\f{1}{2}}\le \f{1}{2}C_2\ez s^{\dz},\quad\quad\sum_{|I|\le N-3}[E^{c,\rm{in}}(s,\Gamma^Iu)]^{\f{1}{2}}\le \f{1}{2}C_2\ez
\end{split}
\ee
for all $s\in[2,\tilde{T}]$.
\end{prop}

In the above proposition the hyperbolic time $\tilde{T}$ is arbitrary, therefore the solution $(u,v)$ exists globally in the interior region $\mathscr{D}^{\rm{in}}$ (i.e., $\tilde{T}^*=\infty$ where $\tilde{T}^*$ is as in \eqref{s4: tildeT^*}) and satisfies the estimates \eqref{s4: inuvL^2N}-\eqref{s4: inucN-2} for all $s\in[2,\infty)$. Below we give the proof of Proposition \ref{s4: maxtildeT}.

Using Proposition \ref{s2: incon}, \eqref{s3: exhuL^2} and \eqref{s4: inuvL^2N}-\eqref{s4: inucN-2}, we obtain
\beq\label{s4: inhuN-1}
\begin{split}
\|st^{-1}\Gamma^Iu\|_{L^2(\mathscr{H}^{\rm{in}}_s)}&\lesssim\l(\int_{|x|=r(s)}s^2r^{-1}|\Gamma^Iu|^2(t(s),x){\rm{d}}\sz\r)^{\f{1}{2}}+[E^{c,\rm{in}}(s,\Gamma^Iu)]^{\f{1}{2}}\\
&\lesssim[\ez+(C_1\ez)^2]s^\dz+C_2\ez s^{\f{1}{2}+\dz}\lesssim C_2\ez s^{\f{1}{2}+\dz},\quad |I|\le N-1,\\
\|st^{-1}\Gamma^Iu\|_{L^2(\mathscr{H}^{\rm{in}}_s)}&\lesssim\l(\int_{|x|=r(s)}s^2r^{-1}|\Gamma^Iu|^2(t(s),x){\rm{d}}\sz\r)^{\f{1}{2}}+[E^{c,\rm{in}}(s,\Gamma^Iu)]^{\f{1}{2}}\\
&\lesssim\ez+(C_1\ez)^2+C_2\ez s^{\dz}\lesssim C_2\ez s^{\dz},\quad|I|\le N-2,\\
\|st^{-1}\Gamma^Iu\|_{L^2(\mathscr{H}^{\rm{in}}_s)}&\lesssim C_2\ez,\quad|I|\le N-3.
\end{split}
\eeq
We recall \eqref{e^h_m}, \eqref{dinhs}, \eqref{dhs0}, \eqref{s2: E^c,in_dz} and \eqref{s2: E^c,in_0} for the definitions of the energy functionals appearing in \eqref{s4: inuvL^2N}-\eqref{s4: inucN-2}. Using \eqref{s3: exhvN}, \eqref{s3: exhuN-1}, \eqref{s3: exhuL^2}, \eqref{s3: exhL_0uL^2} and \eqref{s4: inhuN-1} (and recalling \eqref{s3: st^-1L_0} and \eqref{s3: Omab}), we obtain the following $L^2$-type estimates for $s\in[2,\tilde{T})$:
\beq\label{s4: inL^2N}
\sum_{|I|\le N}\|st^{-1}\big(|\partial\Gamma^Iu|+|\partial\Gamma^Iv|\big)+|\Gamma^Iv|\|_{L^2(\mathscr{H}^{\rm{in}}_s)}\lesssim C_2\ez s^\dz,
\eeq
\beq\label{s4: inL^2N-1}
\sum_{|I|\le N-1}\|st^{-1}\big(|\partial\Gamma^Iu|+|\partial\Gamma^Iv|\big)+|\Gamma^Iv|\|_{L^2(\mathscr{H}_s)}\lesssim C_2\ez,
\eeq
\beq\label{s4: incL^2N-1}
\sum_{|I|\le N-1}\|st^{-1}\big(|\Gamma^Iu|+|L_0\Gamma^Iu|+|L\Gamma^Iu|+|\Omega\Gamma^Iu|\big)\|_{L^2(\mathscr{H}^{\rm{in}}_s)}\lesssim C_2\ez s^{\f{1}{2}+\dz},
\eeq
\beq\label{s4: incL^2N-2}
\sum_{|I|\le N-2}\|st^{-1}\big(|\Gamma^Iu|+|L_0\Gamma^Iu|+|L\Gamma^Iu|+|\Omega\Gamma^Iu|\big)\|_{L^2(\mathscr{H}_s)}\lesssim C_2\ez s^{\dz},
\eeq
\beq\label{s4: incL^2N-3}
\sum_{|I|\le N-3}\|st^{-1}\big(|\Gamma^Iu|+|L_0\Gamma^Iu|+|L\Gamma^Iu|+|\Omega\Gamma^Iu|\big)\|_{L^2(\mathscr{H}_s)}\lesssim C_2\ez.
\eeq
We note that on $\mathscr{H}_s$ we have $r\le t$ and
\ben
&&L_a\l(\f{s}{t}\r)=t\bar{\partial}_a\l(\f{s}{t}\r)=-\f{s}{t}\f{x_a}{t},\\
&&L_bL_a\l(\f{s}{t}\r)=-t\bar{\partial}_b\l(\f{s}{t}\f{x_a}{t}\r)=\l(2\f{x_ax_b}{t^2}-\dz_{ab}\r)\f{s}{t}.
\een
Hence by \eqref{s4: inL^2N-1}, \eqref{s4: incL^2N-3} and Lemma \ref{s2: Sobin}, we obtain the following pointwise estimates for $s\in[2,\tilde{T})$:
\beqn\label{s4: inL^infty}
&&\sum_{|I|\le N-3}\|t^{\f{3}{2}}\big(st^{-1}|\partial\Gamma^Iu|+st^{-1}|\partial\Gamma^Iv|+|\Gamma^Iv|\big)\|_{L^\infty(\mathscr{H}^{\rm{in}}_s)}\nonumber\\
&+&\sum_{|I|\le N-5}\|st^{\f{1}{2}}\big(|\Gamma^Iu|+|L_0\Gamma^Iu|+|L\Gamma^Iu|\big)\|_{L^\infty(\mathscr{H}^{\rm{in}}_s)}\lesssim C_2\ez.
\eeqn

\subsection{Improved estimates for the solution $(u,v)$ in the interior region}\label{ss4.2}

${\bf Step\ 1.}$ Refining the highest order energy estimates of $u$ and $v$: By Proposition \ref{s2: inhs}, for $|I|\le N$, we have
\beq\label{s4.2: S1uv}
\begin{split}
E^{{\rm{in}},h}_0(s,\Gamma^Iu)&\lesssim E^{{\rm{in}},h}_0(2,\Gamma^Iu)+\int_{\mathscr{l}_{[\f{5}{2},\f{s^2+1}{2}]}}|G\Gamma^Iu|^2{\rm{d}}\sz+\int_2^s[E^{{\rm{in}},h}_0(\tau,\Gamma^Iu)]^{\f{1}{2}}\cdot\|\Gamma^IF_u\|_{L^2(\mathscr{H}^{\rm{in}}_\tau)}{\rm{d}}\tau,\\
E^{{\rm{in}},h}_1(s,\Gamma^Iv)&\lesssim E^{{\rm{in}},h}_1(2,\Gamma^Iv)+\int_{\mathscr{l}_{[\f{5}{2},\f{s^2+1}{2}]}}\big(|G\Gamma^Iv|^2+|\Gamma^Iv|^2\big){\rm{d}}\sz+\int_2^s[E^{{\rm{in}},h}_1(\tau,\Gamma^Iv)]^{\f{1}{2}}\cdot\|\Gamma^IF_v\|_{L^2(\mathscr{H}^{\rm{in}}_\tau)}{\rm{d}}\tau.
\end{split}
\eeq
We recall \eqref{s3: equv} for the definitions of $F_u$ and $F_v$. We note the fact that $\tau\le t\lesssim\tau^2$ on $\mathscr{H}^{\rm{in}}_\tau$. For $|I|\le N$ and $\tau\in[2,s]$, we have 
\beqn\label{s4: improN}
&&\||\Gamma^I(uv)|+|\Gamma^I(\partial uv)|+|\Gamma^I(u\partial v)|+|\Gamma^I(\partial u\partial v)|\|_{L^2(\mathscr{H}^{\rm{in}}_\tau)}\nonumber\\
&\lesssim&\sum_{|I_1|+|I_2|\le N}\|\big(|\Gamma^{I_1}u|+|\partial\Gamma^{I_1}u|\big)\cdot\big(|\Gamma^{I_2}v|+|\partial\Gamma^{I_2}v|\big)\|_{L^2(\mathscr{H}^{\rm{in}}_\tau)}\nonumber\\
&\lesssim&\sum_{\substack{|I_1|\le N-5\\|I_2|\le N}}\|\tau t^{\f{1}{2}}\big(|\Gamma^{I_1}u|+|\partial\Gamma^{I_1}u|\big)\|_{L^\infty(\mathscr{H}^{\rm{in}}_\tau)}\cdot\tau^{-\dz}\l\|\f{\tau}{t}\big(|\Gamma^{I_2}v|+|\partial\Gamma^{I_2}v|\big)\r\|_{L^2(\mathscr{H}^{\rm{in}}_\tau)}\cdot\tau^{-1+\dz}\nonumber\\
&+&\!\!\!\!\sum_{\substack{|I_2|\le N-4\\|J|=1,|I'_1|\le N-1\\|I_1|\le N}}\!\!\Big\{\tau^{-\f{1}{2}-\dz}\l\|\f{\tau}{t}\Gamma^J\Gamma^{I'_1}u\r\|_{L^2(\mathscr{H}^{\rm{in}}_\tau)}+\tau^{-\dz}\l\|\f{\tau}{t}\partial\Gamma^{I_1}u\r\|_{L^2(\mathscr{H}^{\rm{in}}_\tau)}\Big\}\cdot\|t^{\f{3}{2}}\big(|\Gamma^{I_2}v|+|\partial\Gamma^{I_2}v|\big)\|_{L^\infty(\mathscr{H}^{\rm{in}}_\tau)}\cdot\tau^{-1+\dz}\nonumber\\
&\lesssim&(C_2\ez)^2\tau^{-1+\dz},
\eeqn
where we use \eqref{s4: inL^infty}, \eqref{s4: inL^2N} and \eqref{s4: incL^2N-1}. For $|I|\le N-1$ and $\tau\in[2,s]$, we have
\beqn\label{s4: improN-1}
&&\||\Gamma^I(uv)|+|\Gamma^I(\partial uv)|+|\Gamma^I(u\partial v)|+|\Gamma^I(\partial u\partial v)|\|_{L^2(\mathscr{H}^{\rm{in}}_\tau)}\nonumber\\
&\lesssim&\sum_{\substack{|I_1|\le N-5\\|I_2|\le N}}\|\tau t^{\f{1}{2}}\big(|\Gamma^{I_1}u|+|\partial\Gamma^{I_1}u|\big)\|_{L^\infty(\mathscr{H}^{\rm{in}}_\tau)}\cdot\tau^{-\dz}\|\Gamma^{I_2}v\|_{L^2(\mathscr{H}^{\rm{in}}_\tau)}\cdot\tau^{-\f{3}{2}+\dz}\nonumber\\
&+&\sum_{\substack{|I_2|\le N-4\\|J|=1,|I'_1|\le N-2\\|I_1|\le N-1}}\Big\{\tau^{-\dz}\l\|\f{\tau}{t}\Gamma^J\Gamma^{I'_1}u\r\|_{L^2(\mathscr{H}^{\rm{in}}_\tau)}+\l\|\f{\tau}{t}\partial\Gamma^{I_1}u\r\|_{L^2(\mathscr{H}^{\rm{in}}_\tau)}\Big\}\cdot\|t^{\f{3}{2}}\big(|\Gamma^{I_2}v|+|\partial\Gamma^{I_2}v|\big)\|_{L^\infty(\mathscr{H}^{\rm{in}}_\tau)}\cdot\tau^{-\f{3}{2}+\dz}\nonumber\\
&\lesssim&(C_2\ez)^2\tau^{-\f{3}{2}+\dz},
\eeqn
where we use \eqref{s4: inL^infty}, \eqref{s4: inL^2N}, \eqref{s4: incL^2N-2} and \eqref{s4: inL^2N-1}.

We point out that the estimates of the boundary integrals on $\mathscr{l}_{[\f{5}{2},\f{s^2+1}{2}]}$ in \eqref{s4.2: S1uv} were given in \eqref{s3: S1gzuN} (see the definition \eqref{s3: gzudz}, \eqref{s3: gzv0} and \eqref{s3: gzu0}); see also \eqref{s3: exhuN}, \eqref{s3: exhvN} and \eqref{s3: exhuN-1}. We also note that on the boundary $\mathscr{l}_{[\f{5}{2},\f{s^2+1}{2}]}$ we have $t\lesssim s^2$ and therefore $s^{-2\dz}\lesssim t^{-\dz}$. Hence, we obtain
\be
\begin{split}
\sum_{|I|\le N}E^{{\rm{in}},h}_0(s,\Gamma^Iu)&\lesssim [\ez^2+(C_1\ez)^3]s^{2\dz}+(C_2\ez)^2\sum_{|I|\le N}\int_2^s[E^{{\rm{in}},h}_0(\tau,\Gamma^Iu)]^{\f{1}{2}}\cdot\tau^{-1+\dz}{\rm{d}}\tau\\
&\lesssim [\ez^2+(C_1\ez)^3]s^{2\dz}+(C_2\ez)^2 s^\dz\sup_{\tau\in[2,s]}\sum_{|I|\le N}[E^{{\rm{in}},h}_0(\tau,\Gamma^Iu)]^{\f{1}{2}},\\
\sum_{|I|\le N}E^{{\rm{in}},h}_1(s,\Gamma^Iv)&\lesssim\ez^2+(C_1\ez)^3+(C_2\ez)^2\sum_{|I|\le N}\int_2^s[E^{{\rm{in}},h}_1(\tau,\Gamma^Iv)]^{\f{1}{2}}\cdot\tau^{-1+\dz}{\rm{d}}\tau\\
&\lesssim\ez^2+(C_1\ez)^3+(C_2\ez)^2 s^\dz\sup_{\tau\in[2,s]}\sum_{|I|\le N}[E^{{\rm{in}},h}_1(\tau,\Gamma^Iv)]^{\f{1}{2}}
\end{split}
\ee
and
\be
\begin{split}
&\quad\sum_{|I|\le N-1}\big\{E^{{\rm{in}},h}_0(s,\Gamma^Iu)+E^{{\rm{in}},h}_1(s,\Gamma^Iv)\big\}\\
&\lesssim\ez^2+(C_1\ez)^3+(C_2\ez)^2\sum_{|I|\le N-1}\int_2^s\big\{[E^{{\rm{in}},h}_0(\tau,\Gamma^Iu)]^{\f{1}{2}}+[E^{{\rm{in}},h}_1(\tau,\Gamma^Iv)]^{\f{1}{2}}\big\}\cdot\tau^{-\f{3}{2}+\dz}{\rm{d}}\tau\\
&\lesssim\ez^2+(C_1\ez)^3+(C_2\ez)^2\sup_{\tau\in[2,s]}\sum_{|I|\le N-1}\big\{[E^{{\rm{in}},h}_0(\tau,\Gamma^Iu)]^{\f{1}{2}}+[E^{{\rm{in}},h}_1(\tau,\Gamma^Iv)]^{\f{1}{2}}\big\}.
\end{split}
\ee
It follows that
\beq\label{s4: S1us}
\sum_{|I|\le N}\big\{[E^{{\rm{in}},h}_0(s,\Gamma^Iu)]^{\f{1}{2}}+[E^{{\rm{in}},h}_1(s,\Gamma^Iv)]^{\f{1}{2}}\big\}\lesssim[\ez+(C_2\ez)^{\f{3}{2}}]s^\dz,
\eeq
\beq\label{s4: S2vs}
\sum_{|I|\le N-1}\big\{[E^{{\rm{in}},h}_0(s,\Gamma^Iu)]^{\f{1}{2}}+[E^{{\rm{in}},h}_1(s,\Gamma^Iv)]^{\f{1}{2}}\big\}\lesssim\ez+(C_2\ez)^{\f{3}{2}}.
\eeq

${\bf Step\ 2.}$ Refining the conformal energy estimates of $u$: Let $U=u+F_u$ be as in \eqref{s3: defU}. By Proposition \ref{s2: incon}, for $|I|\le N-1$, we have
\beq\label{s4: E_cinU0}
E^{c,\rm{in}}(s,\Gamma^IU)\lesssim E^{c,\rm{in}}(2,\Gamma^IU)+\mathcal{\tilde{R}}(s,\Gamma^IU)+\int_2^s[E^{c,\rm{in}}(\tau,\Gamma^IU)]^{\f{1}{2}}\cdot\|\tau\Gamma^I(S+W)\|_{L^2(\mathscr{H}^{\rm{in}}_\tau)}{\rm{d}}\tau,
\eeq
where $S$ and $W$ are as in \eqref{s3: defS} and \eqref{s3: defT}, and
\beq\label{s4: tildeRU}
\begin{split}
&\mathcal{\tilde{R}}(s,\Gamma^IU):=\int_{\mathscr{l}_{[\f{5}{2},\f{s^2+1}{2}]}}H[\Gamma^IU](t,x){\rm{d}}\sz,\\
&H[\Gamma^IU]:=|(t+r)(\partial_t\Gamma^IU+\partial_r\Gamma^IU)+2\Gamma^IU|^2+(t-r)^2r^{-2}|\Omega \Gamma^IU|^2.
\end{split}
\eeq

For $|I|\le N-1$ and $\tau\in[2,s]$, by Lemma \ref{s2: Q_0}, we have
\beqn\label{s4: S3gzudv}
&&\||\Gamma^I(\partial_\gz u\partial^\gz v)|+|\Gamma^I(\partial_\gz u\partial^\gz \partial v)|\|_{L^2(\mathscr{H}^{\rm{in}}_\tau)}\nonumber\\
&\lesssim&\sum_{|I_1|+|I_2|\le N-1}\||(\partial_\gz\Gamma^{I_1}u)(\partial^\gz\Gamma^{I_2}v)|+|(\partial_\gz\Gamma^{I_1}u)(\partial^\gz\Gamma^{I_2}\partial v)|\|_{L^2(\mathscr{H}^{\rm{in}}_\tau)}\nonumber\\
&\lesssim&\sum_{|I_1|+|I_2|\le N-1}\l\|t^{-1}\big\{\big(|\partial_t\Gamma^{I_2}v|+|\partial_t\Gamma^{I_2}\partial v|\big)\cdot|L_0\Gamma^{I_1}u|+\big(|L_a\Gamma^{I_2}v|+|L_a\Gamma^{I_2}\partial v|\big)\cdot|\partial^a\Gamma^{I_1}u|\big\}\r\|_{L^2(\mathscr{H}^{\rm{in}}_\tau)}\nonumber\\
&\lesssim&\sum_{\substack{|I_1|\le N-5\\|I_2|\le N}}\|\tau t^{\f{1}{2}}L_0\Gamma^{I_1}u\|_{L^\infty(\mathscr{H}^{\rm{in}}_\tau)}\cdot\tau^{-\dz}\l\|\f{\tau}{t}\partial\Gamma^{I_2}v\r\|_{L^2(\mathscr{H}^{\rm{in}}_\tau)}\cdot\tau^{-\f{5}{2}+\dz}\nonumber\\
&+&\sum_{\substack{|I_2|\le N-5\\|I_1|\le N-1}}\tau^{-\f{1}{2}-\dz}\l\|\f{\tau}{t}L_0\Gamma^{I_1}u\r\|_{L^2(\mathscr{H}^{\rm{in}}_\tau)}\cdot\|t^{\f{3}{2}}\big(|\partial\Gamma^{I_2}v|+|\partial\Gamma^{I_2}\partial v|\big)\|_{L^\infty(\mathscr{H}^{\rm{in}}_\tau)}\cdot\tau^{-2+\dz}\nonumber\\
&+&\sum_{\substack{|I_1|\le N-3\\|I_2|\le N}}\|\tau t^{\f{1}{2}}\partial\Gamma^{I_1}u\|_{L^\infty(\mathscr{H}^{\rm{in}}_\tau)}\cdot\tau^{-\dz}\l\|\f{\tau}{t}\big(|\Gamma^{I_2}v|+|\partial\Gamma^{I_2}v|\big)\r\|_{L^2(\mathscr{H}^{\rm{in}}_\tau)}\cdot\tau^{-\f{5}{2}+\dz}\nonumber\\
&+&\sum_{\substack{|I_2|\le N-5\\|I_1|\le N-1}}\tau^{-\dz}\l\|\f{\tau}{t}\partial\Gamma^{I_1}u\r\|_{L^2(\mathscr{H}^{\rm{in}}_\tau)}\cdot\|t^{\f{3}{2}}\big(|L\Gamma^{I_2}v|+|L\Gamma^{I_2}\partial v|\big)\|_{L^\infty(\mathscr{H}^{\rm{in}}_\tau)}\cdot\tau^{-\f{5}{2}+\dz}\nonumber\\
&\lesssim&(C_2\ez)^2\tau^{-2+\dz},
\eeqn
where we use \eqref{s4: inL^infty}, \eqref{s4: inL^2N} and \eqref{s4: incL^2N-1}. For $|I|\le N-2$, we use \eqref{s4: incL^2N-2} (instead of \eqref{s4: incL^2N-1}) in \eqref{s4: S3gzudv} and obtain
\beq\label{s4: S3gzudvN-2}
\||\Gamma^I(\partial_\gz u\partial^\gz v)|+|\Gamma^I(\partial_\gz u\partial^\gz \partial v)|\|_{L^2(\mathscr{H}^{\rm{in}}_\tau)}\lesssim(C_2\ez)^2\tau^{-\f{5}{2}+\dz},\quad |I|\le N-2.
\eeq
 We also note that, the proof of \eqref{s4: S3gzudv} gives
\beq\label{s4: S3gzdudvN-2}
\||\Gamma^I(\partial_\gz \partial u\partial^\gz v)|+|\Gamma^I(\partial_\gz \partial u\partial^\gz \partial v)|\|_{L^2(\mathscr{H}^{\rm{in}}_\tau)}\lesssim(C_2\ez)^2\tau^{-2+\dz},\quad |I|\le N-2.
\eeq
Similarly, the proof of \eqref{s4: S3gzudvN-2} implies
\beq\label{s4: S3gzdudvN-3}
\||\Gamma^I(\partial_\gz \partial u\partial^\gz v)|+|\Gamma^I(\partial_\gz \partial u\partial^\gz \partial v)|\|_{L^2(\mathscr{H}^{\rm{in}}_\tau)}\lesssim(C_2\ez)^2\tau^{-\f{5}{2}+\dz},\quad |I|\le N-3.
\eeq
We note that $\f{t-r}{t}\sim\f{\tau^2}{t^2}$ and $t\lesssim\tau^2$ on $\mathscr{H}^{\rm{in}}_\tau$. For $|I|\le N-1$, using Lemma \ref{s2: Q_0}, we also have 
\beqn\label{s4: S3gzdudvN-1}
&&\||\Gamma^I(\partial_\gz \partial u\partial^\gz v)|+|\Gamma^I(\partial_\gz \partial u\partial^\gz \partial v)|\|_{L^2(\mathscr{H}^{\rm{in}}_\tau)}\nonumber\\
&\lesssim&\sum_{|I_1|+|I_2|\le N-1}\||(\partial_\gz\Gamma^{I_1}\partial u)(\partial^\gz\Gamma^{I_2}v)|+|(\partial_\gz\Gamma^{I_1}\partial u)(\partial^\gz\Gamma^{I_2}\partial v)|\|_{L^2(\mathscr{H}^{\rm{in}}_\tau)}\nonumber\\
&\lesssim&\sum_{\substack{|I_1|+|I_2|\le N-1\\|J|,|J'|=1}}\l\|\f{\tau^2}{t^2}|\Gamma^J\Gamma^{I_1}\partial u|\cdot\big(|\Gamma^{J'}\Gamma^{I_2}v|+|\Gamma^{J'}\Gamma^{I_2}\partial v|\big)\r\|_{L^2(\mathscr{H}^{\rm{in}}_\tau)}\nonumber\\
&\lesssim&\sum_{\substack{|I_1|\le N-4, |J|=1\\|I_2|\le N}}\|\tau t^{\f{1}{2}}\Gamma^J\Gamma^{I_1}\partial u\|_{L^\infty(\mathscr{H}^{\rm{in}}_\tau)}\cdot\tau^{-\dz}\l\|\f{\tau}{t}\big(|\Gamma^{I_2}v|+|\partial\Gamma^{I_2}v|\big)\r\|_{L^2(\mathscr{H}^{\rm{in}}_\tau)}\cdot\tau^{-\f{3}{2}+\dz}\nonumber\\
&+&\sum_{\substack{|I_2|\le N-5, |J'|=1\\|I_1|\le N}}\tau^{-\dz}\l\|\f{\tau}{t}\partial\Gamma^{I_1}u\r\|_{L^2(\mathscr{H}^{\rm{in}}_\tau)}\cdot\|t^{\f{3}{2}}\big(|\Gamma^{J'}\Gamma^{I_2}v|+|\Gamma^{J'}\Gamma^{I_2}\partial v|\big)\|_{L^\infty(\mathscr{H}^{\rm{in}}_\tau)}\cdot\tau^{-\f{3}{2}+\dz}\nonumber\\
&\lesssim&(C_2\ez)^2\tau^{-\f{3}{2}+\dz},
\eeqn
where we use \eqref{s4: inL^infty} and \eqref{s4: inL^2N}. Combining \eqref{s4: S3gzdudvN-1} with \eqref{s4: S3gzudv} (here we recall \eqref{s3: defS}), we obtain
\beq\label{s4: S3gzSN-1}
\|\Gamma^IS\|_{L^2(\mathscr{H}^{\rm{in}}_\tau)}\lesssim(C_2\ez)^2\tau^{-\f{3}{2}+\dz},\quad |I|\le N-1,\tau\in[2,s],
\eeq
while \eqref{s4: S3gzudvN-2} and \eqref{s4: S3gzdudvN-2} give
\beq\label{s4: S3gzSN-2}
\|\Gamma^IS\|_{L^2(\mathscr{H}^{\rm{in}}_\tau)}\lesssim(C_2\ez)^2\tau^{-2+\dz},\quad |I|\le N-2,\tau\in[2,s].
\eeq
Also, the estimates \eqref{s4: S3gzudvN-2} and \eqref{s4: S3gzdudvN-3} imply
\beq\label{s4: S3gzSN-3}
\|\Gamma^IS\|_{L^2(\mathscr{H}^{\rm{in}}_\tau)}\lesssim(C_2\ez)^2\tau^{-\f{5}{2}+\dz},\quad |I|\le N-3,\tau\in[2,s].
\eeq

Next we turn to the estimate of $\|\Gamma^IW\|_{L^2(\mathscr{H}^{\rm{in}}_\tau)}$ for $|I|\le N-1$ and $\tau\in[2,s]$. As in ${\bf Step\ 2.1.}$ in Section \ref{ss3.3}, we only consider the terms $u\partial^2v\partial v$, $\partial^2u\partial v\partial v$, $u^2\partial v$ and $u\partial^2uv$ (see \eqref{s3: pic}) in the expression of $W$ (see \eqref{s3: defT}). The estimates of the other terms are similar to or simpler than these ones. For $|I|\le N-1$, we have
\beqn\label{s4: S3gzT}
&&\|\Gamma^I(u\partial^2 v\partial v)\|_{L^2(\mathscr{H}^{\rm{in}}_\tau)}\lesssim\sum_{|I_1|+|I_2|+|I_3|\le N-1}\||\Gamma^{I_1}u|\cdot|\partial\Gamma^{I_2}\partial v|\cdot|\partial\Gamma^{I_3}v|\|_{L^2(\mathscr{H}^{\rm{in}}_\tau)}\nonumber\\
&\lesssim&\sum_{\substack{|I_1|,|I_2|\le N-5\\|I_3|\le N-1}}\|\tau t^{\f{1}{2}}\Gamma^{I_1}u\|_{L^\infty(\mathscr{H}^{\rm{in}}_\tau)}\cdot\|t^{\f{3}{2}}\partial\Gamma^{I_2}\partial v\|_{L^\infty(\mathscr{H}^{\rm{in}}_\tau)}\cdot\l\|\f{\tau}{t}\partial\Gamma^{I_3}v\r\|_{L^2(\mathscr{H}^{\rm{in}}_\tau)}\cdot\tau^{-3}\nonumber\\
&+&\sum_{\substack{|I_1|,|I_3|\le N-5\\|I_2|\le N-1}}\|\tau t^{\f{1}{2}}\Gamma^{I_1}u\|_{L^\infty(\mathscr{H}^{\rm{in}}_\tau)}\cdot\tau^{-\dz}\l\|\f{\tau}{t}\partial\Gamma^{I_2}\partial v\r\|_{L^2(\mathscr{H}^{\rm{in}}_\tau)}\cdot\|t^{\f{3}{2}}\partial\Gamma^{I_3}v\|_{L^\infty(\mathscr{H}^{\rm{in}}_\tau)}\cdot\tau^{-3+\dz}\nonumber\\
&+&\sum_{\substack{|I_2|,|I_3|\le N-5\\|I_1|\le N-1}}\tau^{-\f{1}{2}-\dz}\l\|\f{\tau}{t}\Gamma^{I_1}u\r\|_{L^2(\mathscr{H}^{\rm{in}}_\tau)}\cdot\|t^{\f{3}{2}}\partial\Gamma^{I_2}\partial v\|_{L^\infty(\mathscr{H}^{\rm{in}}_\tau)}\cdot\|t^{\f{3}{2}}\partial\Gamma^{I_3}v\|_{L^\infty(\mathscr{H}^{\rm{in}}_\tau)}\cdot\tau^{-\f{5}{2}+\dz}\nonumber\\
&\lesssim&(C_2\ez)^3\tau^{-\f{5}{2}+\dz},
\eeqn
where we use \eqref{s4: inL^infty}, \eqref{s4: inL^2N-1}, \eqref{s4: inL^2N} and \eqref{s4: incL^2N-1}. Similarly,
\beqn\label{s4: S3gzT'}
&&\|\Gamma^I(\partial^2u\partial v\partial v)\|_{L^2(\mathscr{H}^{\rm{in}}_\tau)}\lesssim\sum_{|I_1|+|I_2|+|I_3|\le N-1}\||\partial\Gamma^{I_1}\partial u|\cdot|\partial\Gamma^{I_2}v|\cdot|\partial\Gamma^{I_3}v|\|_{L^2(\mathscr{H}^{\rm{in}}_\tau)}\nonumber\\
&\lesssim&\sum_{\substack{|I_1|,|I_2|\le N-4\\|I_3|\le N-1}}\|\tau t^{\f{1}{2}}\partial\Gamma^{I_1}\partial u\|_{L^\infty(\mathscr{H}^{\rm{in}}_\tau)}\cdot\|t^{\f{3}{2}}\partial\Gamma^{I_2}v\|_{L^\infty(\mathscr{H}^{\rm{in}}_\tau)}\cdot\l\|\f{\tau}{t}\partial\Gamma^{I_3}v\r\|_{L^2(\mathscr{H}^{\rm{in}}_\tau)}\cdot\tau^{-3}\nonumber\\
&+&\sum_{\substack{|I_2|,|I_3|\le N-4\\|I_1|\le N-1}}\tau^{-\dz}\l\|\f{\tau}{t}\partial\Gamma^{I_1}\partial u\r\|_{L^2(\mathscr{H}^{\rm{in}}_\tau)}\cdot\|t^{\f{3}{2}}\partial\Gamma^{I_2}v\|_{L^\infty(\mathscr{H}^{\rm{in}}_\tau)}\cdot\|t^{\f{3}{2}}\partial\Gamma^{I_3}v\|_{L^\infty(\mathscr{H}^{\rm{in}}_\tau)}\cdot\tau^{-3+\dz}\nonumber\\
&\lesssim&(C_2\ez)^3\tau^{-3+\dz},
\eeqn
where we use \eqref{s4: inL^infty}, \eqref{s4: inL^2N-1} and \eqref{s4: inL^2N}. We also have
\beqn\label{s4: S3gzT''}
&&\|\Gamma^I(u^2\partial v)\|_{L^2(\mathscr{H}^{\rm{in}}_\tau)}\lesssim\sum_{|I_1|+|I_2|+|I_3|\le N-1}\||\Gamma^{I_1}u|\cdot|\Gamma^{I_2}u|\cdot|\partial\Gamma^{I_3}v|\|_{L^2(\mathscr{H}^{\rm{in}}_\tau)}\nonumber\\
&\lesssim&\sum_{\substack{|I_1|,|I_2|\le N-5\\|I_3|\le N-1}}\|\tau t^{\f{1}{2}}\Gamma^{I_1}u\|_{L^\infty(\mathscr{H}^{\rm{in}}_\tau)}\cdot\|\tau t^{\f{1}{2}}\Gamma^{I_2}u\|_{L^\infty(\mathscr{H}^{\rm{in}}_\tau)}\cdot\l\|\f{\tau}{t}\partial\Gamma^{I_3}v\r\|_{L^2(\mathscr{H}^{\rm{in}}_\tau)}\cdot\tau^{-3}\nonumber\\
&+&\sum_{\substack{|I_2|,|I_3|\le N-5\\|I_1|\le N-1}}\tau^{-\f{1}{2}-\dz}\l\|\f{\tau}{t}\Gamma^{I_1}u\r\|_{L^2(\mathscr{H}^{\rm{in}}_\tau)}\cdot\|\tau t^{\f{1}{2}}\Gamma^{I_2}u\|_{L^\infty(\mathscr{H}^{\rm{in}}_\tau)}\cdot\|t^{\f{3}{2}}\partial\Gamma^{I_3}v\|_{L^\infty(\mathscr{H}^{\rm{in}}_\tau)}\cdot\tau^{-\f{5}{2}+\dz}\nonumber\\
&\lesssim&(C_2\ez)^3\tau^{-\f{5}{2}+\dz}
\eeqn
and
\beqn\label{s4: S3gzT'''}
&&\|\Gamma^I(u\partial^2uv)\|_{L^2(\mathscr{H}^{\rm{in}}_\tau)}\lesssim\sum_{|I_1|+|I_2|+|I_3|\le N-1}\||\Gamma^{I_1}u|\cdot|\partial\Gamma^{I_2}\partial u|\cdot|\Gamma^{I_3}v|\|_{L^2(\mathscr{H}^{\rm{in}}_\tau)}\nonumber\\
&\lesssim&\sum_{\substack{|I_1|,|I_2|\le N-5\\|I_3|\le N-1}}\|\tau t^{\f{1}{2}}\Gamma^{I_1}u\|_{L^\infty(\mathscr{H}^{\rm{in}}_\tau)}\cdot\|\tau t^{\f{1}{2}}\partial\Gamma^{I_2}\partial u\|_{L^\infty(\mathscr{H}^{\rm{in}}_\tau)}\cdot\|\Gamma^{I_3}v\|_{L^2(\mathscr{H}^{\rm{in}}_\tau)}\cdot\tau^{-3}\nonumber\\
&+&\sum_{\substack{|I_1|,|I_3|\le N-5\\|I_2|\le N-1}}\|\tau t^{\f{1}{2}}\Gamma^{I_1}u\|_{L^\infty(\mathscr{H}^{\rm{in}}_\tau)}\cdot\tau^{-\dz}\l\|\f{\tau}{t}\partial\Gamma^{I_2}\partial u\r\|_{L^2(\mathscr{H}^{\rm{in}}_\tau)}\cdot\|t^{\f{3}{2}}\Gamma^{I_3}v\|_{L^\infty(\mathscr{H}^{\rm{in}}_\tau)}\cdot\tau^{-3+\dz}\nonumber\\
&+&\sum_{\substack{|I_2|,|I_3|\le N-4\\|I_1|\le N-1}}\tau^{-\f{1}{2}-\dz}\l\|\f{\tau}{t}\Gamma^{I_1}u\r\|_{L^2(\mathscr{H}^{\rm{in}}_\tau)}\cdot\|\tau t^{\f{1}{2}}\partial\Gamma^{I_2}\partial u\|_{L^\infty(\mathscr{H}^{\rm{in}}_\tau)}\cdot\|t^{\f{3}{2}}\Gamma^{I_3}v\|_{L^\infty(\mathscr{H}^{\rm{in}}_\tau)}\cdot\tau^{-\f{5}{2}+\dz}\nonumber\\
&\lesssim&(C_2\ez)^3\tau^{-\f{5}{2}+\dz},
\eeqn
where we use \eqref{s4: inL^infty}, \eqref{s4: inL^2N-1}, \eqref{s4: incL^2N-1} and \eqref{s4: inL^2N}. Combining \eqref{s4: S3gzT}-\eqref{s4: S3gzT'''}, we obtain
\beq\label{s4: S3gzTs}
\|\Gamma^IW\|_{L^2(\mathscr{H}^{\rm{in}}_\tau)}\lesssim(C_2\ez)^3\tau^{-\f{5}{2}+\dz},\quad |I|\le N-1,\tau\in[2,s].
\eeq
By \eqref{s4: S3gzSN-1}, \eqref{s4: S3gzSN-2}, \eqref{s4: S3gzSN-3} and \eqref{s4: S3gzTs}, for any $\tau\in[2,s]$, we have
\beq\label{s4: S4gzS}
\begin{split}
\tau\|\Gamma^I(S+W)\|_{L^2(\mathscr{H}^{\rm{in}}_\tau)}&\lesssim(C_2\ez)^2\tau^{-\f{1}{2}+\dz},\quad |I|\le N-1,\\
\tau\|\Gamma^I(S+W)\|_{L^2(\mathscr{H}^{\rm{in}}_\tau)}&\lesssim(C_2\ez)^2\tau^{-1+\dz},\quad |I|\le N-2,\\
\tau\|\Gamma^I(S+W)\|_{L^2(\mathscr{H}^{\rm{in}}_\tau)}&\lesssim(C_2\ez)^2\tau^{-\f{3}{2}+\dz},\quad |I|\le N-3.
\end{split}
\eeq
We point out that the estimates of the boundary integral $\mathcal{\tilde{R}}(s,\Gamma^IU)$ in \eqref{s4: E_cinU0} (see \eqref{s4: tildeRU}) were given in  \eqref{s3: S4cgzUN-1} and \eqref{s3: S4cgzUN-2} (see \eqref{s3: S4H(t,x)}). Hence, by \eqref{s4: S4gzS}, we have
\be
\begin{split}
\sum_{|I|\le N-1}E^{c,\rm{in}}(s,\Gamma^IU)&\lesssim[\ez^2+(C_1\ez)^4]s^{2\dz}+(C_2\ez)^2\sum_{|I|\le N-1}\int_2^s[E^{c,\rm{in}}(\tau,\Gamma^IU)]^{\f{1}{2}}\cdot\tau^{-\f{1}{2}+\dz}{\rm{d}}\tau\\
&\lesssim[\ez^2+(C_1\ez)^4]s^{2\dz}+(C_2\ez)^2 s^{\f{1}{2}+\dz}\sup_{\tau\in[2,s]}\sum_{|I|\le N-1}[E^{c,\rm{in}}(\tau,\Gamma^IU)]^{\f{1}{2}},\\
\sum_{|I|\le N-2}E^{c,\rm{in}}(s,\Gamma^IU)&\lesssim\ez^2+(C_1\ez)^4+(C_2\ez)^2\sum_{|I|\le N-2}\int_2^s[E^{c,\rm{in}}(\tau,\Gamma^IU)]^{\f{1}{2}}\cdot\tau^{-1+\dz}{\rm{d}}\tau\\
&\lesssim\ez^2+(C_1\ez)^4+(C_2\ez)^2 s^{\dz}\sup_{\tau\in[2,s]}\sum_{|I|\le N-2}[E^{c,\rm{in}}(\tau,\Gamma^IU)]^{\f{1}{2}},\\
\sum_{|I|\le N-3}E^{c,\rm{in}}(s,\Gamma^IU)&\lesssim\ez^2+(C_1\ez)^4+(C_2\ez)^2\sum_{|I|\le N-3}\int_2^s[E^{c,\rm{in}}(\tau,\Gamma^IU)]^{\f{1}{2}}\cdot\tau^{-\f{3}{2}+\dz}{\rm{d}}\tau\\
&\lesssim\ez^2+(C_1\ez)^4+(C_2\ez)^2\sup_{\tau\in[2,s]}\sum_{|I|\le N-3}[E^{c,\rm{in}}(\tau,\Gamma^IU)]^{\f{1}{2}},
\end{split}
\ee
which imply
\beq\label{s4: S4cUN-1}
\begin{split}
\sum_{|I|\le N-1}[E^{c,\rm{in}}(s,\Gamma^IU)]^{\f{1}{2}}&\lesssim\ez s^{\dz}+(C_2\ez)^2s^{\f{1}{2}+\dz},\\
\sum_{|I|\le N-2}[E^{c,\rm{in}}(s,\Gamma^IU)]^{\f{1}{2}}&\lesssim\ez+(C_2\ez)^2s^\dz,\\
\sum_{|I|\le N-3}[E^{c,\rm{in}}(s,\Gamma^IU)]^{\f{1}{2}}&\lesssim\ez+(C_2\ez)^2.
\end{split}
\eeq
We recall the definitions \eqref{s2: E^c,in_dz} and \eqref{s2: E^c,in_0}. We note that on $\mathscr{H}^{\rm{in}}_s$ (recall \eqref{L_0}) we have
\be
|L_0w|+|x^a\bar{\partial}_aw|+|w|+s|\bar{\nabla}w|\lesssim(s^2/t)|\partial w|+|Lw|+|w|
\ee
for any smooth function $w$. We recall that $U-u=F_u$ by \eqref{s3: defU}. For $|I|\le N-1$, we have
\beqn\label{s4: S4E_cd}
[E^{c,\rm{in}}(s, \Gamma^IF_u)]^{\f{1}{2}}&\lesssim&\|s^2t^{-1}\big(|\partial\Gamma^I(uv)|+|\partial\Gamma^I(u\partial v)|+|\partial\Gamma^I(\partial uv)|+|\partial\Gamma^I(\partial u\partial v)|\big)\|_{L^2(\mathscr{H}^{\rm{in}}_s)}\nonumber\\
&+&\||L\Gamma^I(uv)|+|L\Gamma^I(u\partial v)|+|L\Gamma^I(\partial uv)|+|L\Gamma^I(\partial u\partial v)|\|_{L^2(\mathscr{H}^{\rm{in}}_s)}\nonumber\\
&+&\|\Gamma^IF_u\|_{L^2(\mathscr{H}^{\rm{in}}_s)}:=II_1+II_2+II_3.
\eeqn
We have
\beqn\label{s4: S4E_cd'}
II_1&\lesssim&\sum_{\substack{|I_1|\le N-4\\|I_2|\le N-1}}\|st^{\f{1}{2}}\big(|\partial\Gamma^{I_1}u|+|\partial\Gamma^{I_1}\partial u|\big)\|_{L^\infty(\mathscr{H}^{\rm{in}}_s)}\cdot\l\|\f{s}{t}\big(|\Gamma^{I_2}v|+|\Gamma^{I_2}\partial v|\big)\r\|_{L^2(\mathscr{H}^{\rm{in}}_s)}\cdot s^{-\f{1}{2}}\nonumber\\
&+&\sum_{\substack{|I_2|\le N-4\\|I_1|\le N}}s^{-\dz}\l\|\f{s}{t}\partial\Gamma^{I_1}u\r\|_{L^2(\mathscr{H}^{\rm{in}}_s)}\cdot\|t^{\f{3}{2}}\big(|\Gamma^{I_2}v|+|\Gamma^{I_2}\partial v|\big)\|_{L^\infty(\mathscr{H}^{\rm{in}}_s)}\cdot s^{-\f{1}{2}+\dz}\nonumber\\
&+&\sum_{\substack{|I_1|\le N-5\\|I_2|\le N-1}}\|st^{\f{1}{2}}\big(|\Gamma^{I_1}u|+|\partial\Gamma^{I_1}u|\big)\|_{L^\infty(\mathscr{H}^{\rm{in}}_s)}\cdot s^{-\dz}\l\|\f{s}{t}\big(|\partial\Gamma^{I_2}v|+|\partial\Gamma^{I_2}\partial v|\big)\r\|_{L^2(\mathscr{H}^{\rm{in}}_s)}\cdot s^{-\f{1}{2}+\dz}\nonumber\\
&+&\!\!\!\sum_{\substack{|I_2|\le N-5\\|J|=1,|I'_1|\le N-2\\|I_1|\le N-1}}\Big\{s^{-\dz}\l\|\f{s}{t}\Gamma^J\Gamma^{I'_1}u\r\|_{L^2(\mathscr{H}^{\rm{in}}_s)}+\l\|\f{s}{t}\partial\Gamma^{I_1}u\r\|_{L^2(\mathscr{H}^{\rm{in}}_s)}\Big\}\cdot\|t^{\f{3}{2}}\big(|\partial\Gamma^{I_2}v|+|\partial\Gamma^{I_2}\partial v|\big)\|_{L^\infty(\mathscr{H}^{\rm{in}}_s)}\cdot s^{-\f{1}{2}+\dz},\nonumber\\
&\lesssim&(C_2\ez)^2,
\eeqn
\beqn\label{s4: S4E_cd''}
II_2&\lesssim&\sum_{\substack{|I_1|\le N-5\\|I_2|\le N-1}}\|st^{\f{1}{2}}\big(|L\Gamma^{I_1}u|+|L\Gamma^{I_1}\partial u|\big)\|_{L^\infty(\mathscr{H}^{\rm{in}}_s)}\cdot\l\|\f{s}{t}\big(|\Gamma^{I_2}v|+|\Gamma^{I_2}\partial v|\big)\r\|_{L^2(\mathscr{H}^{\rm{in}}_s)}\cdot s^{-1}\nonumber\\
&+&\sum_{\substack{|I_2|\le N-4\\|I_1|\le N-1}}\Big\{s^{-\f{1}{2}-\dz}\l\|\f{s}{t}L\Gamma^{I_1}u\r\|_{L^2(\mathscr{H}^{\rm{in}}_s)}+s^{-\dz}\l\|\f{s}{t}L\Gamma^{I_1}\partial u\r\|_{L^2(\mathscr{H}^{\rm{in}}_s)}\Big\}\cdot\|t^{\f{3}{2}}\big(|\Gamma^{I_2}v|+|\Gamma^{I_2}\partial v|\big)\|_{L^\infty(\mathscr{H}^{\rm{in}}_s)}\cdot s^{-1+\dz}\nonumber\\
&+&\sum_{\substack{|I_1|\le N-5\\|I_2|\le N-1}}\|st^{\f{1}{2}}\big(|\Gamma^{I_1}u|+|\partial\Gamma^{I_1}u|\big)\|_{L^\infty(\mathscr{H}^{\rm{in}}_s)}\cdot s^{-\dz}\l\|\f{s}{t}\big(|L\Gamma^{I_2}v|+|L\Gamma^{I_2}\partial v|\big)\r\|_{L^2(\mathscr{H}^{\rm{in}}_s)}\cdot s^{-1+\dz}\nonumber\\
&+&\!\!\!\sum_{\substack{|I_2|\le N-5\\|I_1|\le N-1}}\Big\{s^{-\f{1}{2}-\dz}\l\|\f{s}{t}\Gamma^{I_1}u\r\|_{L^2(\mathscr{H}^{\rm{in}}_s)}+\l\|\f{s}{t}\partial\Gamma^{I_1}u\r\|_{L^2(\mathscr{H}^{\rm{in}}_s)}\Big\}\cdot\|t^{\f{3}{2}}\big(|L\Gamma^{I_2}v|+|L\Gamma^{I_2}\partial v|\big)\|_{L^\infty(\mathscr{H}^{\rm{in}}_s)}\cdot s^{-1+\dz},\nonumber\\
&\lesssim&(C_2\ez)^2,
\eeqn
and $II_3\lesssim(C_2\ez)^2$, where we use \eqref{s4: inL^infty}, \eqref{s4: inL^2N-1}, \eqref{s4: inL^2N}, \eqref{s4: incL^2N-2}, \eqref{s4: incL^2N-1} and \eqref{s4: improN-1}. Combining \eqref{s4: S4E_cd}-\eqref{s4: S4E_cd''} with \eqref{s4: S4cUN-1}, we obtain 
\beq\label{s4: S4us}
\begin{split}
\sum_{|I|\le N-1}[E^{c,\rm{in}}(s,\Gamma^Iu)]^{\f{1}{2}}&\lesssim[\ez+(C_2\ez)^2]s^{\f{1}{2}+\dz},\\
\sum_{|I|\le N-2}[E^{c,\rm{in}}(s,\Gamma^Iu)]^{\f{1}{2}}&\lesssim[\ez+(C_2\ez)^2]s^\dz,\\
\sum_{|I|\le N-3}[E^{c,\rm{in}}(s,\Gamma^Iu)]^{\f{1}{2}}&\lesssim\ez+(C_2\ez)^2.
\end{split}
\eeq

Combining ${\bf Steps\ 1-2}$ (see \eqref{s4: S1us}, \eqref{s4: S2vs} and \eqref{s4: S4us}), we have strictly improved the bootstrap estimates \eqref{s4: inuvL^2N}-\eqref{s4: inucN-2} for $C_2$ sufficiently large and $0<\ez\ll C_2^{-1}$ sufficiently small. Hence, the proof of Proposition \ref{s4: maxtildeT} is completed. We conclude that the solution $(u,v)$ exists globally in time in the interior region $\mathscr{D}^{\rm{in}}$ and satisfies \eqref{s4: inuvL^2N}-\eqref{s4: inucN-2} for all $s\in[2,\infty)$.

\section{Scattering Results}\label{s5}

In this section we give the proof of Theorem \ref{thm2}. 

Let $(u,v)$ be the global solution to \eqref{s3: equv}-\eqref{s1: ini} given by Theorem \ref{thm1}. We denote 
\beq\label{s5: vecu}
\vec{u}=(u,\partial_tu)'=\l(\begin{array}{c}
u\\
\partial_tu
\end{array}\r),\quad\quad\vec{v}=(v,\partial_tv)'=\l(\begin{array}{c}
v\\
\partial_tv
\end{array}\r),  
\eeq
where $(a_1,a_2)'$ denotes the transpose of a vector $\vec{a}=(a_1,a_2)$ in $\mathbb{R}^2$. We also set
\beq\label{s5: vecF}
\vec{F}_u=(0,F_u)',\quad\quad\vec{F}_v=(0,F_v)',\quad\quad\vec{u}_0=(u_0,u_1)',\quad\quad\vec{v}_0=(v_0,v_1)'.
\eeq 
By the linear theory of wave and Klein-Gordon equations, we can write
\beq\label{s5: vecu=}
\vec{u}=\mathcal{S}(t-2)\vec{u}_0+\int_2^t\mathcal{S}(t-\tau)\vec{F}_u(\tau){\rm{d}}\tau,
\eeq
\beq\label{s5: vecv=}
\vec{v}=\mathcal{\tilde{S}}(t-2)\vec{v}_0+\int_2^t\mathcal{\tilde{S}}(t-\tau)\vec{F}_v(\tau){\rm{d}}\tau,
\eeq
where
\be
\mathcal{S}(t)=\l(\begin{array}{cc}
\cos(t\sqrt{-\Delta})&\f{\sin(t\sqrt{-\Delta})}{\sqrt{-\Delta}}\\
-\sqrt{-\Delta}\sin(t\sqrt{-\Delta})&\cos(t\sqrt{-\Delta})
\end{array}\r),\quad\quad\mathcal{\tilde{S}}(t)=\l(\begin{array}{cc}
\cos(t\langle\nabla\rangle)&\f{\sin(t\langle\nabla\rangle)}{\langle\nabla\rangle}\\
-\langle\nabla\rangle\sin(t\langle\nabla\rangle)&\cos(t\langle\nabla\rangle)
\end{array}\r).
\ee

Let $l\in\mathbb{N}$. We denote ${\bf H}^l(\mathbb{R}^3):=H^{l+1}(\mathbb{R}^3)\times H^l(\mathbb{R}^3)$ and ${\bf{H}}_{l}(\mathbb{R}^3):=\big(\dot{H}^{l+1}(\mathbb{R}^3)\cap\dot{H}^1(\mathbb{R}^3)\big)\times H^l(\mathbb{R}^3)$, where $H^k(\mathbb{R}^3), \dot{H}^{k}(\mathbb{R}^3), k\in\mathbb{N}$ denote the Sobolev spaces and homogeneous Sobolev spaces respectively.

\begin{lem}\label{Sca}
The following statements hold:
\begin{itemize}
\item[i)] Let $l\in\mathbb{N}$ and $\dz>0$. For any $\mathbb{R}^2$-valued function $\vec{f}(\tau,x)=(f_1,f_2)'$ which is defined in $[2,\infty)\times\mathbb{R}^3$ and satisfies $\vec{f}(\tau,\cdot)\in{\bf{H}}^l(\mathbb{R}^3)$ for any fixed $\tau\in[2,\infty)$, and any $t\ge 2$, $4\le T_1<T_2<\infty$, we have
\ben
&&\l\|\int_{T_1}^{T_2}\mathcal{S}(t-\tau)(\vec{f}(\tau)){\rm{d}}\tau\r\|_{{\bf{H}}_{l}(\mathbb{R}^3)}\\
&\lesssim& T_1^{-\f{\dz}{2}}\sum_{k=0}^l\Bigg\{\l(\int_{T_1}^{T_2}\||\vec{\nabla}^k\vec{f}|\|_{L^2_x(\Sigma^{\rm{ex}}_\tau)}^2 \cdot\tau^{1+\dz}{\rm{d}}\tau\r)^{\f{1}{2}}+\l(\int_{T_1^{\f{1}{2}}}^{T_2}\||\vec{\nabla}^k\vec{f}|\|_{L^2(\mathscr{H}^{\rm{in}}_s)}^2\cdot s^{1+2\dz}{\rm{d}}s\r)^{\f{1}{2}}\Bigg\},\\
&&\l\|\int_{T_1}^{T_2}\mathcal{\tilde{S}}(t-\tau)(\vec{f}(\tau)){\rm{d}}\tau\r\|_{{\bf{H}}^l(\mathbb{R}^3)}\\
&\lesssim& T_1^{-\f{\dz}{2}}\sum_{k=0}^l\Bigg\{\l(\int_{T_1}^{T_2}\||\vec{\nabla}^k\vec{f}|+|f_1|\|_{L^2_x(\Sigma^{\rm{ex}}_\tau)}^2 \cdot\tau^{1+\dz}{\rm{d}}\tau\r)^{\f{1}{2}}+\l(\int_{T_1^{\f{1}{2}}}^{T_2}\||\vec{\nabla}^k\vec{f}|+|f_1|\|_{L^2(\mathscr{H}^{\rm{in}}_s)}^2\cdot s^{1+2\dz}{\rm{d}}s\r)^{\f{1}{2}}\Bigg\},
\een
where $|\vec{\nabla}^k\vec{f}|:=|\nabla^{k+1}f_1|+|\nabla^kf_2|$ and $\nabla=(\partial_1,\partial_2,\partial_3)$.

\item[ii)] Let $l\in\mathbb{N}$ and $\vec{u},\vec{v},\vec{F}_u,\vec{F}_v,\vec{u}_0,\vec{v}_0$ be as in \eqref{s5: vecu}-\eqref{s5: vecv=} with $\vec{u}_0\in{\bf{H}}_{l}(\mathbb{R}^3), \vec{v}_0\in{\bf{H}}^{l}(\mathbb{R}^3)$ and $F_u(\tau),F_v(\tau)\in H^l(\mathbb{R}^3)$ for any fixed $\tau\in[2,\infty)$. If for some $\dz>0$, it holds that
\be
M:=\int_2^4\|F(\tau,\cdot)\|_{L^2_x(\mathbb{R}^3)}{\rm{d}}\tau+\l(\int_{4}^{+\infty}\|F(\tau,\cdot)\|_{L^2_x(\Sigma^{\rm{ex}}_\tau)}^2 \cdot\tau^{1+\dz}{\rm{d}}\tau+\int_{2}^{+\infty}\|F\|_{L^2(\mathscr{H}^{\rm{in}}_s)}^2\cdot s^{1+2\dz}{\rm{d}}s\r)^{\f{1}{2}}<\infty,
\ee
where $F:=\sum_{k=0}^{l}\big(|\nabla^{k}F_u|+|\nabla^kF_v|\big)$, then the solution $(\vec{u},\vec{v})$ scatters to a free solution in ${\bf{H}}_{l}(\mathbb{R}^3)\times {\bf{H}}^{l}(\mathbb{R}^3)$, i.e., there exist $\vec{u}^*_0=(u^*_0,u^*_1)'\in{\bf{H}}_{l}(\mathbb{R}^3)$ and $\vec{v}^*_0=(v^*_0,v^*_1)'\in{\bf{H}}^{l}(\mathbb{R}^3)$ such that
\ben
\lim_{t\to+\infty}\|\vec{u}-\vec{u}^*\|_{{\bf{H}}_{l}(\mathbb{R}^3)}=0,\quad\quad\lim_{t\to+\infty}\|\vec{v}-\vec{v}^*\|_{{\bf{H}}^{l}(\mathbb{R}^3)}=0,
\een
where $\vec{u}^*=(u^*,\partial_tu^*)'$, $\vec{v}^*=(v^*,\partial_tv^*)'$, and $(u^*,v^*)$ is the solution to the $3D$ linear homogeneous wave-Klein-Gordon system with the initial data $(u^*_0,u^*_1,v^*_0,v^*_1)$.
\end{itemize}

\begin{proof}
$i)$ We only consider the case $l=0$. The conclusion for $l\ge 1$ follows from this and the definitions of the spaces ${\bf{H}}_l(\mathbb{R}^3)$ and ${\bf{H}}^l(\mathbb{R}^3)$. For any fixed $t\ge 2$, let $\vec{U}(\tau):=\mathcal{S}(t-\tau)(\vec{f}(\tau))=(U_1,U_2)'$. For any $4\le T_1<T_2<\infty$, by standard energy inequalities, we have
\ben
\l\|\int_{T_1}^{T_2}\vec{U}(\tau){\rm{d}}\tau\r\|_{{\bf{H}}_0(\mathbb{R}^3)}&\lesssim&\l(\int_{\mathbb{R}^3}\Bigg\{\l|\int_{T_1}^{T_2}\nabla U_1(\tau){\rm{d}}\tau\r|^2+\l|\int_{T_1}^{T_2}U_2(\tau){\rm{d}}\tau\r|^2\Bigg\}{\rm{d}}x\r)^{\f{1}{2}}\\
&\lesssim&\l(\int_{\mathbb{R}^3}\l\{\int_{T_1}^{T_2}|\nabla U_1(\tau)\cdot \tau^{\f{1+\dz}{2}}|^2{\rm{d}}\tau+\int_{T_1}^{T_2}|U_2(\tau)\cdot\tau^{\f{1+\dz}{2}}|^2{\rm{d}}\tau\r\}\cdot\int_{T_1}^{T_2}\tau^{-(1+\dz)}{\rm{d}}\tau{\rm{d}}x\r)^{\f{1}{2}}\\
&\lesssim&T_1^{-\f{\dz}{2}}\l(\int_{T_1}^{T_2}\int_{\mathbb{R}^3}\l\{|\nabla U_1(\tau)|^2+|U_2(\tau)|^2\r\}\cdot \tau^{1+\dz}{\rm{d}}x{\rm{d}}\tau\r)^{\f{1}{2}}\\
&\lesssim&II:=T_1^{-\f{\dz}{2}}\l(\int_{T_1}^{T_2}\l(\int_{r\ge \tau-1}+\int_{r<\tau-1}\r)\l\{|\nabla f_1(\tau)|^2+|f_2(\tau)|^2\r\}\cdot \tau^{1+\dz}{\rm{d}}x{\rm{d}}\tau\r)^{\f{1}{2}}.
\een
By a change of variables $(\tau,x)\to(s,x)$ with $s=\sqrt{\tau^2-|x|^2}$, we obtain
\ben
II&\lesssim&T_1^{-\f{\dz}{2}}\l(\int_{T_1}^{T_2}\int_{\Sigma^{\rm{ex}}_\tau}\l\{|\nabla f_1|^2+|f_2|^2\r\}(\tau,x)\cdot \tau^{1+\dz}{\rm{d}}x{\rm{d}}\tau\r)^{\f{1}{2}}\\
&+&T_1^{-\f{\dz}{2}}\l(\int_{T_1^{\f{1}{2}}}^{T_2}\int_{r<\f{s^2-1}{2}}\l\{|\nabla f_1|^2+|f_2|^2\r\}(\sqrt{s^2+|x|^2},x)\cdot \tau^{1+\dz}\f{s}{\tau}{\rm{d}}x{\rm{d}}s\r)^{\f{1}{2}}\\
&\lesssim&T_1^{-\f{\dz}{2}}\Bigg\{\l(\int_{T_1}^{T_2}\||\nabla f_1|+|f_2|\|_{L^2_x(\Sigma^{\rm{ex}}_\tau)}^2 \cdot\tau^{1+\dz}{\rm{d}}\tau\r)^{\f{1}{2}}+\l(\int_{T_1^{\f{1}{2}}}^{T_2}\||\nabla f_1|+|f_2|\|_{L^2(\mathscr{H}^{\rm{in}}_s)}^2\cdot s^{1+2\dz}{\rm{d}}s\r)^{\f{1}{2}}\Bigg\}.
\een
Similarly,
\ben
&&\l\|\int_{T_1}^{T_2}\mathcal{\tilde{S}}(t-\tau)(\vec{f}(\tau)){\rm{d}}\tau\r\|_{{\bf{H}}^0(\mathbb{R}^3)}\\
&\lesssim& T_1^{-\f{\dz}{2}}\Bigg\{\l(\int_{T_1}^{T_2}\||\nabla f_1|+|f_2|+|f_1|\|_{L^2_x(\Sigma^{\rm{ex}}_\tau)}^2 \cdot\tau^{1+\dz}{\rm{d}}\tau\r)^{\f{1}{2}}+\l(\int_{T_1^{\f{1}{2}}}^{T_2}\||\nabla f_1|+|f_2|+|f_1|\|_{L^2(\mathscr{H}^{\rm{in}}_s)}^2\cdot s^{1+2\dz}{\rm{d}}s\r)^{\f{1}{2}}\Bigg\}.
\een
$ii)$ Let $\vec{u}$, $\vec{v}$, $\vec{u}_0$, $\vec{v}_0$, $\vec{F}_u$, $\vec{F}_v$ be as in \eqref{s5: vecu}-\eqref{s5: vecv=} and 
\ben
\vec{u}^*_0:&=&\vec{u}_0+\int_2^{+\infty}\mathcal{S}(2-\tau)\vec{F}_u(\tau){\rm{d}}\tau,\quad\quad \vec{u}^*=\mathcal{S}(t-2)\vec{u}^*_0,\\
\vec{v}^*_0:&=&\vec{v}_0+\int_2^{+\infty}\mathcal{\tilde{S}}(2-\tau)\vec{F}_v(\tau){\rm{d}}\tau,\quad\quad \vec{v}^*=\mathcal{\tilde{S}}(t-2)\vec{v}^*_0.
\een
For any $4\le T_1<T_2<\infty$, by $i)$, we have
\ben
&&\l\|\int_{T_1}^{T_2}\mathcal{S}(2-\tau)\vec{F}_u(\tau){\rm{d}}\tau\r\|_{{\bf{H}}_{l}(\mathbb{R}^3)}+\l\|\int_{T_1}^{T_2}\mathcal{\tilde{S}}(2-\tau)\vec{F}_v(\tau){\rm{d}}\tau\r\|_{{\bf{H}}^{l}(\mathbb{R}^3)}\\
&\lesssim& T_1^{-\f{\dz}{2}}\sum_{k=0}^{l}\Bigg\{\l(\int_{T_1}^{T_2}\||\nabla^{k}F_u|+|\nabla^kF_v|\|_{L^2_x(\Sigma^{\rm{ex}}_\tau)}^2 \cdot\tau^{1+\dz}{\rm{d}}\tau\r)^{\f{1}{2}}+\l(\int_{T_1^{\f{1}{2}}}^{T_2}\||\nabla^kF_u|+|\nabla^kF_v|\|_{L^2(\mathscr{H}^{\rm{in}}_s)}^2\cdot s^{1+2\dz}{\rm{d}}s\r)^{\f{1}{2}}\Bigg\}\\
&\lesssim&MT_1^{-\f{\dz}{2}}\to 0\quad\mathrm{as}\ \  T_1\to+\infty.
\een
Hence, $\vec{u}^*_0$ and $\vec{v}^*_0$ are well-defined in ${\bf{H}}_{l}(\mathbb{R}^3)$ and ${\bf{H}}^{l}(\mathbb{R}^3)$ respectively. Similarly, using $i)$ again, we have
\ben
&&\|\vec{u}-\vec{u}^*\|_{{\bf{H}}_{l}(\mathbb{R}^3)}=\l\|\int_t^{+\infty}\mathcal{S}(t-\tau)\vec{F}_u(\tau){\rm{d}}\tau\r\|_{{\bf{H}}_{l}(\mathbb{R}^3)}\lesssim Mt^{-\f{\dz}{2}}\to 0\quad\mathrm{as}\ \ t\to+\infty,\\
&&\|\vec{v}-\vec{v}^*\|_{{\bf{H}}^{l}(\mathbb{R}^3)}=\l\|\int_t^{+\infty}\mathcal{\tilde{S}}(t-\tau)\vec{F}_v(\tau){\rm{d}}\tau\r\|_{{\bf{H}}^{l}(\mathbb{R}^3)}\lesssim Mt^{-\f{\dz}{2}}\to 0\quad\mathrm{as}\ \ t\to+\infty.
\een
\end{proof}
\end{lem}

$\textit{Proof\ of\ Theorem\ \ref{thm2}}.$

By the exterior estimates \eqref{ScaexvN} and \eqref{ScaexdvN-1}, for any $t\in[2,\infty)$, we have
\beq\label{s5: exFN}
\||\Gamma^I(uv)|+|\Gamma^I(\partial uv)|\|_{L^2_x(\Sigma^{\rm{ex}}_\tau)}\lesssim (C_1\ez)^2\big\{\sqrt{l(\tau)}+\tau^{-1}\big\}\tau^{-1+\f{\dz}{2}},\quad |I|\le N, \tau\in[2,t],
\eeq
\beqn\label{s5: exFN-1}
\||\Gamma^I(uv)|+|\Gamma^I(\partial uv)|&+&|\Gamma^I(u\partial v)|+|\Gamma^I(\partial u\partial v)|\|_{L^2_x(\Sigma^{\rm{ex}}_\tau)}\nonumber\\
&\lesssim& (C_1\ez)^2\big\{\sqrt{l(\tau)}+\tau^{-1}\big\}\tau^{-1+\f{\dz}{2}},\quad |I|\le N-1, \tau\in[2,t].
\eeqn

If $F_u=P_0uv+P_1^\az\partial_\az uv$, then the highest order energy of $u$ in the exterior region is bounded (which improves \eqref{s3: S1gzuN}), i.e.,
\beq\label{s5: exF_u=uv}
\|\Gamma^Iu\|_{X^{{\rm{ex}},\lz}_{0,0,t}}\lesssim\ez+(C_1\ez)^2,\quad|I|\le N.
\eeq
Hence the conformal energy of $u$ up to order $N-1$ is indeed bounded (which improves \eqref{s3: S4cgzuN-12}), i.e.,
\beq\label{s5: cgzu}
[E^{c,\rm{ex}}(t,\Gamma^Iu)]^{\f{1}{2}}\lesssim\ez+(C_1\ez)^2,\quad|I|\le N-1.
\eeq
Moreover, in this case, using \eqref{s3: exuvL^2N'} in the proof of \eqref{s3: cgzudv}, we find that \eqref{s3: cgzS} can also be improved, i.e.,
\beq\label{s5: cgzS}
\int_2^t\|r\Gamma^IS\|_{L^2_x(\Sigma^{\rm{ex}}_\tau)}{\rm{d}}\tau\lesssim (C_1\ez)^2,\quad |I|\le N-1
\eeq
and therefore we can improve \eqref{s3: S4cgzUN-1}, i.e.,
\beq\label{s5: cgzU}
[E^{c,\rm{ex}}(t,\Gamma^IU)]^{\f{1}{2}}+\mathcal{R}_0(t,\Gamma^IU)\lesssim\ez+(C_1\ez)^2,\quad|I|\le N-1,
\eeq
where $\mathcal{R}_0(t,\Gamma^IU)$ is as in \eqref{s3: S4H(t,x)}. Then the energy estimates on exterior hyperboloids \eqref{s3: exhuN} and \eqref{s3: exEch} can be improved:
\beq\label{s5: exhu}
\begin{split}
&E^{{\rm{ex}},h}_{0}(s,\Gamma^Iu)+\int_{\mathscr{l}_{[2,\f{s^2+1}{2}]}}|G\Gamma^Iu|^2{\rm{d}}\sz
\lesssim\ez^2+(C_1\ez)^3,\quad |I|\le N,\\
&\l(\int_{\mathscr{H}^{\rm{ex}}_s}\big(|(L_0+x^a\bar{\partial}_a+2)\Gamma^Iu|^2+s^2|\bar{\nabla}\Gamma^Iu|^2\big){\rm{d}}x\r)^{\f{1}{2}}\lesssim\ez+(C_1\ez)^2,\quad |I|\le N-1.
\end{split}
\eeq

For the interior estimates, we recall from \eqref{s4: improN-1} that for any $s\in[2,\infty)$, we have
\beq\label{s5: inFN-1}
\||\Gamma^I(uv)|+|\Gamma^I(\partial uv)|+|\Gamma^I(u\partial v)|+|\Gamma^I(\partial u\partial v)|\|_{L^2(\mathscr{H}^{\rm{in}}_\tau)}\lesssim(C_2\ez)^2\tau^{-\f{3}{2}+\dz},\quad |I|\le N-1, \tau\in[2,s].
\eeq
In particular, if $F_u=P_0uv$, we improve the bootstrap assumptions \eqref{s4: inuvL^2N} and \eqref{s4: inucN-1} as follows
\beq\label{s5: inbsimpro}
\sum_{|I|\le N}\big\{[E^{{\rm{in}},h}_0(s,\Gamma^Iu)]^{\f{1}{2}}+[E^{{\rm{in}},h}_1(s,\Gamma^Iv)]^{\f{1}{2}}\big\}+\sum_{|I|\le N-1}[E^{c,\rm{in}}(s,\Gamma^Iu)]^{\f{1}{2}}\le C_2\ez.
\eeq
Then using \eqref{s5: exF_u=uv}, \eqref{s5: cgzU} and \eqref{s5: exhu}, we can improve the estimates \eqref{s4: improN} and \eqref{s4: S3gzudv} as follows
\beqn\label{s5: inuvimpro}
&&\||\Gamma^I(uv)|+|\Gamma^I(\partial uv)|\|_{L^2_x(\Sigma^{\rm{ex}}_\tau)}\lesssim(C_2\ez)^2\tau^{-\f{3}{2}},\quad |I|\le N,\tau\in[2,s],\\
&&\|\Gamma^I(\partial_\gz u\partial^\gz v)\|_{L^2_x(\Sigma^{\rm{ex}}_\tau)}\lesssim(C_2\ez)^2\tau^{-\f{5}{2}},\quad |I|\le N-1,\tau\in[2,s].\nonumber
\eeqn
Hence by similar arguments as in Section \ref{ss4.2}, the bootstrap assumption \eqref{s5: inbsimpro} is closed, which implies that \eqref{s5: inbsimpro}-\eqref{s5: inuvimpro} hold for all $s\in[2,\infty)$.

For any $4\le T_1<\infty$, by \eqref{s5: exFN-1} and \eqref{s5: inFN-1}, we have
\ben
\sum_{k=0}^{N-1}\l(\int_{T_1}^{+\infty}\||\nabla^{k}F_u|+|\nabla^kF_v|\|_{L^2_x(\Sigma^{\rm{ex}}_\tau)}^2 \cdot\tau^{1+\dz}{\rm{d}}\tau+\int_{T_1^{\f{1}{2}}}^{+\infty}\||\nabla^kF_u|+|\nabla^kF_v|\|_{L^2(\mathscr{H}^{\rm{in}}_s)}^2\cdot s^{1+2\dz}{\rm{d}}s\r)^{\f{1}{2}}\lesssim T_1^{-\f{1}{4}+\dz}.
\een
If we assume further $F_u=P_0uv$, by \eqref{s5: exFN} and \eqref{s5: inuvimpro}, we obtain
\ben
\sum_{k=0}^{N}\l(\int_{T_1}^{+\infty}\||\nabla^{k}F_u|+|\nabla^kF_v|\|_{L^2_x(\Sigma^{\rm{ex}}_\tau)}^2 \cdot\tau^{1+\dz}{\rm{d}}\tau+\int_{T_1^{\f{1}{2}}}^{+\infty}\||\nabla^kF_u|+|\nabla^kF_v|\|_{L^2(\mathscr{H}^{\rm{in}}_s)}^2\cdot s^{1+2\dz}{\rm{d}}s\r)^{\f{1}{2}}\lesssim T_1^{-\f{1}{4}+\f{\dz}{2}}.
\een
By Lemma \ref{Sca}, we obtain the conclusion. $\hfill\square$

\begin{appendices}

\section{}\label{sA}

$\textit{Proof\ of\ Lemma\ \ref{s2: Sobex}.}$ 

Let $(r,\varphi,\tz)$ be the spherical coordinates in $\mathbb{R}^3$, i.e., $x=(r\sin\varphi\cos\tz,r\sin\varphi\sin\tz,r\cos\varphi)$. We can write the surface element on $\mathbb{S}^2$ as ${\rm{d}}\sz=\sin\varphi{\rm{d}}\varphi{\rm{d}}\tz$. By Sobolev imbedding theorem,
\beq\label{s2: Sexeq0}
\sup_{\mathbb{S}^2}|u(t,r,\varphi,\tz)|^2\lesssim\sum_{0\le k+l\le 2}\int_{\mathbb{S}^2}|\partial^k_\varphi\partial^l_\tz u(t,r,\varphi,\tz)|^2{\rm{d}}\sz.
\eeq
For any fixed $t\ge 2$ we denote $\xi(r)=2+r-t$ and $u_{k,l}(r,\varphi,\tz):=\partial^k_\varphi\partial^l_\tz u(t,r,\varphi,\tz)$. Let $v$ be any smooth function which is compactly supported in $x$, then in the region $\{r\ge t-1\}$ we have
\ben
\partial_r\big((\xi(r))^\Lz r^2v^2\big)\ge\Lz(\xi(r))^{\Lz-1}r^2 v^2+2(\xi(r))^{\Lz}r^2v\partial_rv.
\een
Hence for each $x\in\Sigma^{\rm{ex}}_t$, we have
\beq\label{s2: HS1}
(\xi(r))^{\Lz} r^2v^2(t,x)\lesssim\int_r^\infty\l\{(\xi(\rz))^{\Lz-1}\rz^2v^2+(\xi(\rz))^{\Lz+1}\rz^2|\partial_\rz v|^2\r\}{\rm{d}}\rz.
\eeq
If $u$ is compactly supported in $x$, we obtain
\ben
(\xi(r))^{\Lz} r^2|u_{k,l}(r,\varphi,\tz)|^2
\lesssim\int_r^\infty\!\!\rz^2\!\l\{(\xi(\rz))^{\Lz-1}|u_{k,l}(\rz,\varphi,\tz)|^2+(\xi(\rz))^{\Lz+1}|\partial_\rz u_{k,l}(\rz,\varphi,\tz)|^2\r\}{\rm{d}}\rz,
\een
for any $0\le k+l\le 2$. Integrating over $\mathbb{S}^2$ and using \eqref{s2: Sexeq0}, we derive
\beq\label{s2: Sexeq1}
\sup_{\Sigma^{\rm{ex}}_t}(\xi(r))^{\Lz} r^2|u(t,x)|^2
\lesssim\sum_{0\le k+l\le 2}\int_{\Sigma^{\rm{ex}}_t}\l\{(\xi(\rz))^{\Lz-1}|\partial^k_\varphi\partial^l_\tz u(t,x)|^2+(\xi(\rz))^{\Lz+1}|\partial_\rz\partial^k_\varphi\partial^l_\tz u(t,x)|^2\r\}{\rm{d}}x.
\eeq
In the general case where $u$ is not compactly supported in $x$, we choose a cut-off function $\chi\in C^\infty_0(\mathbb{R})$ and apply \eqref{s2: Sexeq1} to $\chi(\ez r)u$ for any $\ez>0$, then we derive
\ben
&&\sup_{\Sigma^{\rm{ex}}_t}(\xi(r))^{\Lz} r^2 |\chi(\ez r)u(t,x)|^2\\
&\lesssim& \sum_{0\le k+l\le 2}\int_{\Sigma^{\rm{ex}}_t}\big\{(\xi(\rz))^{\Lz-1}|\chi(\ez\rz)\partial^k_\varphi\partial^l_\tz u|^2+(\xi(\rz))^{\Lz+1}\big(|\chi(\ez\rz)\partial_\rz\partial^k_\varphi\partial^l_\tz u|^2+\ez^2|\chi'(\ez\rz)|^2\cdot|\partial^k_\varphi\partial^l_\tz u|^2\big)\big\}{\rm{d}}x.
\een
Note that on $\Sigma^{\rm{ex}}_t$ we have $|\ez\chi'(\ez \rz)|^2\le(\xi(\rz))^{2}|\ez\chi'(\ez \rz)|^2\le|\ez \rz\chi'(\ez \rz)|^2\lesssim 1$. Let $\ez\to 0$, we obtain \eqref{s2: Sexeq1} as well. Taking the same weight $(\xi(\rz))^{\Lz}$ on the right-hand side of \eqref{s2: HS1} and arguing as above, we also obtain
\be
\sup_{\Sigma^{\rm{ex}}_t}(\xi(r))^{\Lz} r^2 |u(t,x)|^2 \lesssim\sum_{0\le k+l\le 2}\int_{\Sigma^{\rm{ex}}_t}(\xi(\rz))^{\Lz}\l\{|\partial^k_\varphi\partial^l_\tz u(t,x)|^2+|\partial_\rz\partial^k_\varphi\partial^l_\tz u(t,x)|^2\r\}{\rm{d}}x.
\ee
By direct calculation,
\ben
&&\partial_\tz=\Omega_{12},\quad\quad\partial_\varphi=-\sin\tz\ \Omega_{23}-\cos\tz\ \Omega_{13},\\
&&\partial_\tz\partial_\tz=\Omega_{12}\Omega_{12},\quad\quad\partial_\tz\partial_\varphi=-\sin\tz\ \Omega_{23}\Omega_{12}-\cos\tz\ \Omega_{13}\Omega_{12},\\
&&\partial_\varphi\partial_\varphi=\sin^2\tz\ \Omega_{23}\Omega_{23}+\cos^2\tz\ \Omega_{13}\Omega_{13}+\sin\tz\cos\tz(\Omega_{23}\Omega_{13}+\Omega_{13}\Omega_{23}),
\een
which implies the conclusion of the lemma. $\hfill\square$

\end{appendices}

Qian Zhang

School of Mathematical Sciences, Beijing Normal University

Laboratory of Mathematics and Complex Systems, Ministry of Education Beijing 100875, China

Email: qianzhang@bnu.edu.cn

\end{document}